\documentclass[a4paper]{article}
\usepackage{arxiv}
\usepackage{comment}
\usepackage[utf8]{inputenc} 
\usepackage[T1]{fontenc}    
\usepackage{hyperref}       
\usepackage{url}            
\usepackage{booktabs}       
\usepackage{amsmath}
\usepackage{amssymb}
\usepackage{amsthm}
\usepackage{amsfonts}       
\usepackage{nicefrac}       
\usepackage{microtype}      
\usepackage{xcolor}
\usepackage{graphicx}
\usepackage{stmaryrd}
\usepackage{bm}
\usepackage[normalem]{ulem}
\usepackage[tickmarkheight=0.1cm,textsize=tiny]{todonotes}
\usepackage{algorithmicx}
\usepackage{algorithm}
\usepackage{algpseudocode}
\usepackage{caption}
\usepackage{subcaption}
\usepackage{lineno}
\usepackage{orcidlink}
\usepackage[inline,final]{showlabels}

\usepackage{multirow}
\usepackage{comment}
\usepackage{enumitem}
\usepackage{mathrsfs}

\hypersetup{
    colorlinks=true, 
    linkcolor=blue, 
    urlcolor=red, 
    citecolor=magenta,
    linktoc=all 
}

\newcommand{\Uad}{\mathcal{U}_{\textrm{ad}}}

\newtheorem{remark}{Remark}[section]

\DeclareUnicodeCharacter{0301}{\hspace{-1ex}\'{ }} 

\newcommand{\mb}[1]{\boldsymbol{#1}}
\newcommand{\p}{\partial}
\newcommand{\contra}[1]{\widetilde{#1}}

\newcommand{\mN}{\mathbb{N}}
\newcommand{\ud}{\textrm{d}}

\newcommand{\dd}{\mathrm{d}}
\newcommand{\iq}{Q}
\makeatletter
\newcommand*{\mint}[1]{%
	\mint@l{#1}{}%
}
\newcommand*{\mint@l}[2]{%
	\@ifnextchar\limits{%
		\mint@l{#1}%
	}{%
		\@ifnextchar\nolimits{%
			\mint@l{#1}%
		}{%
			\@ifnextchar\displaylimits{%
				\mint@l{#1}%
			}{%
				\mint@s{#2}{#1}%
			}%
		}%
	}%
}
\newcommand*{\mint@s}[2]{%
	\@ifnextchar_{%
		\mint@sub{#1}{#2}%
	}{%
		\@ifnextchar^{%
			\mint@sup{#1}{#2}%
		}{%
			\mint@{#1}{#2}{}{}%
		}%
	}%
}
\def\mint@sub#1#2_#3{%
	\@ifnextchar^{%
		\mint@sub@sup{#1}{#2}{#3}%
	}{%
		\mint@{#1}{#2}{#3}{}%
	}%
}
\def\mint@sup#1#2^#3{%
	\@ifnextchar_{%
		\mint@sup@sub{#1}{#2}{#3}%
	}{%
		\mint@{#1}{#2}{}{#3}%
	}%
}
\def\mint@sub@sup#1#2#3^#4{%
	\mint@{#1}{#2}{#3}{#4}%
}
\def\mint@sup@sub#1#2#3_#4{%
	\mint@{#1}{#2}{#4}{#3}%
}
\newcommand*{\mint@}[4]{%
	\mathop{}%
	\mkern-\thinmuskip
	\mathchoice{%
		\mint@@{#1}{#2}{#3}{#4}%
		\displaystyle\textstyle\scriptstyle
	}{%
		\mint@@{#1}{#2}{#3}{#4}%
		\textstyle\scriptstyle\scriptstyle
	}{%
		\mint@@{#1}{#2}{#3}{#4}%
		\scriptstyle\scriptscriptstyle\scriptscriptstyle
	}{%
		\mint@@{#1}{#2}{#3}{#4}%
		\scriptscriptstyle\scriptscriptstyle\scriptscriptstyle
	}%
	\mkern-\thinmuskip
	\int#1%
	\ifx\\#3\\\else_{#3}\fi
	\ifx\\#4\\\else^{#4}\fi
}
\newcommand*{\mint@@}[7]{%
	\begingroup
	\sbox0{$#5\int\m@th$}%
	\sbox2{$#5\int_{}\m@th$}%
	\dimen2=\wd0 %
	\let\mint@limits=#1\relax
	\ifx\mint@limits\relax
	\sbox4{$#5\int_{\kern1sp}^{\kern1sp}\m@th$}%
	\ifdim\wd4>\wd2 %
	\let\mint@limits=\nolimits
	\else
	\let\mint@limits=\limits
	\fi
	\fi
	\ifx\mint@limits\displaylimits
	\ifx#5\displaystyle
	\let\mint@limits=\limits
	\fi
	\fi
	\ifx\mint@limits\limits
	\sbox0{$#7#3\m@th$}%
	\sbox2{$#7#4\m@th$}%
	\ifdim\wd0>\dimen2 %
	\dimen2=\wd0 %
	\fi
	\ifdim\wd2>\dimen2 %
	\dimen2=\wd2 %
	\fi
	\fi
	\rlap{%
		$#5%
		\vcenter{%
			\hbox to\dimen2{%
				\hss
				$#6{#2}\m@th$%
				\hss
			}%
		}%
		$%
	}%
	\endgroup
}

\usepackage{esint}
\newcommand{\qint}{\mint{-}}

\newcommand{\intm}{\mb{\mathcal{V}}}
\newcommand{\projm}{\mb{\mathcal{P}}}

\definecolor{darkgreen}{RGB}{50, 200, 50}
\newcommand{\rev}[1]{{\color{black}#1\color{black}}} 
\newcommand{\reva}[1]{{\color{black}#1\color{black}}} 
\newcommand{\revb}[1]{{\color{black}#1\color{black}}} 

\title{Admissible Lax-Wendroff Flux Reconstruction Method with Automatic Differentiation on Adaptive Curved Meshes for Relativistic Hydrodynamics}

\author{
Sujoy~Basak \orcidlink{0009-0009-0612-6361}\thanks{Corresponding author}\\
Department of Mathematics\\
Indian Institute of Technology Delhi\\
New Delhi -- 110016, India\\
\texttt{sujoybasak42@gmail.com} \\
\texttt{maz228082@maths.iitd.ac.in} \\
\And
Arpit~Babbar \orcidlink{0000-0002-9453-370X} \\
Institute of Mathematics\\
Johannes Gutenberg University Mainz\\
Staudingerweg 9, 55122 Mainz, Germany\\
\texttt{ababbar@uni-mainz.de} \\
\And
Harish~Kumar \orcidlink{0000-0003-4746-2336}\\
Dept. of Mathematics, IIT Delhi, India--110016\\
\texttt{hkumar@iitd.ac.in}\\
\&\\ IITD-Abu Dhabi Campus, Abu Dhabi, UAE\\
\texttt{hkumar@iitdabudhabi.ac.ae} \\
\And
Praveen~Chandrashekar \orcidlink{0000-0003-1903-4107}\\
Centre for Applicable Mathematics\\
Tata Institute of Fundamental Research\\
Bangalore -- 560065, India\\
\texttt{praveen@math.tifrbng.res.in}
}
\begin{document}
\maketitle
\begin{abstract}
The relativistic hydrodynamics (RHD) equations can give rise to solutions which have shocks, contact discontinuities, and other sharp structures, which interact and evolve over time. Capturing these sharp waves effectively requires a mesh with high resolution, making the scheme computationally expensive. In this work, adaptive mesh refinement is used with the high-order Lax-Wendroff flux reconstruction (LWFR) method to solve the system of RHD equations, which is closed with general \revb{Synge-type} equations of state. To make the scheme Jacobian-free, the idea of automatic differentiation is incorporated for computing the temporal derivatives in the time average flux approximations. The high-order method is blended with an admissible low-order method at the subcell level to control the Gibbs oscillations and maintain the physical admissibility of the solution. Finally, several test cases involving high Lorentz factors, low densities, low pressures, strong shock waves, and other discontinuities are used to demonstrate the robustness, accuracy, and effectiveness of the proposed method. These simulations are performed with AMR using various linear and curved meshes to show the scheme's efficiency and ability to handle complex geometries.
\end{abstract}
\keywords{Relativistic hydrodynamics \and Admissibility preservation \and General \revb{Synge-type} equation of state \and Lax-Wendroff flux reconstruction \and Automatic differentiation \and Adaptive mesh refinement \and Curved mesh}
\section{Introduction}
The relativistic hydrodynamics (RHD) equations have important applications in astrophysics and high-energy physics~\cite{begelman1984theory,bottcher2012relativistic,mirabel1999sources,zensus1997parsec}. In the form of a conservation law, the system of RHD equations is given by~\cite{anile2005relativistic,landau1987see,synge1965relativity}
\begin{align}\label{eq: RHD_equation}
    \frac{\p \mb{u}}{\p t}+\sum_{i=1}^d\frac{\p \mb{f}_i(\mb{u})}{\p x_i}=\mb{0},
\end{align}
\revb{where} the vector of conservative variables is
\[
    \mb{u}=(D, m_1,m_2,\dots,m_d, E)^\top,
\]
and flux vector in $i^{\text{th}}$ direction is
\[
    \mb{f}_i(\mb{u})=(Dv_i, m_1v_i+p\delta_{1,i}, m_2v_i+p\delta_{2,i},\dots, m_dv_i+p\delta_{d,i}, m_i)^\top,\quad \forall i=1,\dots, d.
\]
The conservative variables $(D, m_1,m_2,\dots,m_d, E)$ can be expressed in terms of the primitive variables $(\rho, v_1,v_2,\dots,v_d, p)$ as
\[
    D=\rho\Gamma, \quad m_i=\rho h v_i\Gamma^2, \quad E={\rho h\Gamma^2}-p,
\]
where $\Gamma$ is the Lorentz factor and $h$ is the specific enthalpy. The Lorentz factor in terms of the velocity vector $\mb{v}=(v_1,v_2,\dots,v_d)^\top$ is given by
\[
    \Gamma=\frac{1}{\sqrt{1-|\mb{v}|^2}}.
\]
We have adopted a normalization of the unknowns, so that the speed of light is unity. The system of RHD equations~\eqref{eq: RHD_equation} can be made closed with an equation of state. The ideal equation of state (ID-EOS) is frequently used in the literature due to its simplicity; it is given by
\begin{align}\label{eq: ID_eos}
    h = 1+\frac{\gamma}{\gamma -1}\frac{p}{\rho},
\end{align}
where $\gamma$ is the specific heat ratio. However, because of its poor approximation in the relativistic cases~\cite{ryu2006equation}, we consider three more equations of state: TM-EOS~\cite{mathews1971hydromagnetic,mignone2005piecewise}, IP-EOS~\cite{sokolov2001simple}, and RC-EOS~\cite{ryu2006equation}, given by
\begin{align}
    h &= \frac{5p}{2\rho} + \sqrt{\frac{9p^2}{4\rho^2}+1},\label{eq: TM_eos}\\
    h &= \frac{2p}{\rho}+\sqrt{\frac{4p^2}{\rho^2} + 1},\label{eq: IP_eos}\\
    \text{and} \quad h &= \frac{2(6p^2 + 4p\rho +\rho^2)}{\rho (3p+ 2\rho)},\label{eq: RC_eos}
\end{align}
respectively. For more details on these equations of state, one can refer to Section~2.1 of~\cite{basak2025constraints}. \reva{Consequently, we need to use four different methods for conversion from conservative variables $(D, \mb{m}, E)$ to primitive variables $(\rho, \mb{v}, p)$, as it is dependent on the equation of state in use. For the ID-EOS~\eqref{eq: ID_eos}, we follow the method of conversion from Section~2.3 in~\cite{cai2024provably}. For the TM-EOS~\eqref{eq: TM_eos}, IP-EOS~\eqref{eq: IP_eos}, and RC-EOS~\eqref{eq: RC_eos}, we follow the methods from~\cite{ryu2006equation},~\cite{sokolov2001simple}, and~\cite{basak2025constraints}, respectively.} For the physical consistency of the RHD equations~\eqref{eq: RHD_equation} the solution vector $\revb{\mb{u}=(D, \mb{m}, E)^\top}$ should be in the admissible set
\begin{equation}\label{eq: ad_region_1}
    \Uad =\{\mb{u} \in \mathbb{R}^{d+2}: \rho (\mb{u}) >0, p(\mb{u})>0, \epsilon(\mb{u})>0, 1 -|\mb{v}(\mb{u})|>0\},
\end{equation}
where the specific internal energy $\epsilon = h - \frac{p}{\rho} - 1$. The set $\Uad$  can be reformulated as~\cite{wu2015high,wu2016physical},
\begin{equation}\label{eq: ad_region_2}
    \Uad' = \{\mb{u} \in \mathbb{R}^{d+2}: D(\mb{u}) = D >0,\ q(\mb{u}):= E-\sqrt{D^2 + |\mb{m}|^2}>0\}.
\end{equation}
The solution $\mb{u}$, which belongs to the admissible set, is called an admissible solution. The admissibility of the solution is also necessary to maintain the hyperbolicity of the system.

The development of efficient numerical schemes for solving RHD equations has a long history, starting from the finite difference scheme with an artificial viscosity method to capture the shock waves in~\cite{wilson1972numerical} to designing high-order schemes like essentially non-oscillatory (ENO) and weighted ENO schemes~\cite{dolezal1995relativistic,del2002efficient,tchekhovskoy2007wham,chen2022physical}, discontinuous Galerkin schemes~\cite{radice2011discontinuous}, entropy stable schemes~\cite{bhoriya2020entropy, biswas2022entropy, xu2024high}, and approximate piece-wise parabolic reconstruction schemes~\cite{marti1996extension,aloy1999genesis,mignone2005piecewise}. More works to numerically solve the RHD equations can be found in~\cite{marti1991numerical,marti1994analytical,BALSARA1994,dai1997iterative,ibanez1999riemann,wu2014finite,wu2014third,wu2021minimum}. 
However, most of these methods do not guarantee the physical admissibility of the solution. For finding the admissible solution, one can use a very small CFL number or highly diffusive schemes with a very fine mesh~\cite{zhang2006ram,hughes2002three}. But this way can be very expensive, making this an active area of research; several admissibility preserving schemes for the RHD equations are studied in~\cite{wu2015high,wu2016physical,qin2016bound,wu2017design, basak2025bound, basak2025constraints}. In this work, we develop a \revb{high-order scheme to solve the RHD equations} such that the element-averages (or cell-averages) of the solution are admissible, then using the scaling limiter from~\cite{zhang2010maximum} the solution is made admissible throughout the domain.

The solution of the RHD equations often has sharp and small-scale structures, along with shock and contact discontinuities, which evolve and interact with time in a complex manner. Capturing these wave structures needs a mesh with high resolution. On the other hand, in the smooth regions, a mesh with low resolution is sufficient. Since using a uniform mesh with high resolution throughout the domain substantially increases the computational cost, it is highly preferable to consider dynamic mesh adaptivity. It enables efficient resolution where and when it is required, providing a balance between accuracy and computational cost. In~\cite{brandt1977multi, berger1984adaptive}, the authors showed the importance of adaptive mesh refinement (AMR) in the field of science and engineering. In~\cite{berger1989local}, the authors have shown the benefit of AMR in the presence of shocks with Euler equations. One can also refer to~\cite{o2005adaptive, plewa2001amra, balsara2001divergence} for more literature on AMR and~\cite{hughes2002three,donmez2002general,anninos2005cosmos++, zhang2006ram, wang2008relativistic} for works with AMR in relativistic cases. We solve RHD equations incorporating the idea of AMR along with linear and curved meshes. The curved mesh plays an important role in adopting a computational domain having curved boundaries. It also helps in resolving the curved flow features in the solution more effectively.

In this work, we use the Lax-Wendroff flux reconstruction (LWFR) method~\cite{BABBAR2022111423, babbar2024admissibility}, which is a combination of the Lax-Wendroff (LW) method and flux reconstruction (FR) method, for solving the RHD equations. The LW method was initially proposed in~\cite{LW1} as a second-order finite difference method and later used in~\cite{qiu2005discontinuous,lou2020flux} to combine with finite element frameworks, giving an arbitrarily high-order accurate scheme. The LW method is known to be computationally efficient for the time update in finite element frameworks because of having a single stage. Whereas the Runge-Kutta methods of high order which is often used in the literature, have multiple stages, and each of the stages involves a change of data between different cores in modern parallel architectures. The multi-step methods~\cite{gottlieb2009high} of high-order, which can also be used for the time update, require information from several preceding time levels. The FR method was introduced in~\cite{huynh2007flux}, where the idea was to get a continuous reconstruction of the flux function using some correction functions. The quadrature-free FR method is well-suited for parallel implementation on modern vector processors~\cite{vincent2016towards,lopez2014verification,vandenhoeck2019implicit}. The LWFR method with AMR is recently implemented in~\cite{babbar2025lax} and used for the solution of the compressible Euler equations. To solve the RHD equations, the LWFR method with non-adaptive Cartesian grids is developed in~\cite{basak2025bound, basak2025constraints}, together with subcell limiting and guaranteed admissibility of the solution.

LWFR method needs the computation of \textit{time-average} flux functions~\cite{BABBAR2022111423}, which involves the temporal derivatives of the flux function. Usually, one can use the flux Jacobians and the high-order flux derivatives, which are high-dimensional tensors to change the temporal derivatives to spatial derivatives, as done in~\cite{qiu2005discontinuous, qiu2003finite,lou2020flux}. But since this approach needs high computational resources, in~\cite{burger2017approximate, zorio2017approximate}, the authors have proposed a finite difference method to approximate the temporal derivatives. The finite difference approach is shown to be faster~\cite{burger2017approximate} than the previous approach and used for successful computation of the time average fluxes for conservation laws~\cite{BABBAR2022111423, babbar2024admissibility}. However, in~\cite{basak2025constraints} the authors have shown that the finite difference approach fails for the case of conservation laws, where the flux computation of a quantity requires the quantity to be in the admissible set~\eqref{eq: ad_region_2}. In~\cite{basak2025constraints}, an additional scaling is used to deal with this difficulty. In this work, we use automatic differentiation (AD) to compute the temporal flux derivatives, which omits the need to apply the finite difference approximations and thereby eliminates the need for additional scaling in calculating the time average flux. In fact, this process is equivalent to the process of applying the chain rule as in~\cite{qiu2005discontinuous, qiu2003finite}, since it computes the directional derivatives of the flux functions. For more details on AD, one can refer to~\cite{griewank2008evaluating, babbar2025automatic}.

The rest of the paper is organized as follows. We discuss the numerical scheme on curved meshes in Section~\ref{sec: numerical_scheme}, along with the implementation of automatic differentiation. In Section~\ref{sec: boundary_treatments}, we discuss the boundary treatments used in this work depending on the flow direction. In Section~\ref{sec: mesh_adaptivity}, we discuss the mesh adaptivity, along with a discussion on the AMR indicator used in this work, and treatment of the non-conformal faces. Section~\ref{sec: numerical_results} presents the validation of the numerical scheme with various test cases, along with demonstrations of the benefit of using AMR. Finally, the paper ends with a summary in Section~\ref{sec: conclusion}.

\section{Numerical Scheme}\label{sec: numerical_scheme}
We explain the scheme in two dimensions ($d=2$), for which the system of conservation laws~\eqref{eq: RHD_equation} is given by
\begin{equation}\label{eq: RHD_equation_2d}
    \frac{\p \mb{u}}{\p t}+\frac{\p \mb{f}(\mb{u})}{\p x} + \frac{\p \mb{g}(\mb{u})}{\p y}=\mb{0},
\end{equation}
\begin{figure}
    \centering
    \includegraphics[width=0.9\textwidth]{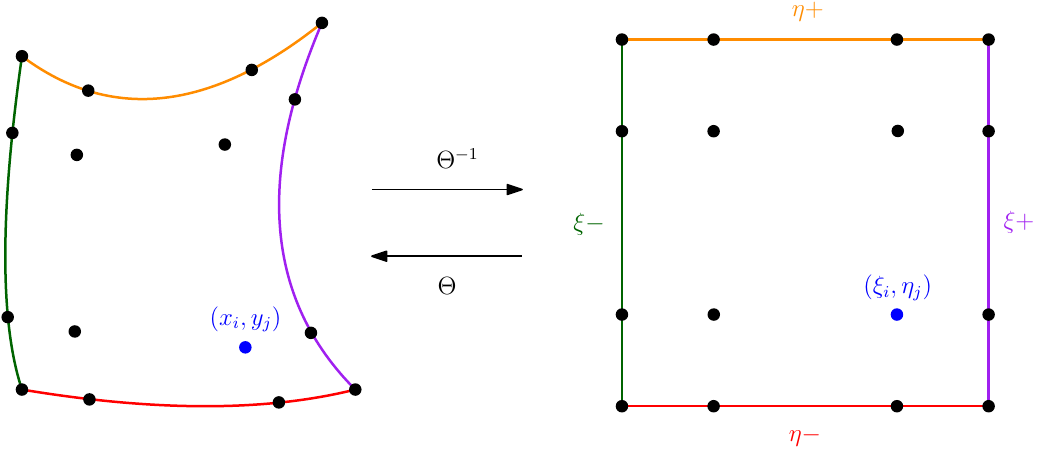}
    \caption{Illustration of the mapping between the physical and the reference element.}
    \label{fig: map_phys_to_ref}
\end{figure}
where we have set $(x_1, x_2) = (x, y)$, and the flux in $x$ and $y$- directions as $\mb{f_1} = \mb{f}$ and $\mb{f_2} = \mb{g}$ respectively in~\eqref{eq: RHD_equation} for notational simplicity. We begin by partitioning the computational domain $\Omega$ into $M$ disjoint quadrilaterals $\Omega_e$ such that $\Omega = \bigcup_{e=1}^M \Omega_e$.
Each element $\Omega_e$ is obtained by mapping the reference element $\contra{\Omega} = [-1,1]\times [-1, 1]$ (see Figure~\ref{fig: map_phys_to_ref}) by
\begin{equation}\label{eq: transform map}
    \Theta(\xi, \eta) = \sum_{(p,q)\in \mN^2_N} (x_p, y_q) \ell_p(\xi) \ell_q(\eta).
\end{equation}
\revb{The solution points in the reference element $\contra{\Omega}$ are denoted as} $(\xi_i, \eta_j)$, and $(x_i, y_j)$ are the solution points in the physical element $\Omega_e$. We have taken $(N+1)\times(N+1)$ solution points inside the reference element $\contra{\Omega}$ so that the solution is approximated by $N$-degree polynomials locally in each direction, which is same as the degree of the reference map.
The Lagrange's polynomials are denoted by
\[
    \ell_p(\xi) = \prod_{\substack{r=0\\ r\neq p}}^N \frac{\xi - \xi_r}{\xi_p - \xi_r},\qquad \ell_q(\eta) = \prod_{\substack{r=0\\ r\neq q}}^N \frac{\eta - \eta_r}{\eta_q - \eta_r},
\]
and the set of multi-indices is defined as
\[
    \mN^2_N = \big\{(p,q): p,q\in \{0,1,2,\dots,N\} \big\}.
\]
The solution points $-1\leq\xi_0\leq \xi_1 \leq \cdots \leq \xi_N\leq 1$ and $-1\leq\eta_0\leq \eta_1 \leq \cdots \eta_N\leq 1$ are taken as the Gauss-Legendre-Lobatto (GLL) nodes, \rev{whose quadrature rule performs exact integration for polynomials up to degree $(2N-1)$}.

Defining the Jacobian of the transformation map~\eqref{eq: transform map}
\begin{equation}\label{eq: jacobian}
    \mb{J} = \frac{\p (x,y)}{\p(\xi, \eta)},
\end{equation}
the system of conservation laws~\eqref{eq: RHD_equation_2d} is transformed into the reference element $\contra{\Omega}$
\begin{align}\label{eq: RHD_equation_2d_transformed}
    \frac{\p \contra{\mb{u}}}{\p t} + \frac{\p\contra{\mb{f}}}{\p \xi} + \frac{\p\contra{\mb{g}}}{\p \eta} = \mb{0},
\end{align}
where
\begin{equation}\label{eq: reference_sol_vs_physical_sol}
    \contra{\mb{u}} = |\mb{J}| \mb{u}, \quad 
    \contra{\mb{f}} = \frac{\p y}{\p \eta} \mb{f} - \frac{\p x}{\p \eta} \mb{g}, \quad 
    \contra{\mb{g}} = - \frac{\p y}{\p \xi} \mb{f} + \frac{\p x}{\p \xi} \mb{g}.
\end{equation}
A detailed mathematical derivation of this transformation is given in Appendix~\ref{sec: appendix}.

\subsection{Lax-Wendroff flux reconstruction (LWFR)}
The solution of RHD equations~\eqref{eq: RHD_equation_2d} inside the element $\Omega_e$  is expressed in terms of reference co-ordinates as
\begin{equation}\label{eq: sol_poly}
    \mb{u}_{e,h} (\xi, \eta, t) = \sum_{(i,j)\in\mN_N^2} \mb{u}_{e,i,j}(t) \ell_i(\xi) \ell_j(\eta),
\end{equation}
\revb{where} $\mb{u}_{e,i,j}$ are the unknowns that we need to find with the scheme, and are in fact the solutions at the solution points (also called node points), which are GLL nodes.

As already discussed in the introduction, the LWFR method is a combination of Lax-Wendroff and flux reconstruction method, and the key idea of the Lax-Wendroff method is to apply Taylor's expansion on the solution function $\mb{u}$ around time $t=t_n$ to find the solution at the next time level $t_{n+1}$
\[
    \mb{u}^{n+1} = \mb{u}^{n} + \sum_{m=1}^{N+1} \frac{\Delta t^m}{m!} \frac{\p^m \mb{u}^n}{\p t^m} + O(\Delta t^{N+2}).
\]
We have retained the terms up to $O(\Delta t^{N+1})$, since we desire the order of the scheme to be $N+1$ in space and time. Changing the temporal derivatives to spatial derivatives using~(\ref{eq: RHD_equation_2d_transformed},\ref{eq: reference_sol_vs_physical_sol}), the solution update is given by
\begin{equation}\label{eq: update_1}
    \mb{u}^{n+1} = \mb{u}^{n} - \frac{1}{|\mb{J}|}\sum_{m=1}^{N+1} \frac{\Delta t^m}{m!} \frac{\p^{m-1}}{\p t^{m-1}}\frac{\p\contra{\mb{f}}^n}{\p \xi} - \frac{1}{|\mb{J}|}\sum_{m=1}^{N+1} \frac{\Delta t^m}{m!} \frac{\p^{m-1}}{\p t^{m-1}}\frac{\p\contra{\mb{g}}^n}{\p \eta}.
\end{equation}
Defining the \textit{contravariant time average fluxes}
\begin{equation}\label{eq: ta_flux}
    \contra{\mb{F}} = \sum_{m=0}^{N} \frac{\Delta t^m}{(m+1)!}\frac{\p^m \contra{\mb{f}}^n}{\p t^m}, \qquad \contra{\mb{G}} = \sum_{m=0}^{N} \frac{\Delta t^m}{(m+1)!}\frac{\p^m \contra{\mb{g}}^n}{\p t^m},
\end{equation}
the coefficients $\mb{u}_{e,i,j}$ in~\eqref{eq: sol_poly} can be found using~\eqref{eq: update_1} as~\cite{babbar2025lax, basak2025bound}
\begin{equation}\label{eq: LWFR_update}
\revb{
    \mb{u}^{n+1}_{e,i,j} = \mb{u}^{n}_{e,i,j} - \frac{\Delta t}{|\mb{J}_{e,i,j}|} \left[ \frac{\p\contra{\mb{F}}_h}{\p \xi}(\xi_i, \eta_j) + \frac{\p\contra{\mb{G}}_h}{\p \eta}(\xi_i, \eta_j) \right], \qquad \forall (\xi_i, \eta_j) \in \contra{\Omega}_e,}
\end{equation}
where $\contra{\mb{F}}_h, \contra{\mb{G}}_h$ are the global continuous approximations of the contravariant time average fluxes $\contra{\mb{F}}, \contra{\mb{G}}$, found by incorporating the flux reconstruction idea~\cite{babbar2025lax, basak2025bound} as
\begin{align}\label{eq: FR}
\begin{split}
\contra{\mb{F}}_h(\xi,\eta_j) &= \Big[\contra{\mb{F}}^*_{\xi-}(\eta_j) - \contra{\mb{F}}_h^\delta(-1,\eta_j) \Big] c_L(\xi) + \contra{\mb{F}}_h^\delta(\xi,\eta_j) + \Big[\contra{\mb{F}}^*_{\xi+}(\eta_j) - \contra{\mb{F}}_h^\delta(1,\eta_j) \Big] c_R(\xi),\\
\contra{\mb{G}}_h(\xi_i,\eta) &= \Big[\contra{\mb{G}}^*_{\eta-}(\xi_i) - \contra{\mb{G}}_h^\delta(\xi_i,-1) \Big] c_L(\eta) + \contra{\mb{G}}_h^\delta(\xi_i,\eta) + \Big[\contra{\mb{G}}^*_{\eta+}(\xi_i) - \contra{\mb{G}}_h^\delta(\xi_i,1) \Big] c_R(\eta).
\end{split}
\end{align}
Here, $\contra{\mb{F}}_h^\delta, \contra{\mb{G}}_h^\delta$ are some approximations of the contravariant time average fluxes,
which are possibly discontinuous on the element boundaries, as will be discussed in Section~\ref{sec: TA_flux_ad}. \revb{The}
$\contra{\mb{F}}^*_{\xi\pm}$, $\contra{\mb{G}}^*_{\eta\pm}$, are the numerical fluxes at the faces $\xi\pm, \eta\pm$, respectively (see Figure~\ref{fig: map_phys_to_ref}), computed with Rusanov type approximation
\begin{align}\label{eq: numerical_flux}
\begin{split}
        \contra{\mb{F}}^*_{\xi\pm}(\cdot) &= \frac{1}{2}\big[\contra{\mb{F}}^-_{\xi\pm}(\cdot) + \contra{\mb{F}}^+_{\xi\pm}(\cdot) \big] - \frac{1}{2}\lambda_{\xi\pm}\big[\mb{U}^+_{\xi\pm}(\cdot) - \mb{U}^-_{\xi\pm}(\cdot) \big],\\
        \contra{\mb{G}}^*_{\eta\pm}(\cdot) &= \frac{1}{2}\big[\contra{\mb{G}}^-_{\eta\pm}(\cdot) + \contra{\mb{G}}^+_{\eta\pm}(\cdot) \big] - \frac{1}{2}\lambda_{\eta\pm}\big[\mb{U}^+_{\eta\pm}(\cdot) - \mb{U}^-_{\eta\pm}(\cdot) \big],
\end{split}
\end{align}
where $\contra{\mb{F}}^\pm_{\xi\pm}, \contra{\mb{G}}^\pm_{\eta\pm}$ are the trace values of the fluxes $\contra{\mb{F}}, \contra{\mb{G}}$ at the faces $\xi\pm, \eta\pm$ depending on the direction~\cite{basak2025bound}. The quantities $\mb{U}^\pm_{\xi\pm}, \mb{U}^\pm_{\eta\pm}$ are the traces for the time average solution
\[
    \mb{U} = \sum_{k=0}^N \frac{\Delta t ^k}{(k+1)!}\frac{\p^k \mb{u}}{\p t^k}
\]
and $\lambda_{\xi\pm},\lambda_{\eta\pm}$ are the local wave speed estimates at the faces $\xi\pm$, $\eta\pm$, respectively. For the RHD equations, the wave speed estimates are given by~\cite{basak2025bound},
\[
    \lambda = \max \left\{ \lambda^-_{\text{max}}, \lambda^+_{\text{max}} \right\},
\]
where 
\[
    \lambda^\pm_{\text{max}} = \max \left\{\frac{\left(1-{(s^\pm)}^2\right)v^\pm_{\perp} + s^\pm \sqrt{(1 - {|\mb{v}^\pm|}^2)Q^\pm}}{1 - (s^\pm)^2 |\mb{v}^\pm|^2}, v^\pm_{\perp}, \frac{\left(1-{(s^\pm)}^2\right)v^\pm_{\perp} \revb{-} s^\pm \sqrt{(1 - {|\mb{v}^\pm|}^2)Q^\pm}}{1 - (s^\pm)^2 |\mb{v}^\pm|^2}  \right\},
\]
with $Q^\pm = 1 - (v_{\perp}^\pm)^2 - (s^\pm)^2\left[|\mb{v}^\pm|^2 - (v_{\perp}^\pm)^2\right]$, and $v^\pm_{\perp}$ as the velocities of the fluid perpendicular to the boundary. \revb{The quantities} $s^\pm$ are the sound speeds~\cite{basak2025constraints} at the boundary. Here, $(\cdot)^-$ and $(\cdot)^+$ denote the trace values from left/lower and right/upper elements, respectively. For more details, one can refer to Section~3.2 in~\cite{basak2025bound}. The correction functions $c_L, c_R$ are taken as $g_2$ correction functions from~\cite{huynh2007flux}, which are degree $N+1$ polynomials having the properties
\[
    c_L(-1) = c_R(1) =1, \quad c_L(1) = c_R(-1)=0.
\]
It is worth noting that the FR scheme with GLL nodes as solution points is equivalent to a DG scheme with the solution points and the quadrature points taken as GLL nodes~\cite{babbar2025lax}. A detailed discussion on the equivalence of FR and DG schemes, depending on the choice of correction function, can be found in~\cite{huynh2007flux}.

\subsubsection{Computation of the contravariant time average flux with AD}\label{sec: TA_flux_ad}
Computation of contravariant time average fluxes $\contra{\mb{F}}(\xi_i, \eta_j), \contra{\mb{G}}(\xi_i, \eta_j)$~\eqref{eq: ta_flux} at solution points is one of the crucial steps in the LWFR method and can be done by shifting temporal derivative from the flux $\contra{\mb{f}}$ to the solution $\mb{u}$ using the Faà di Bruno’s formula~\cite{di1857note}
\begin{equation}\label{eq: faa_di_bruno}
    \frac{\p^m \contra{\mb{f}}}{\p t ^m} = \sum_{\pi \in \Pi} \left(\frac{\p^{|\pi|} \contra{\mb{f}}}{\p \mb{u} ^{|\pi|}} \prod_{B \in \pi} \frac{\p^{|B|} \mb{u}}{\p t^{|B|}}\right),
\end{equation}
where $\Pi$ denotes the set of all the partitions of the set $\{1, 2, \dots, m\}$. Here, $|\pi|$ denotes the number of blocks in the partition $\pi$ and $|B|$ denotes the number of elements in the block $B$, respectively. Using the formula~\eqref{eq: faa_di_bruno} to compute the temporal derivatives as done in~\cite{qiu2005discontinuous, qiu2003finite,lou2020flux}, needs the computation of flux Jacobians and higher derivatives of fluxes, which are high-dimensional tensors, increasing the computational cost. One remedy is to use the finite difference approach as used in~\cite{burger2017approximate, zorio2017approximate, BABBAR2022111423}, but it fails for the RHD equations~\cite{basak2025constraints} as flux evaluation for some non-admissible quantity has to be carried out. We use the automatic differentiation (AD) approach to compute the temporal derivatives in~\eqref{eq: ta_flux}. The AD approach computes the directional derivatives while avoiding the explicit construction of the flux Jacobians, and thus does not suffer from the performance issue of using~\eqref{eq: faa_di_bruno}~\cite{babbar2025automatic}. 

With the notations
\begin{equation}\label{eq: diff_not}
\contra{\mb{f}}^{(m)} = \Delta t^m \frac{\p^m \contra{\mb{f}}}{\p t^m}, \quad \contra{\mb{g}}^{(m)} = \Delta t^m \frac{\p^m \contra{\mb{g}}}{\p t^m}, \quad \mb{u}^{(m)} = \Delta t^m \frac{\p^m \mb{u}}{\p t^m}, \qquad m=1,2,\ldots,
\end{equation}
the equations~\eqref{eq: ta_flux} can be rewritten as
\[
    \contra{\mb{F}} = \sum_{m=0}^N \frac{\contra{\mb{f}}^{(m)}}{(m+1)!}, \quad \contra{\mb{G}} = \sum_{m=0}^N \frac{\contra{\mb{g}}^{(m)}}{(m+1)!}.
\]
For first-order derivative ($m=1$), $\contra{\mb{f}}^{(1)}$  can be expressed as
\begin{equation}\label{eq:f_der_1}
    \contra{\mb{f}}^{(1)} = \contra{\mb{f}}'(\mb{u})\mb{u}^{(1)},
\end{equation}
which is the directional derivative of $\contra{\mb{f}}$ with respect to $\mb{u}$ along the direction of $\mb{u}^{(1)}$. We find this directional derivative with the help of an AD library, as AD can be used to find the directional derivative of a function along any direction (for more details, one can refer to~\cite{babbar2025automatic}). Some of such libraries are \texttt{ForwardDiff.jl}~\cite{revels2016forward}, which exponentially increases the computational cost with the increase in the order of the derivatives~\cite{tan2023higher} and \texttt{TaylorDiff.jl}~\cite{tan2022taylordiff}, which uses the Taylor arithmetic (Chapter~13 in~\cite{griewank2008evaluating}). In \texttt{TaylorDiff.jl}, the rules for high-order derivatives are directly implemented~\cite{tan2023higher}. In this work, we use the library \texttt{Enzyme.jl}~\cite{moses2020instead}, which is fastest for this kind of application as reported in~\cite{babbar2025automatic}. The \texttt{Enzyme.jl} library works at the LLVM compiler level, which is used by Julia, and it does not need any additional storage for holding the derivatives. For more details, one can refer to Appendix~A in~\cite{babbar2025automatic}. For the high-order derivative term $\contra{\mb{f}}^{(m)}$, we introduce two derivative bundles
\begin{equation}\label{eq:der_bundles}
    \mathcal{B}_{\mb{u}}^{(m)} = \big(\mb{u}, \mb{u}^{(1)}, \dots , \mb{u}^{(m)}\big), \qquad 
     \tilde{\mathcal{B}}_{\mb{u}}^{(m)} = \big(\mb{u}^{(1)}, \mb{u}^{(2)} \dots , \mb{u}^{(m)}\big).
\end{equation}
From~\eqref{eq: faa_di_bruno}, it can be seen that $\contra{\mb{f}}^{(m)}$~\eqref{eq: diff_not} can be expressed as a function of $\mathcal{B}_{\mb{u}}^{(m)}$, and hence
\begin{equation}\label{eq: flux_der_high_order}
    \contra{\mb{f}}^{(m)} := \contra{\mb{f}}^{(m)}\left(\mathcal{B}_{\mb{u}}^{(m)}\right) 
    = \Delta t \frac{\p}{\p t} \contra{\mb{f}}^{(m-1)}\left(\mathcal{B}_{\mb{u}}^{(m-1)}\right) 
    = \Delta t \frac{\p \contra{\mb{f}}^{(m-1)}}{\p \mathcal{B}_{\mb{u}}^{(m-1)}} \frac{\p \mathcal{B}_{\mb{u}}^{(m-1)}}{\p t} 
    = \frac{\p \contra{\mb{f}}^{(m-1)}}{\p \mathcal{B}_{\mb{u}}^{(m-1)}} \tilde{\mathcal{B}}_{\mb{u}}^{(m)},
\end{equation}
as $\Delta t \frac{\p \mathcal{B}_{\mb{u}}^{(m-1)}}{\p t} = \tilde{\mathcal{B}}_{\mb{u}}^{(m)}$. So, $\contra{\mb{f}}^{(m)}$ is the directional derivative of $\contra{\mb{f}}^{(m-1)}$ with respect to $\mathcal{B}_{\mb{u}}^{(m-1)}$ in the direction $\tilde{\mathcal{B}}_{\mb{u}}^{(m)}$, and hence can be computed with AD in a recursive manner. It is worth mentioning that, since the equation~\eqref{eq: flux_der_high_order} computes $\contra{\mb{f}}^{(m)}$ exactly, this approach is equivalent to applying the Faà di Bruno’s formula~\eqref{eq: faa_di_bruno}. \rev{The same process is followed for finding the flux derivative term $\contra{\mb{g}}^{(m)}$}.

\revb{
\begin{remark}
    In practice, the quantity $\mb{u}^{(1)}$ in~(\ref{eq:f_der_1},\ref{eq:der_bundles}) is found using the conservation law~\eqref{eq: RHD_equation_2d_transformed} as 
    \begin{equation*}
        \mb{u}^{(1)} = -\frac{\Delta t}{|\mb{J}|}\left( \contra{\mb{f}}_{\xi} + \contra{\mb{g}}_{\eta}\right).
    \end{equation*}
   For $m\geq 2$, the $m$-order quantity in the derivative bundles~\eqref{eq:der_bundles} is found recursively using the flux derivative terms $\contra{\mb{f}}_{\xi}^{(m-1)}, \contra{\mb{g}}_{\xi}^{(m-1)}$ as 
    \begin{equation*}
        \mb{u}^{(m)} = -\frac{\Delta t}{|\mb{J}|} \left( \contra{\mb{f}}_{\xi}^{(m-1)} + \contra{\mb{g}}_{\eta}^{(m-1)}\right).
    \end{equation*}
\end{remark}
}

After finding the approximations of the contravariant time average fluxes $\contra{\mb{F}}(\xi_i, \eta_j), \contra{\mb{G}}(\xi_i, \eta_j)$ with this method, we compute the required quantities $\contra{\mb{F}}_h^\delta(\xi, \eta_j), \contra{\mb{G}}_h^\delta(\xi_i, \eta)$ in~\eqref{eq: FR} as
\[
    \contra{\mb{F}}_h^\delta(\xi, \eta_j) = \sum_{i=0}^N \contra{\mb{F}}(\xi_i, \eta_j) \ell_i(\xi), \quad \contra{\mb{G}}_h^\delta(\xi_i, \eta) = \sum_{j=0}^N \contra{\mb{G}}(\xi_i, \eta_j) \ell_j(\eta).
\]

\subsubsection{Blending limiter}\label{sec: blending}
According to Godunov's order barrier theorem, any linear high-order method to solve conservation laws creates Gibbs oscillations near discontinuities~\cite{godunov1959finite}. Thereby, the concept of using a non-linear limiting strategy comes into the picture. We employ the blending limiter~\cite{BABBAR2022111423}, where a high-order scheme is blended at the subcell level with a low-order scheme. \reva{ We now give a brief review of subcell based blending limiters. In~\cite{dumbser2014posteriori,dumbser2016simple} the authors have applied a subcell based approach in an \'a posteriori fashion by dividing an element into $2N+1$ subcells for $N$-degree polynomial approximation of the solution per space dimension. This idea is also implemented in relativistic regime in~\cite{zanotti2015solving} for solving the relativistic magnetohydrodynamics equations.
The subcell based limiter used in this work makes use of an a priori blending scheme. It was first used for LWFR schemes in~\cite{babbar2024admissibility}, and is originally inspired from the subcell limiter of~\cite{hennemann2021provably} for the discontinuous Galerkin method. One can find work on this priori approach for the subcell based blending limiter to solve the RHD equations using LWFR scheme in~\cite{basak2025bound, basak2025constraints}. In this work,} the low-order method is taken as a first-order finite volume method with the Rusanov flux~\cite{rusanov1962calculation}, whose proof for the admissibility preservation with ID-EOS~\eqref{eq: ID_eos} can be found in~\cite{basak2025bound, wu2017design}, and with the general EOSs~(\ref{eq: TM_eos}-\ref{eq: RC_eos}) can be found in~\cite{basak2025constraints}. This ensures the admissibility of the \textit{element-average} of the blended update
\begin{equation}\label{eq: blended_scheme}
    \mb{u}^{n+1}_{e,h} = (1-\alpha_e) \mb{u}^{H,n+1}_{e,h} + \alpha_e \mb{u}^{L,n+1}_{e,h},
\end{equation}
with $\mb{u}^{H,n+1}_{e,h}$ as the high-order update~\eqref{eq: LWFR_update} with a blended inter-element numerical flux~\cite{babbar2024admissibility, basak2025bound} in element $\Omega_e$. \revb{The} $\mb{u}^{L,n+1}_{e,h}$ is the update with the low-order scheme, which is performed after dividing each of the elements $\Omega_e$ into $(N+1) \times (N+1)$ sub-elements (also called as subcells). In~\eqref{eq: blended_scheme}, $\alpha_e$ is the indicator parameter, found from a smoothness indicator model as introduced in Section~4 in~\cite{hennemann2021provably}, using $\rho p \Gamma$ as the indicator quantity. Finally, using the Zhang-Shu scaling limiter~\cite{zhang2010maximum}, and using the admissibility of the element average, the temporal update is ensured to be admissible. The whole process of blending the scheme is explained thoroughly in~\cite{babbar2024admissibility} for the one-dimensional case, and in~\cite{basak2025bound} for the two-dimensional case, hence we do not repeat it here.

\begin{remark}
    The element-average of a quantity $\mb{u}_e$ in $\Omega_e$ can be found as
    \begin{equation}\label{eq: element-average_sol}
        \Bar{\mb{u}}_e = \frac{1}{|\Omega_e|} \int\limits_{\Omega_e} \mb{u}_e \dd x \dd y = \frac{1}{|\Omega_e|} \int\limits_{\contra{\Omega}} \contra{\mb{u}} \dd \xi \dd \eta = \frac{1}{|\Omega_e|} \sum_{(i,j)\in \mN_N^2} \contra{\mb{u}}_{i,j} w_i w_j =\frac{1}{|\Omega_e|} \sum_{(i,j)\in \mN_N^2} \mb{u}_{e,i,j} \mb{J}_{e,i,j} w_i w_j,
    \end{equation}
    where $\mb{J}_{e,i,j}$ is the Jacobian contribution of the transformation map~\eqref{eq: transform map} at $(x_i, y_j)$ in $\Omega_e$, and $w_i,w_j$ are the GLL quadrature weights in $[-1, 1]$.
\end{remark}

\section{Boundary treatments}\label{sec: boundary_treatments}
When simulating hyperbolic conservation laws, it is important to consider the nature of the flow near the boundaries of the computational domain, primarily when one considers artificial boundaries. Boundary conditions need to be prescribed for all the incoming components; however, imposing conditions on the outgoing components would result in overestimation. More literature in this direction can be found in~\cite{kreiss1970initial,strikwerda1976initial, higdon1986initial,kim2004integrated}. \revb{In this section}, we propose an idea to treat such boundaries in the LWFR framework.

In this work, the boundary conditions are handled in a weak sense by calculating the flux contributions at the boundary according to different boundary conditions. The implementation details of the boundary conditions are given below.

\subsection{Artificial boundary}\label{sec: artificial_boundary}
The boundary conditions are specified depending on the wave propagation across these boundaries. More precisely, the boundary conditions specify some solution variables at the boundary, depending on the number of incoming components, which are determined by the sign of the eigenvalues of the flux Jacobians~\cite{godlewski1991hyperbolic}. For the RHD equations, the expressions of these eigenvalues are given by
\begin{equation}\label{eq: eigen_value}
    \lambda_1 = \frac{(1-s^2)v_{\perp} + s \sqrt{(1 - |\mb{v}|^2)Q}}{1 - s^2 |\mb{v}|^2},\quad \lambda_2 = \lambda_3 = v_{\perp},\quad
    \lambda_4 = \frac{(1-s^2)v_{\perp} - s \sqrt{(1 - |\mb{v}|^2)Q}}{1 - s^2 |\mb{v}|^2},
\end{equation}
where $Q = 1 - v_{\perp}^2 - s^2(|\mb{v}|^2 - v_{\perp}^2)$, and $v_{\perp}$ is the velocity of the fluid perpendicular to the boundary from the interior of the domain. The sound speed $s$ is computed as in Section~2.1 of~\cite{basak2025constraints}, depending on the EOS in use. A detailed theoretical analysis regarding the boundary conditions depending on the eigenvalues can be found in~\cite{strikwerda1976initial}. One can also refer to~\cite{kreiss1970initial, higdon1986initial} for a characteristic analysis regarding the boundary conditions for some hyperbolic systems. 

\begin{figure}
    \centering
    \includegraphics[width=0.45\textwidth]{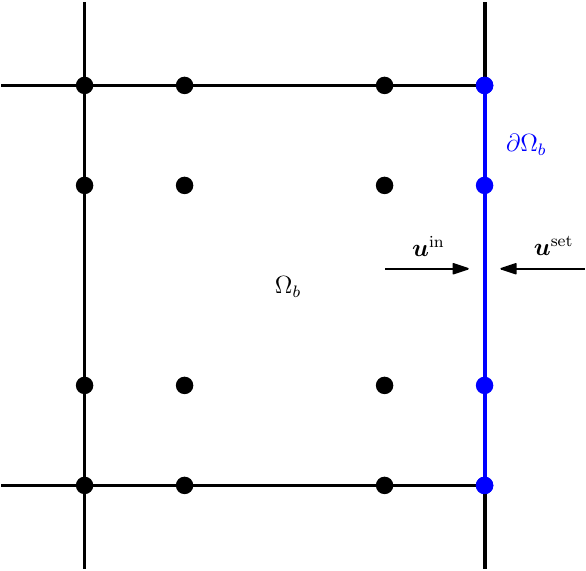}
    \caption{Element having the right face on the boundary.}
    \label{fig:boundary}
\end{figure}

We divide these types of artificial boundaries into three sub-categories based on the sign of the eigenvalues.
\begin{enumerate}[label=(\alph*)]
    \item Outflow: $\lambda_1, \lambda_2, \lambda_3, \lambda_4 > 0$.
    \item Inflow: $\lambda_1, \lambda_2, \lambda_3, \lambda_4 < 0$.
    \item Mixed flow: $\lambda_1, \lambda_2, \lambda_3, \lambda_4$ are not of same sign.
\end{enumerate}
The flux at the outer face of the elements neighboring a boundary, denoted by $\p\Omega_b$ (see Figure~\ref{fig:boundary}) is calculated depending on the type of the boundary as explained below. Here, $b$ is an element index for the elements having a face on a boundary.

\subsubsection{Outflow}\label{sec: outflow}
For outflow boundaries, all the eigenvalues~\eqref{eq: eigen_value} at the boundary are positive, that is, all the characteristics go outside the domain through the boundary. Hence, the flux at the outer face of boundary elements $\contra{\mb{F}}_{\p\Omega_b}$ is calculated purely from the internal solution
\[
    \contra{\mb{F}}_{\p\Omega_b} = \sum_{m=0}^N \frac{\contra{\mb{f}}^{(m)}(\mb{u}_{\Omega_b}^{\text{in}})}{(m+1)!},
\]
where, $\mb{u}_{\Omega_b}^{\text{in}}$ is the solution in the element $\Omega_b$ calculated from the scheme. \revb{The} $\contra{\mb{f}}^{(m)}(\cdot)$ is found with the AD procedure as explained in Section~\ref{sec: TA_flux_ad}. The computation of $\contra{\mb{G}}_{\p\Omega_b}$ follows the same procedure.

\subsubsection{Inflow}\label{sec: inflow}
For inflow boundaries, all the eigenvalues~\eqref{eq: eigen_value} at the boundary are negative, so all the characteristics come inside the domain across the boundary. Hence, to calculate the flux $\contra{\mb{F}}_{\p\Omega_b}$ at the outer face of the boundary elements, we need to set the solution at the boundary depending on the problem. Then the flux $\contra{\mb{F}}_{\p\Omega_b}$ at that face is calculated as
\[
    \contra{\mb{F}}_{\p\Omega_b} = \frac{1}{\Delta t}\int_{t_n}^{t_{n+1}}\contra{\mb{f}}\big(\mb{u}_{\p \Omega_b}^{\text{set}}(t) \big)\ud t,
\]
\revb{where} $\mb{u}_{\p \Omega_b}^{\text{set}}(t)$ is the solution that we set at the boundary $\p \Omega_b$. We evaluate the integral with the GLL quadrature rule. The computation of $\contra{\mb{G}}_{\p\Omega_b}$ follows the same procedure.

\subsubsection{Mixed flow}\label{sec: mixed_flow}
When the eigenvalues~\eqref{eq: eigen_value} at the boundary have mixed signs, then some family of characteristics enter the domain through the boundary, while the others go outside. We refer to these boundaries as mixed flow. In this case, it is necessary to specify some components at the boundary~\cite{kim2004integrated}.

Let us first consider the case, 
\[
    \lambda_1, \lambda_2, \lambda_3 \geq 0,\quad \lambda_4 < 0,
\]
then one family of characteristics enters the domain through this boundary. We specify the pressure $p$ near the boundary for such cases, as done for the Euler equations in~\cite{carlson2011inflow} and calculate the flux. Specifically, the flux $\contra{\mb{F}}_{\p\Omega_b}$ is calculated as
\begin{equation}\label{eq: mixed_bdry_flux}
    \contra{\mb{F}}_{\p\Omega_b} = \sum_{m=0}^N \frac{\contra{\mb{f}}^{(m)}(\mb{u}_{\Omega_b}^{*})}{(m+1)!},
\end{equation}
where the solution quantity $\mb{u}_{\Omega_b}^{*}$ is given by
\[
    \mb{u}_{\Omega_b}^{*} = \rev{\Psi}(\rho^{\text{in}}, \mb{v}_1^{\text{in}}, \mb{v}_2^{\text{in}}, p^{\text{set}}).
\]
The quantities with superscript $(\cdot)^{\text{in}}$ are calculated in the element $\Omega_b$ as a part of the LWFR scheme, and the quantity with superscript $(\cdot)^{\text{set}}$ is the value that we set, depending on the problem. \rev{The map $\Psi:\mathbb{R}^4 \to \mathbb{R}^4$ denotes the conversion from primitive to conservative variables.}

Let
\rev{\[
    \lambda_1> 0, \quad \lambda_2, \lambda_3, \lambda_4 < 0,
\]}
then one family of characteristics goes outside the domain, while the others come inside. For this type of boundary, following~\cite{carlson2011inflow} we specify the value of $\rho, v_1, v_2$ near the boundary. More specifically, we take
\[
    \mb{u}_{\Omega_b}^{*} = \rev{\Psi}(\rho^{\text{set}}, \mb{v}_1^{\text{set}}, \mb{v}_2^{\text{set}}, p^{\text{in}}),
\]
and calculate the flux $\contra{\mb{F}}_{\p\Omega_b}$ using~\eqref{eq: mixed_bdry_flux}. 

The same process is followed for the computation of the flux $\contra{\mb{G}}_{\p\Omega_b}$

In addition to the above artificial boundaries, we use two more boundary conditions in our work, namely periodic boundary and solid wall/reflective boundary. These boundaries are already discussed in Section~3.4 in~\cite{basak2025bound}, and hence we do not repeat them here.

\section{Mesh adaptivity} \label{sec: mesh_adaptivity}
\begin{figure}
    \centering
    \includegraphics[width=0.9\textwidth]{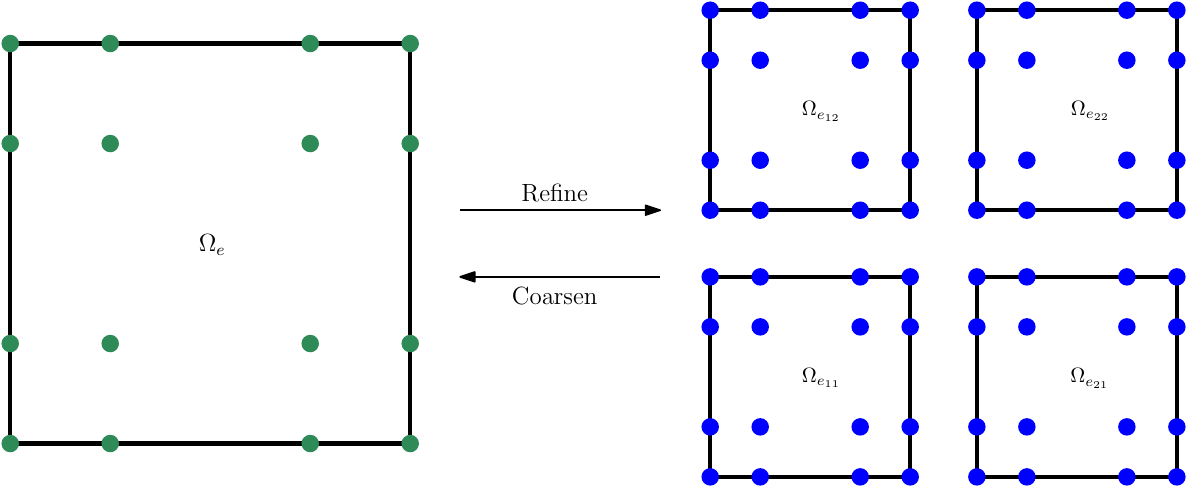}
    \caption{Refinement and coarsening of an element with the parent element on the left and child elements on the right.}
    \label{fig: refine_coarsen}
\end{figure}

When solving the RHD equations, to capture the small scale structures in the domain, the mesh needs to have high resolution. However, using a uniform mesh with high resolution throughout the domain is computationally expensive. As a remedy, we want to have high resolution only in the regions with sharp structures or other non-smooth features. Consequently, a solution smoothness indicator is needed to detect the regions where refining or coarsening is needed. We use a local indicator from~\cite{lohner1987adaptive} which is based on a central difference approximation for the second-order derivative of an indicator variable $Q= \rho p \Gamma$. For each element $\Omega_e$, we define the indicator as
\begin{equation}\label{eq: amr_indicator}
    \beta_e = \max_{(p,q) \in \mN_N^2} \max(\beta_p, \beta_q),
\end{equation}
where 
\begin{align*}
    \beta_p &= \frac{|\iq(\mb{u}_{p+, q}) - 2 \iq(\mb{u}_{pq}) + \iq(\mb{u}_{p-, q}) |}{|\iq(\mb{u}_{p+, q}) - \iq(\mb{u}_{pq})| + |\iq(\mb{u}_{p q}) - \iq(\mb{u}_{p-,q})| + \mathfrak{f}_{\text{w}} (|\iq(\mb{u}_{p+, q})| + 2 |\iq(\mb{u}_{pq})| + |\iq(\mb{u}_{p-, q}) |)},\\
    \beta_q &= \frac{|\iq(\mb{u}_{p, q+}) - 2 \iq(\mb{u}_{pq}) + \iq(\mb{u}_{p, q-}) |}{|\iq(\mb{u}_{p, q+}) - \iq(\mb{u}_{pq})| + |\iq(\mb{u}_{p q}) - \iq(\mb{u}_{p,q-})| + \mathfrak{f}_{\text{w}} (|\iq(\mb{u}_{p, q+})| + 2 |\iq(\mb{u}_{pq})| + |\iq(\mb{u}_{p, q-}) |)},
\end{align*}
\revb{and}
\begin{align*}
    p+ = \begin{cases}
            p, \quad \text{if } p = N\\
            p+1, \quad \text{otherwise}
        \end{cases},\qquad 
    p- = \begin{cases}
            p, \quad \text{if } p = 0\\
            p-1, \quad \text{otherwise}
        \end{cases},\\
    q+ = \begin{cases}
            q, \quad \text{if } q = N\\
            q+1, \quad \text{otherwise}
        \end{cases},\qquad 
    q- = \begin{cases}
            q, \quad \text{if } q = 0\\
            q-1, \quad \text{otherwise}
        \end{cases},
\end{align*}
and $\mathfrak{f}_{\text{w}}$ is taken as $0.2$ following~\cite{lohner1987adaptive}.
Then, the local refinement level of the element $\Omega_e$ is found using the three-level controller from~\texttt{Trixi.jl}~\cite{schlottkelakemper2025trixi}
\begin{equation}\label{eq: amr_controller}
    \text{level}_e = \begin{cases}
        \texttt{base\_level} \quad &\text{if}\ \beta_e < \epsilon_1\\
        \texttt{med\_level} \quad &\text{if}\ \epsilon_1 < \beta_e < \epsilon_2\\
        \texttt{max\_level} \quad &\text{if}\ \epsilon_2 < \beta_e.
    \end{cases}
\end{equation}
For an element $\Omega_e$, if $\text{level}_e$ is greater than the current refinement level of $\Omega_e$, it is marked for refinement and if $\text{level}_e$ is less than the current refinement level of $\Omega_e$, it is marked for coarsening. Initially, an initial refinement level is set throughout the whole mesh for a simulation. Then the elements, which are marked by the AMR controller~\eqref{eq: amr_controller} for refinement, get refined. For the coarsening, all four child elements (Figure~\ref{fig: refine_coarsen}) need to be marked for coarsening by the AMR controller~\eqref{eq: amr_controller} to get merged into a parent element. The mesh is refined further to ensure the neighboring elements differ only by one refinement level. The thresholds $\epsilon_1, \epsilon_2$ are taken according to the needs of the problem.

Let us denote the four child elements of the parent element $\Omega_e$ (Figure~\ref{fig: refine_coarsen}) during a refinement/coarsening as \revb{$\Omega_{e_{12}}, \Omega_{e_{22}}, \Omega_{e_{11}}, \Omega_{e_{21}}$}. 
We again take some notations as
\[
    I = [-1,1],\quad I_\revb{1} = [-1,0], \quad I_\revb{2} = [0,1],
\]
and the corresponding bijective maps $\phi_r:I_r\to I$, $r=\revb{1,2}$ as
\[
    \phi_\revb{1}(\xi) = 2\xi +1, \quad \phi_\revb{2}(\xi) = 2\xi -1.
\]

The four reference elements corresponding to the child elements \revb{$\Omega_{e_{k, s}},\ k,s \in \{1,2\}$} are denoted by 
\[
\contra{\Omega}_{\revb{12}} = I_\revb{1}\times I_\revb{2},\quad \contra{\Omega}_{\revb{22}} = I_\revb{2}\times I_\revb{2}, \quad \contra{\Omega}_{\revb{11}} = I_\revb{1}\times I_\revb{1}, \quad \contra{\Omega}_{\revb{21}} = I_\revb{2}\times I_\revb{1}.
\]
The solution inside the reference element $\contra{\Omega}$ can be expressed in the following form
\begin{equation}\label{eq: sol_in_parent_reference}
    \contra{\mb{u}}(\xi, \eta) = \sum_{(i,j)\in \mN_N^2} \contra{\mb{u}}_{i,j} \ell_i(\xi) \ell_j(\eta) = \sum_{(i,j)\in \mN_N^2} \mb{u}_{e,i,j} \mb{J}_{e,i,j} \ell_i(\xi) \ell_j(\eta),\qquad (\xi, \eta) \in \contra{\Omega},
\end{equation}
and inside $\contra{\Omega}_{\revb{k,s}}$, \revb{$k,s\in\{1,2\}$},
\begin{equation}\label{eq: sol_in_child_reference}
    \contra{\mb{u}}_{\revb{k,s}}(\xi, \eta) = \sum_{(i,j)\in \mN_N^2} \contra{\mb{u}}_{\revb{k,s},i,j} \ell_i\big(\phi_{\revb{k}}(\xi)\big) \ell_j\big(\phi_{\revb{s}}(\eta)\big) = \sum_{(i,j)\in \mN_N^2} \mb{u}_{e_\revb{k,s},i,j} \mb{J}_{e_\revb{k,s},i,j} \ell_i\big(\phi_{\revb{k}}(\xi)\big) \ell_j\big(\phi_{\revb{s}}(\eta)\big), \qquad (\xi, \eta) \in \contra{\Omega}_\revb{k,s}.
\end{equation}

\subsection{Refinement of element}\label{sec: refinement}
When the element $\Omega_e$ is refined into four child elements $\Omega_{e_\revb{k,s}}$ (Figure~\ref{fig: refine_coarsen}), the corresponding reference element $\contra{\Omega}$ is refined into four child reference elements $\contra{\Omega}_{\revb{k,s}}$, $\revb{k,s\in \{1,2\}}$. From~(\ref{eq: sol_in_parent_reference},~\ref{eq: sol_in_child_reference}), we can see that, to find the solution in the child elements, it is enough to express $\contra{\mb{u}}_{\revb{k,s},i,j}$ in terms of $\contra{\mb{u}}_{i,j}$, $\forall\ i,j \in \{1,2, \dots, N\}$. As the child elements will introduce new solution points in the parent element, one can use interpolation to find the solutions at all the newly introduced solution points. We use the Lagrange interpolation as
\begin{align*}
    \contra{\mb{u}}_{\revb{k,s},i,j} = \sum_{(p,q)\in \mN_N^2} \contra{\mb{u}}_{p,q} \ell_p\big(\phi_{\revb{k}}^{-1}(\xi_i)\big) \ell_q\big(\phi_{\revb{s}}^{-1}(\eta_j)\big)
        = \sum_{(p,q)\in \mN_N^2} \contra{\mb{u}}_{p,q} (\intm_{\revb{k}})_{i}^{p} (\intm_{\revb{s}})_{j}^{q},
\end{align*}
where the interpolation matrices $(\intm_{\revb{k}})_{i}^{p}, (\intm_{\revb{s}})_{j}^{q}$ are defined as
\begin{equation}\label{eq: 1d_interpol_matrices}
    (\intm_{\revb{k}})_{i}^{p} = \ell_p\big(\phi_{\revb{k}}^{-1}(\xi_i)\big), \quad     (\intm_{\revb{s}})_{j}^{q} = \ell_q\big(\phi_{\revb{s}}^{-1}(\eta_j)\big).
\end{equation}

The solution at $(x_i,y_j)$ in the element $\Omega_{e_s}$ is computed as
\[
    \mb{u}_{e_{\revb{k,s}},i,j} = \frac{1}{|\mb{J}_{e_{\revb{k,s}},i,j}|} \contra{\mb{u}}_{\revb{k,s},i,j}.
\]

Moreover, we have
\begin{align}
    \sum_{\revb{k,s}=1}^\revb{2} \qint\limits_{\Omega_{e_\revb{k,s}}} \mb{u}_{e_\revb{k,s}} (x,y) \dd x \dd y &= \sum_{\revb{k,s}=1}^\revb{2} \qint\limits_{\contra{\Omega}_\revb{k,s}} \contra{\mb{u}}_\revb{k,s} (\xi, \eta) \dd \xi \dd \eta \nonumber
    = \sum_{\revb{k,s}=1}^\revb{2} \qint\limits_{\contra{\Omega}_\revb{k,s}} \contra{\mb{u}} (\xi, \eta) \dd \xi \dd \eta \nonumber\\
    &= \qint\limits_{\contra{\Omega}} \contra{\mb{u}} (\xi, \eta) \dd \xi \dd \eta 
    = \qint\limits_{\Omega} \mb{u} (x, y) \dd x \dd y, \label{eq: int_eqn_refinement}
\end{align}
and hence, the conservation property of the scheme is maintained~\cite{babbar2025lax}. Here and in the rest of the paper, $\qint$ is used to denote the integral approximation computed with GLL quadrature.

\subsection{Coarsening of elements}\label{sec: coarsening}
Suppose, during the coarsening, the four child elements $\Omega_{e_\revb{k,s}}$, \revb{$k,s\in\{1,2\}$} merge to form the parent element $\Omega_e$ (Figure~\ref{fig: refine_coarsen}). Then the four reference elements $\contra{\Omega}_{\revb{k,s}}$, \revb{$k,s\in\{1,2\}$} are merged to form the parent reference element $\contra{\Omega}$. The solution $\contra{\mb{u}}$ in $\contra{\Omega}$ is expressed in terms of the solution $\contra{\mb{u}}_\revb{k,s}$ in $\contra{\Omega}_{\revb{k,s}}$ with the $L_2$ projection, that is, the following equation should be satisfied
\[
    \qint\limits_{\contra{\Omega}} \contra{\mb{u}}(\xi, \eta) \ell_i(\xi)\ell_j(\eta) \dd\xi \dd\eta = \sum_{\revb{k,s}=1}^\revb{2}\qint\limits_{\contra{\Omega}_{\revb{k,s}}} \contra{\mb{u}}_{\revb{k,s}}(\xi, \eta) \ell_i(\xi)\ell_j(\eta) \dd\xi \dd\eta, \quad \forall\ (i,j) \in \mN_N^2.
\]
Using the equations~(\ref{eq: sol_in_parent_reference},~\ref{eq: sol_in_child_reference}) we get
\begin{align}\label{eq: projection_int_eqn}
    \sum_{(p,q)\in\mN_N^2}\qint\limits_{\contra{\Omega}} \contra{\mb{u}}_{p,q} \ell_p(\xi) \ell_q(\eta) \ell_i(\xi)\ell_j(\eta) \dd\xi \dd\eta = \sum_{\revb{k,s}=1}^\revb{2} \sum_{(p,q)\in\mN_N^2}\qint\limits_{\contra{\Omega}_{\revb{k,s}}} \contra{\mb{u}}_{\revb{k,s},p,q} \ell_p\big( \phi_{\revb{k}}(\xi)\big) \ell_q\big( \phi_{\revb{s}}(\eta)\big) \ell_i(\xi)\ell_j(\eta)& \dd\xi \dd\eta,\\
    \quad &\forall (i,j) \in \mN_N^2.\nonumber
\end{align}
We again have
\begin{align*}
    \qint\limits_{\contra{\Omega}} \ell_p(\xi) \ell_q(\eta) \ell_i(\xi)\ell_j(\eta) \dd\xi \dd\eta &= \delta_{pi}\delta_{qj} w_p w_q,\\
    \qint\limits_{\contra{\Omega}_s} \ell_p\big( \phi_{\revb{k}}(\xi)\big)  \ell_q\big( \phi_{\revb{s}}(\eta)\big) \ell_i(\xi)\ell_j(\eta) \dd\xi \dd\eta &= \frac{1}{4} \big[\ell_i \big(\phi_{\revb{k}}^{-1} (\xi_p) \big) \ell_j \big(\phi_{\revb{s}}^{-1} (\eta_q) \big) w_p w_q \big],
\end{align*}
and using these two relations, the equation \eqref{eq: projection_int_eqn} becomes
\[
    \contra{\mb{u}}_{ij} = \frac{1}{4}\sum_{\revb{k,s}=1}^\revb{2} \sum_{(p,q)\in\mN_N^2} \contra{\mb{u}}_{\revb{k,s},p,q} \frac{w_p}{w_i}\ell_i \big(\phi_{\revb{k}}^{-1} (\xi_p) \big) \frac{w_q}{w_j}\ell_j \big(\phi_{\revb{s}}^{-1} (\eta_q) \big) = \sum_{\revb{k,s}=1}^\revb{2} \sum_{(p,q)\in\mN_N^2} \contra{\mb{u}}_{\revb{k,s},p,q} (\revb{\projm_{k}})_{p}^{i} (\revb{\projm_{\revb{s}}})_{q}^{j},
\]
where the projection matrices $(\revb{\projm_{\revb{k}}})_{p}^{i}, (\revb{\projm_{\revb{s}}})_{q}^{j}$ are defined as
\begin{equation}\label{eq: 1d_projection_matrices}
    (\revb{\projm_{k}})_{p}^{i} = \frac{1}{2} \frac{w_p}{w_i}\ell_i \big(\phi_{\revb{k}}^{-1} (\xi_p) \big),\qquad (\revb{\projm_{s}})_{q}^{j} = \frac{1}{2} \frac{w_q}{w_j}\ell_j \big(\phi_{\revb{s}}^{-1} (\eta_q) \big).
\end{equation}
Finally the solution at $(x_i, y_j)$ in the element $\Omega_e$ can be found by
\[
    \mb{u}_{e,i,j} = \frac{1}{|\mb{J}_{e,i,j}|} \contra{\mb{u}}_{i,j}.
\]
Again, summing the equation~\eqref{eq: projection_int_eqn} over $(i,j)\in \mN_N^2$ we get
\begin{align}
    \sum_{\revb{k,s}=1}^\revb{2} \qint\limits_{\contra{\Omega}_{\revb{k,s}}} \contra{\mb{u}}_{\revb{k,s}} (\xi,\eta) \dd \xi \dd \eta &= \qint\limits_{\contra{\Omega}} \contra{\mb{u}} (\xi, \eta) \dd \xi \dd \eta\nonumber\\
    \sum_{\revb{k,s}=1}^\revb{2} \qint\limits_{\Omega_{e_\revb{k,s}}} \mb{u}_{e_\revb{k,s}} (x,y) \dd x \dd y &= \qint\limits_{\Omega} \mb{u} (x, y) \dd x \dd y, \label{eq: int_eqn_coarsening}
\end{align}
which ensures the conservation property of the scheme.
\begin{remark}
    The transformation of the solution between the elements of different refinement levels is carried out with the reference solutions, which ensures that the equations~(\ref{eq: int_eqn_refinement},~\ref{eq: int_eqn_coarsening}) are satisfied. These equations are important for the conservation property of the scheme~\cite{babbar2025lax}.
\end{remark}

\subsection{Flux at the non-conformal face}
\begin{figure}
    \centering
    \includegraphics[width=0.6\textwidth]{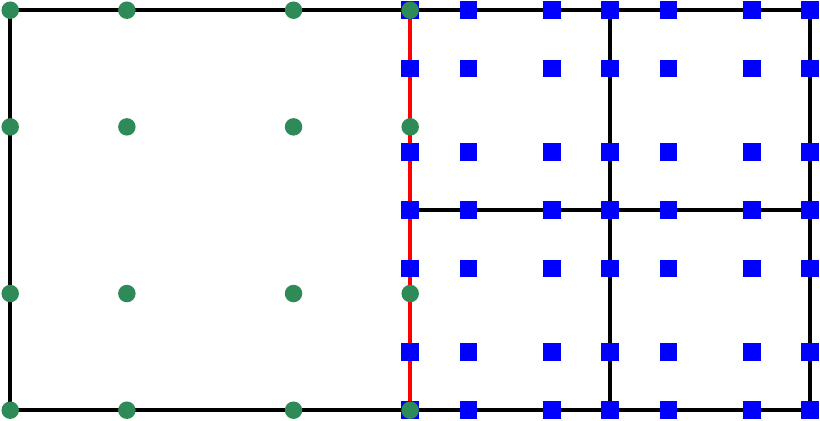}
    \caption{Neighboring elements with non-conformal face in red.}
    \label{fig: non-conformmal_face}
\end{figure}
During the adaptivity of the mesh, when the neighboring elements are of different refinement levels, the solution points on the faces do not match and are called non-conformal faces. In this work, the difference of refinement level of two neighboring elements is at most one, and hence it looks like Figure~\ref{fig: non-conformmal_face}. For the FR type schemes, the flux at these non-conformal faces needs to be computed at a refined face, and here we use the mortar element method~\cite{kopriva1996conservative, kopriva2002computation} for this purpose. 
\begin{figure}
    \centering
    \includegraphics[width=0.97\textwidth]{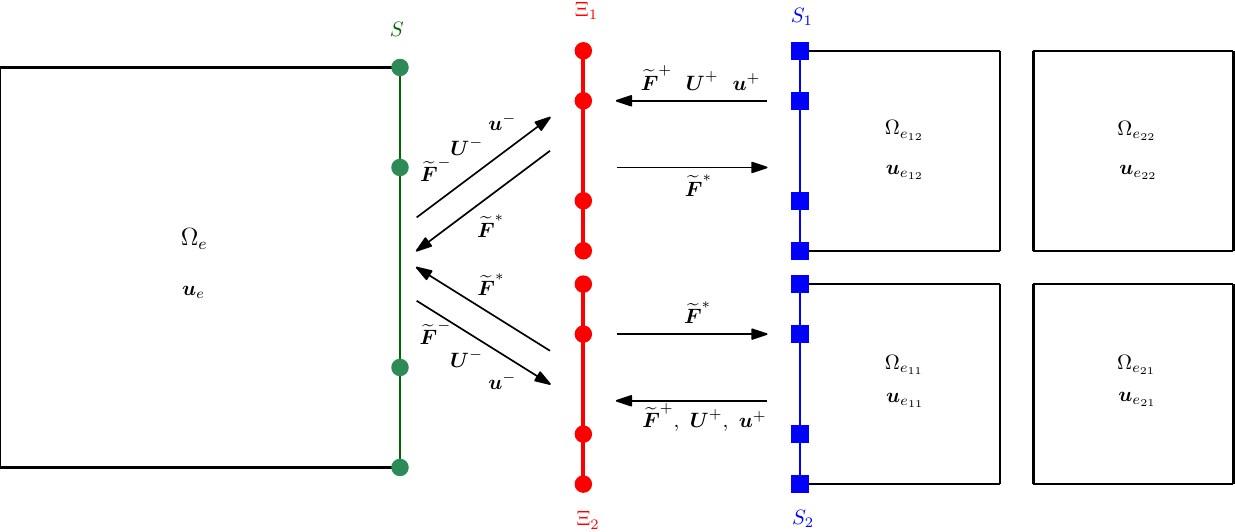}
    \caption{Visualization of numerical flux calculation using mortars.}
    \label{fig: mortars}
\end{figure}
From~\eqref{eq: numerical_flux}, the numerical flux at the face is a function of the trace values of the contravariant time average flux $\contra{\mb{F}}^\pm$ or $\contra{\mb{G}}^\pm$, time average solution $\mb{U}^\pm$, and solution polynomial $\mb{u}^\pm$. As a part of the mortar element method, we first prolong these values from both neighboring elements to common mortars $\Xi_1, \Xi_2$ (Figure~\ref{fig: mortars}), then compute the numerical flux $\contra{\mb{F}}^*$ there and map it back to the element faces $S$ and $S_1, S_2$.

Suppose $\contra{\mb{F}}^-_{\Xi_r} \mb{U}^-_{\Xi_r}, \mb{u}^-_{\Xi_r}$ and $\contra{\mb{F}}^+_{\Xi_r} \mb{U}^+_{\Xi_r}, \mb{u}^+_{\Xi_r}$ denote the prolonged values at the mortars $\Xi_r$ for $r=1,2$ from the elements having lower and higher refinement levels, respectively. Since the solution points at the mortars $\Xi_r$ and the face of higher refined elements $S_r$ are the same, the prolonged values $(\cdot)_{\Xi_r}^+$ are found with the identity map. 
For the prolonged values $(\cdot)_{\Xi_r}^-$, we use interpolation. More specifically, the prolonged value of the contravariant time average flux at the mortar can be expressed as
\begin{align*}
    \contra{\mb{F}}^-_{\Xi_r}(\eta) = \sum_{j=0}^N \ell_j \big(\phi_r(\eta)\big) \contra{\mb{F}}_{\Xi_r, j}^- , \qquad \eta \in \Xi_r
\end{align*}
where the coefficients $\contra{\mb{F}}_{\Xi_r, j}^-$ are found with interpolation as
\begin{align*}
    \contra{\mb{F}}_{\Xi_r, j}^- = \contra{\mb{F}}_S\big(\phi_r^{-1}(\eta_j) \big) = \sum_{q=0}^N \contra{\mb{F}}_S(\eta_q) (\intm_r)_j^q,\qquad \forall\ j=0,1,\dots, N.
\end{align*}
The interpolation matrix $(\intm_r)_j^q$ is defined in~\eqref{eq: 1d_interpol_matrices}.
Similarly, we can find $\mb{U}^-_{\Xi_r}, \mb{u}^-_{\Xi_r}$. Then the numerical flux at the mortars $\contra{\mb{F}}^*_{\Xi_r}$, $r=1,2$ are computed as in~\eqref{eq: numerical_flux} using the prolonged values. 

Since the solution points on the faces $S_r$ are the same as those on the mortars $\Xi_r$, we use the identity map to transfer the fluxes and denote them as $\contra{\mb{F}}^*_{S_r}$. To find the flux $\contra{\mb{F}}^*_{S}$ on the face of the element with lower refined level $S$, we use the $L_2$ projection, that is the following relation should be satisfied
\[\revb{
    \int\limits_S \contra{\mb{F}}^*_{S} \ell_j(\eta) \dd\eta = \sum_{r=1}^2 \int\limits_{\Xi_r} \contra{\mb{F}}^*_{\Xi_r} \ell_j(\eta) \dd\eta,\qquad \forall\ j = 0,1,\dots,N.}
\]
Now, doing the process similar to the Section~\ref{sec: coarsening}, the flux at the required face can be found as
\[
    \contra{\mb{F}}^*_{S,j} = \sum_{r=1}^2 \sum_{q=0}^N \contra{\mb{F}}^*_{\Xi_r,q} (\revb{\projm_r})_q^j,
\]
where $(\revb{\projm_r})_q^j$ is the projection matrix defined in~\eqref{eq: 1d_projection_matrices}.

The computation of the numerical fluxes in the other direction $\contra{\mb{G}}^*_{S_r}$, $\contra{\mb{G}}^*_{S}$ follows the same procedure, and hence we do not repeat it.

\begin{remark}
    As discussed in Section~\ref{sec: blending}, for the admissibility of the solution, we incorporate the flux limiting method from~\cite{babbar2024admissibility}, which is also used for the two-dimensional case in~\cite{basak2025bound}, and then use a scaling limiter~\cite{zhang2010maximum}. However, this strategy only ensures that the solutions at the solution points are admissible. Since the refinement of the mesh from $\Omega_e$ to $\Omega_{e_{k,s}}$ introduces additional node points, the scaling limiter is used with respect to the element-average $\Bar{\mb{u}}_e$~\eqref{eq: element-average_sol} to ensure the admissibility of the element-averages $\Bar{\mb{u}}_{e_{k,s}}$ in the refined elements, or the prolonged solutions with interpolation $\mb{u}^-_{\Xi_r}$ at the mortars.
\end{remark}

\section{Numerical results}\label{sec: numerical_results}
To demonstrate the robustness of the method, various numerical test cases with high Lorentz factors, strong shocks, rarefactions, and other sharp structures are simulated. Along with the linear meshes, we have also used curved meshes and meshes generated with Gmsh~\cite{geuzaine2009gmsh} to verify the capability of the scheme. For the computation of time step $\Delta t$, we have used the \textit{error-based time stepping} (Section~5.2 in~\cite{babbar2025lax}), unless mentioned otherwise. Note that the CFL-based time step computation used in~\cite{BABBAR2022111423, basak2025bound}, which involves the Fourier stability analysis on Cartesian meshes, does not guarantee the stability of the scheme on curved meshes~\cite{babbar2025lax}. The boundaries of the computational domains are treated as discussed in Section~\ref{sec: boundary_treatments}, considering the nature of the flow near the boundary. For all the simulations with AMR, the initial refinement level is taken equal to the $\texttt{max\_level}$~\eqref{eq: amr_controller}, set for that particular simulation. Along with the AMR results, we have also shown the results with a uniform mesh equivalent to the maximum refinement level in AMR, to show that the AMR can capture the same flow structures in the solution with less computational cost.

All the simulations are carried out using \texttt{RHDTrixiLW.jl}~\cite{RHDTrixiLW} on $32$ CPU threads using shared memory based parallelization. The codes are written in Julia using \texttt{TrixiLW.jl}~\cite{trixilw} and \texttt{Trixi.jl}~\cite{ranocha2021adaptive, schlottkelakemper2025trixi, schlottke2021purely} as libraries. The computations are performed on a system equipped with an Intel\textregistered\ Xeon\textregistered\ Gold 5220R CPU (2.20GHz), featuring 48 physical cores and running a Linux operating system.

\subsection{Sinusoidal smooth test}\label{sec: sinusoidal smooth test}
\begin{figure}[]
    \centering
    \begin{subfigure}{0.3\textwidth}
    \includegraphics[width=\linewidth]{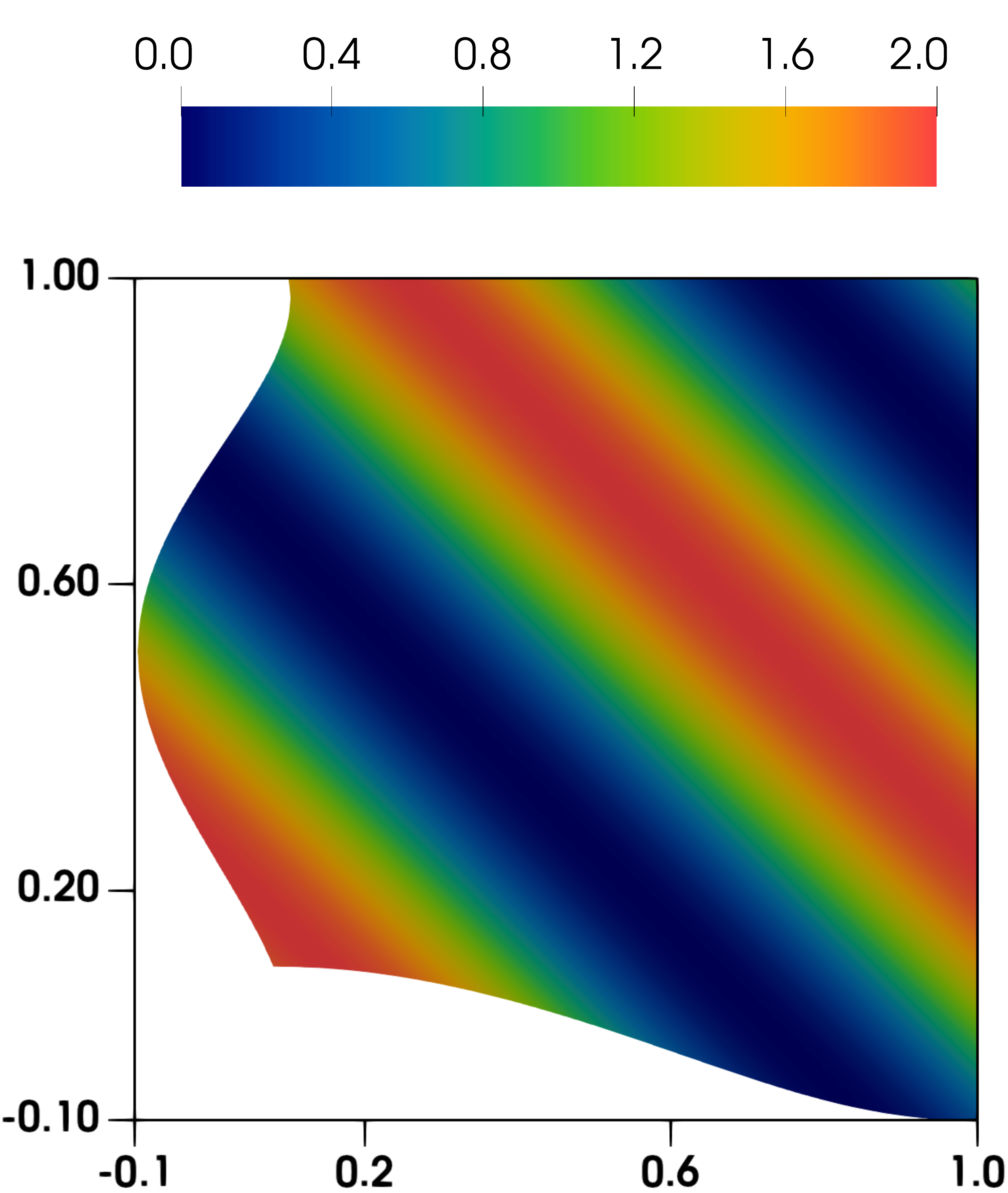}
    \caption{$\rho$ at $t=0.0$.}
    \end{subfigure}
    \quad
        \begin{subfigure}{0.3\textwidth}
    \includegraphics[width=\linewidth]{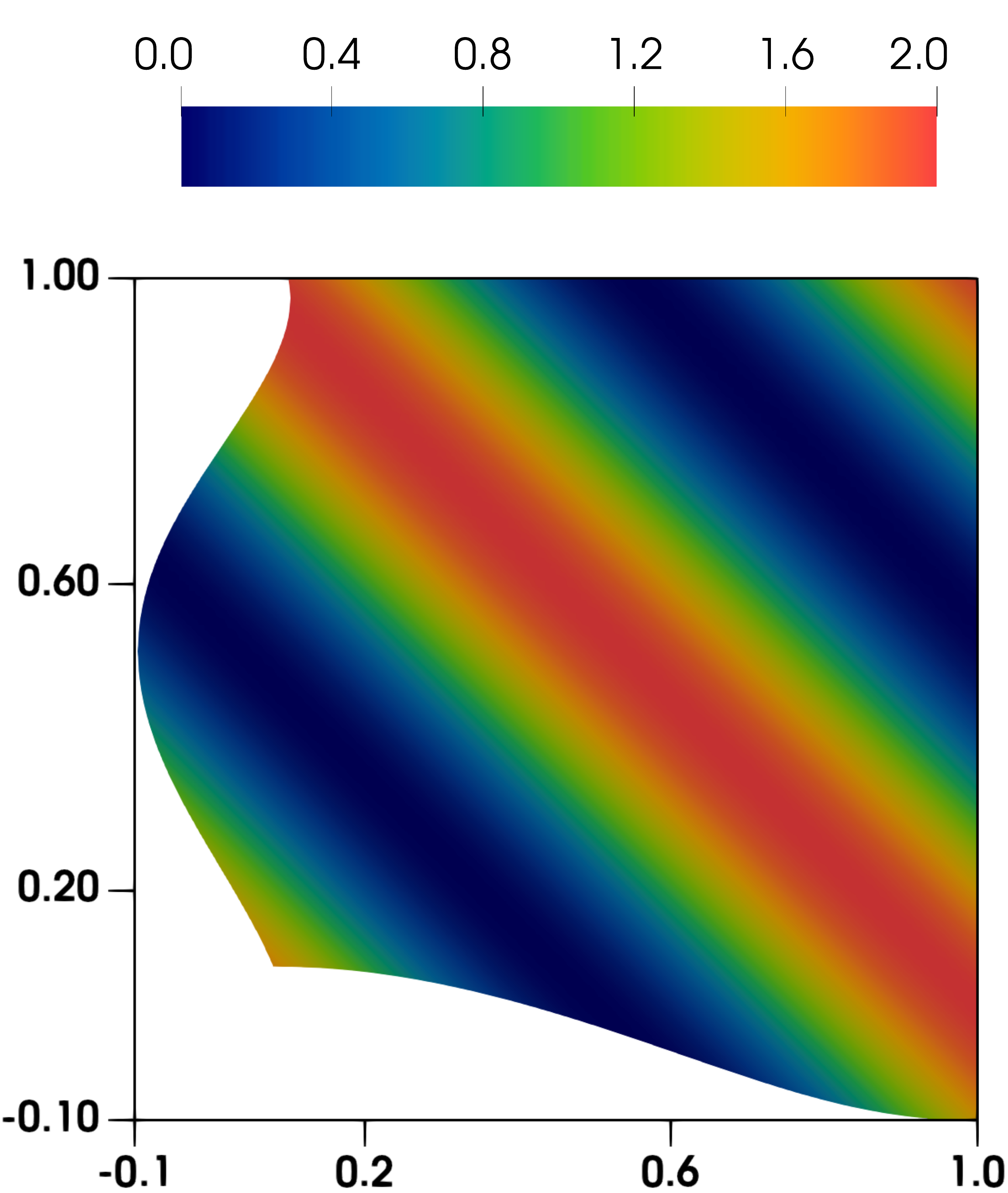}
    \caption{$\rho$ at $t=2.0$.}
    \label{smooth2_rho_16}
    \end{subfigure}
    \quad
    \begin{subfigure}{0.3\textwidth}
    \includegraphics[width=\linewidth]{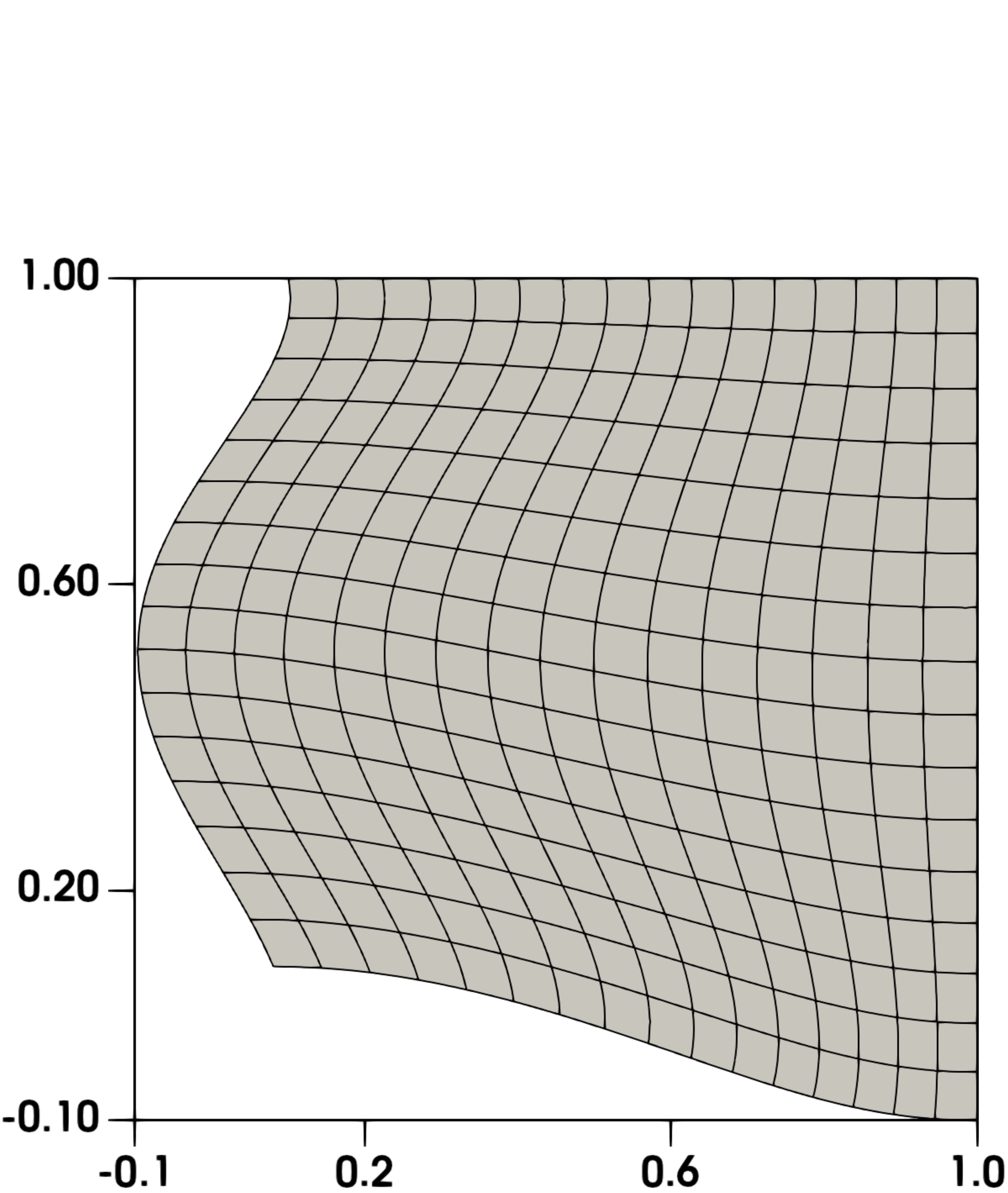}
    \caption{Curved mesh.}
    \label{smooth2_mesh_16}
    \end{subfigure}
    \caption{Sinusoidal smooth test: Results with $16^2$ elements.}
\end{figure}
We consider a test case in which the solution remains smooth for all time. The domain of computation (Figure~\ref{smooth2_mesh_16}) is a distorted square with two curved boundaries and two linear boundaries, obtained by the transformation map $[0,1]^2 \to \Omega$ given by
\[
    \begin{pmatrix}
        x\\
        y
    \end{pmatrix} = 
    \begin{pmatrix}
    \hat{x} + 0.1 \cos\left(\frac{\pi \hat{x}}{2} \right) \cos(2 \pi y)\\
    \hat{y} + 0.1 \cos\left(\frac{\pi \hat{y}}{2} \right) \cos\left(\pi \hat{x} \right)
    \end{pmatrix}.
\]
The fluid density, velocities, and pressure at the initial time are given by
\[
    \rho(x,y) = 1 + 0.999 \sin\big(2\pi (x+y)\big), \quad v_1(x,y)=v_2(x,y) = \frac{0.99}{\sqrt{2}}, \quad p(x,y) = 0.01.
\]
With time, the density advects as
\[
    \rho(x,y,t) = 1 + 0.999 \sin\left(2\pi \left((x+y) -\frac{0.99t}{\sqrt2}\right)\right),
\]
while the other quantities remain unchanged. The simulations are \revb{performed over the time interval $[0, 2]$}, with a uniform mesh and periodic boundaries using the equations of state~(\ref{eq: ID_eos},\ref{eq: TM_eos},\ref{eq: IP_eos},\ref{eq: RC_eos}). The grid convergence study is presented in Tables~\ref{table: sm1_N3_ID43}-\ref{table: sm1_N4_RC}. \reva{It is observed that as we refine the mesh, the numerical order of accuracy of the scheme for degree $N=3,4$ approaches $4,5$, respectively, which is the optimal rate we can expect.}

\begin{table}[!htbp]
    \centering
         \begin{tabular}{|l|c|c|c|c|}
    \hline
\multirow{2}{*}{Elements}
& \multicolumn{2}{|c|}{Error} 
& \multicolumn{2}{|c|}{Order} \\
\cline{2-5}
& $L_2$ & $L_\infty$ & $L_2$ & $L_\infty$ \\
\hline
$16^2$     &       $1.62251e-04$       &             $1.95286e-03$       &    -           &           -    \\ 
$32^2$     &       $1.01168e-05$       &           $1.26982e-04$       &   $4.00340$        &        $3.94290$ \\ 
$64^2$     &       $6.31309e-07$       &           $8.28910e-06$       &    $4.00226$        &       $3.93726$ \\ 
$128^2$     &       $3.94111e-08$       &           $5.31797e-07$       &   $4.00167$        &        $3.96227$ \\ 
$256^2$     &       $2.46514e-09$       &            $3.37867e-08$       &  $3.99886$       &        $3.97635$ \\
    \hline
    \end{tabular}
    \vspace{0.5em}
    \centering
    \caption{$\rho$ using $N=3$ and ID-EOS~\eqref{eq: ID_eos} with $\gamma = \frac{4}{3}$.}
    \label{table: sm1_N3_ID43}
    
        \begin{tabular}{|l|c|c|c|c|}
    \hline
\multirow{2}{*}{Elements}
& \multicolumn{2}{|c|}{Error} 
& \multicolumn{2}{|c|}{Order} \\
\cline{2-5}
& $L_2$ & $L_\infty$ & $L_2$ & $L_\infty$ \\
\hline
$16^2$     &       $1.61925e-04$       &              $1.95411e-03$       &    -           &          -    \\ 
$32^2$     &       $1.00844e-05$       &           $1.26979e-04$       &   $4.00513$        &        $3.94385$ \\
$64^2$     &       $6.29380e-07$       &           $8.30514e-06$       &   $4.00205$        &        $3.93444$ \\ 
$128^2$     &       $3.92832e-08$       &           $5.32753e-07$       &   $4.00195$        &        $3.96247$ \\ 
$256^2$     &       $2.45720e-09$       &            $3.38335e-08$       &  $3.99882$        &        $3.97694$ \\ 
    \hline
    \end{tabular}
    \vspace{0.5em}
    \centering
    \caption{$\rho$ using $N=3$ and ID-EOS~\eqref{eq: ID_eos} with $\gamma = \frac{5}{3}$.}
    \label{table: sm1_N3_ID53}

       \begin{tabular}{|l|c|c|c|c|}
    \hline
\multirow{2}{*}{Elements}
& \multicolumn{2}{|c|}{Error} 
& \multicolumn{2}{|c|}{Order} \\
\cline{2-5}
& $L_2$ & $L_\infty$ & $L_2$ & $L_\infty$ \\
\hline
$16^2$     &       $3.36359e-04$       &             $3.87423e-03$       &     -           &          -    \\ 
$32^2$     &       $1.61591e-04$       &          $1.99372e-03$       &    $1.05766$        &       $0.95845$ \\ 
$64^2$     &      $ 3.08720e-05$       &          $4.65147e-04$       &    $2.38797$        &        $2.09970$ \\ 
$128^2$     &       $1.98873e-06$       &           $3.91838e-05$       &   $3.95638$        &        $3.56936$ \\ 
$256^2$     &       $9.19875e-08$       &           $1.58153e-06$       &   $4.43427$        &       $4.63086$ \\ 
    \hline
    \end{tabular}
    \vspace{0.5em}
    \centering
    \caption{$\rho$ using $N=3$ and TM-EOS~\eqref{eq: TM_eos}.}
    \label{table: sm1_N3_TM}

    \centering
       \begin{tabular}{|l|c|c|c|c|}
    \hline
\multirow{2}{*}{Elements}
& \multicolumn{2}{|c|}{Error} 
& \multicolumn{2}{|c|}{Order} \\
\cline{2-5}
& $L_2$ & $L_\infty$ & $L_2$ & $L_\infty$ \\
\hline
$16^2$     &       $4.97481e-04$       &           $5.28327e-03$       &     -           &            -    \\ 
$32^2$     &       $2.09684e-04$       &            $2.55048e-03$       &  $1.24643$        &        $1.05066$ \\ 
$64^2$     &       $2.99267e-05$       &          $4.72164e-04$       &    $2.80871$        &        $2.43341$ \\ 
$128^2$     &       $1.75484e-06$       &           $3.33983e-05$       &   $4.09203$        &        $3.82144$ \\ 
$256^2$     &       $8.01745e-08$       &          $1.36774e-06$       &    $4.45205$        &        $4.60991$ \\ 
   \hline
    \end{tabular}
    \vspace{0.5em}
    \centering
    \caption{$\rho$ using $N=3$ and IP-EOS~\eqref{eq: IP_eos}.}
    \label{table: sm1_N3_IP}
       \begin{tabular}{|l|c|c|c|c|}
    \hline
\multirow{2}{*}{Elements}
& \multicolumn{2}{|c|}{Error} 
& \multicolumn{2}{|c|}{Order} \\
\cline{2-5}
& $L_2$ & $L_\infty$ & $L_2$ & $L_\infty$ \\
\hline
$16^2$     &       $2.76776e-04$       &            $3.23536e-03$       &     -           &           -    \\ 
$32^2$     &       $1.07216e-04$       &            $1.34685e-03$       &   $1.36820$        &       $1.26433$ \\ 
$64^2$     &       $1.91623e-05$       &            $3.11851e-04$       &    $2.48417$        &      $2.11066$ \\ 
$128^2$     &       $8.64451e-07$       &           $1.70217e-05$       &   $4.47034$        &        $4.19541$ \\ 
$256^2$     &       $3.79193e-08$       &           $6.94005e-07$       &    $4.51078$        &      $4.61628$ \\  
    \hline
    \end{tabular}
    \vspace{0.5em}
    \centering
    \caption{$\rho$ using $N=3$ and RC-EOS~\eqref{eq: RC_eos}.}
    \label{table: sm1_N3_RC}
\end{table}

\begin{table}[!htbp]
    \centering
       \begin{tabular}{|l|c|c|c|c|}
    \hline
\multirow{2}{*}{Elements}
& \multicolumn{2}{|c|}{Error} 
& \multicolumn{2}{|c|}{Order} \\
\cline{2-5}
& $L_2$ & $L_\infty$ & $L_2$ & $L_\infty$ \\
\hline
$16^2$     &       $1.23356e-05$       &             $1.42907e-04$       &    -           &           -    \\ 
$32^2$     &       $3.92584e-07$       &            $4.62840e-06$       &   $4.97368$        &       $4.94842$ \\ 
$64^2$     &       $1.23486e-08$       &            $1.46927e-07$       &   $4.99058$        &       $4.97735$ \\ 
$128^2$     &       $3.86735e-10$       &           $4.64005e-09$       &   $4.99686$        &        $4.98481$ \\  
    \hline
    \end{tabular}
    \vspace{0.5em}
    \centering
    \caption{$\rho$ using $N=4$ and ID-EOS~\eqref{eq: ID_eos} with $\gamma = \frac{4}{3}$.}
    \label{table: sm1_N4_ID43}

    \begin{tabular}{|l|c|c|c|c|}
    \hline
\multirow{2}{*}{Elements}
& \multicolumn{2}{|c|}{Error} 
& \multicolumn{2}{|c|}{Order} \\
\cline{2-5}
& $L_2$ & $L_\infty$ & $L_2$ & $L_\infty$ \\
\hline
$16^2$     &       $1.24027e-05$       &             $1.44156e-04$       &      -           &          -    \\ 
$32^2$     &       $3.95277e-07$       &             $4.67175e-06$       &  $4.97165$        &       $4.94752$ \\ 
$64^2$     &       $1.24212e-08$       &             $1.47995e-07$       &  $4.99199$        &       $4.98034$ \\ 
$128^2$     &       $3.89725e-10$       &             $4.72965e-09$       &  $4.99420$        &      $4.96767$ \\ 
    \hline
    \end{tabular}
    \vspace{0.5em}
    \centering
    \caption{$\rho$ using $N=4$ and ID-EOS~\eqref{eq: ID_eos} with $\gamma = \frac{5}{3}$.}
    \label{table: sm1_N4_ID53}
    
   \begin{tabular}{|l|c|c|c|c|}
    \hline
\multirow{2}{*}{Elements}
& \multicolumn{2}{|c|}{Error} 
& \multicolumn{2}{|c|}{Order} \\
\cline{2-5}
& $L_2$ & $L_\infty$ & $L_2$ & $L_\infty$ \\
\hline
$16^2$     &       $2.22749e-04$       &              $2.14626e-03$       &    -           &          -    \\ 
$32^2$     &       $7.52951e-05$       &             $1.02060e-03$       &  $1.56479$        &       $1.07241$ \\ 
$64^2$     &       $6.69363e-06$       &            $1.42431e-04$       &   $3.49170$        &       $2.84108$ \\ 
$128^2$     &       $2.31529e-07$       &            $5.89641e-06$       &   $4.85352$        &       $4.59428$ \\ 
$256^2$     &       $5.13373e-09$       &             $1.24935e-07$       &  $5.49504$        &       $5.56059$ \\ 
    \hline
    \end{tabular}
    \vspace{0.5em}
    \centering
    \caption{$\rho$ using $N=4$ and TM-EOS~\eqref{eq: TM_eos}.}
    \label{table: sm1_N4_TM}

\begin{tabular}{|l|c|c|c|c|}
    \hline
\multirow{2}{*}{Elements}
& \multicolumn{2}{|c|}{Error} 
& \multicolumn{2}{|c|}{Order} \\
\cline{2-5}
& $L_2$ & $L_\infty$ & $L_2$ & $L_\infty$ \\
\hline
$16^2$     &       $3.26051e-04$       &              $2.93905e-03$       &     -           &          -    \\ 
$32^2$     &       $8.40564e-05$       &             $1.16407e-03$       &   $1.95567$        &       $1.33618$ \\ 
$64^2$     &       $5.70429e-06$       &           $1.20648e-04$       &    $3.88124$        &        $3.27030$ \\ 
$128^2$     &       $1.69327e-07$       &            $4.19917e-06$       &   $5.07417$        &        $4.84456$ \\ 
$256^2$     &       $3.92909e-09$       &             $8.98232e-08$       &   $5.42947$        &       $5.54687$ \\ 
    \hline
    \end{tabular}
    \vspace{0.5em}
    \centering
    \caption{$\rho$ using $N=4$ and IP-EOS~\eqref{eq: IP_eos}.}
    \label{table: sm1_N4_IP}

\begin{tabular}{|l|c|c|c|c|}
    \hline
\multirow{2}{*}{Elements}
& \multicolumn{2}{|c|}{Error} 
& \multicolumn{2}{|c|}{Order} \\
\cline{2-5}
& $L_2$ & $L_\infty$ & $L_2$ & $L_\infty$ \\
\hline
$16^2$     &       $1.57772e-04$       &              $1.51796e-03$       &    -           &          -    \\ 
$32^2$     &       $4.87677e-05$       &             $6.65998e-04$       &  $1.69384$        &       $1.18855$ \\ 
$64^2$     &       $4.77744e-06$       &            $1.05038e-04$       &   $3.35162$        &       $2.66461$ \\ 
$128^2$     &       $9.52441e-08$       &             $2.71482e-06$       &  $5.64846$        &       $5.27391$ \\ 
$256^2$     &       $2.02339e-09$       &             $5.15864e-08$       &   $5.55678$        &      $5.71772$ \\ 
    \hline
    \end{tabular}
    \vspace{0.5em}
    \centering
    \caption{$\rho$ using $N=4$ and RC-EOS~\eqref{eq: RC_eos}.}
    \label{table: sm1_N4_RC}
\end{table}

\subsection{Isentropic vortex test}\label{sec: isentropic_vortex}
\begin{figure}[]
    \centering
    \begin{subfigure}{0.3\textwidth}
    \includegraphics[width=\linewidth]{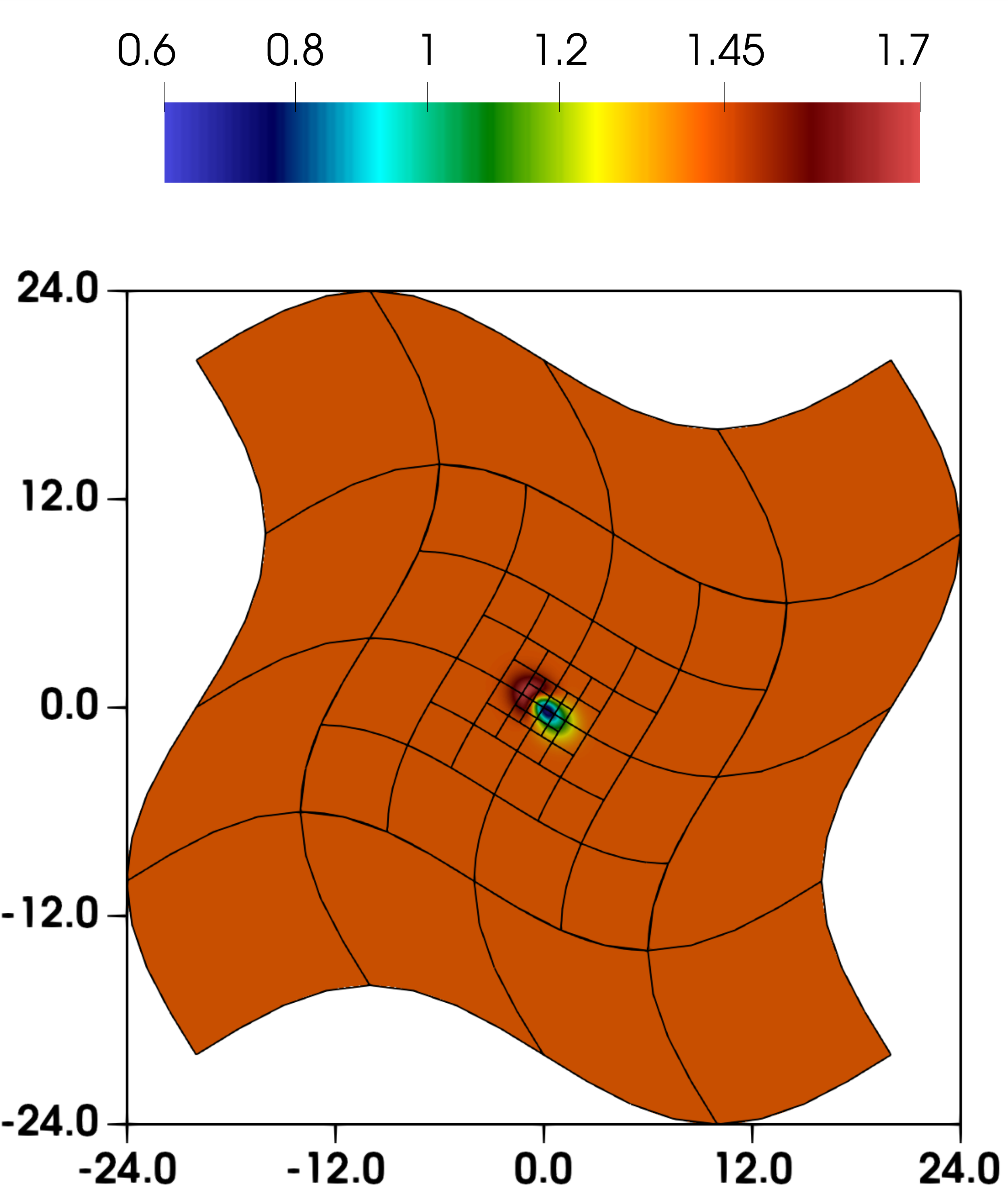}
    \vspace*{0.42cm}
    \caption{$D$ with AMR.}
    \label{fig: IV_domain}
    \end{subfigure}
    \quad
        \begin{subfigure}{0.3\textwidth}
    \includegraphics[width=\linewidth]{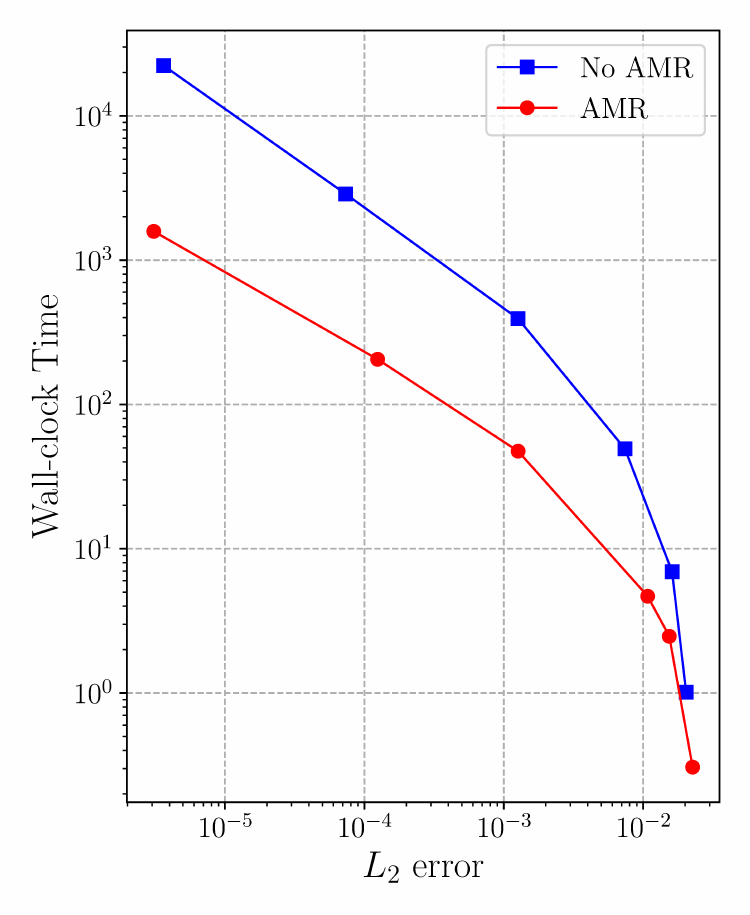}
    \caption{$L_2$ error versus wall-clock time in seconds.}
    \label{fig: l2_error_wctime}
    \end{subfigure}
    \quad
    \begin{subfigure}{0.3\textwidth}
    \includegraphics[width=\linewidth]{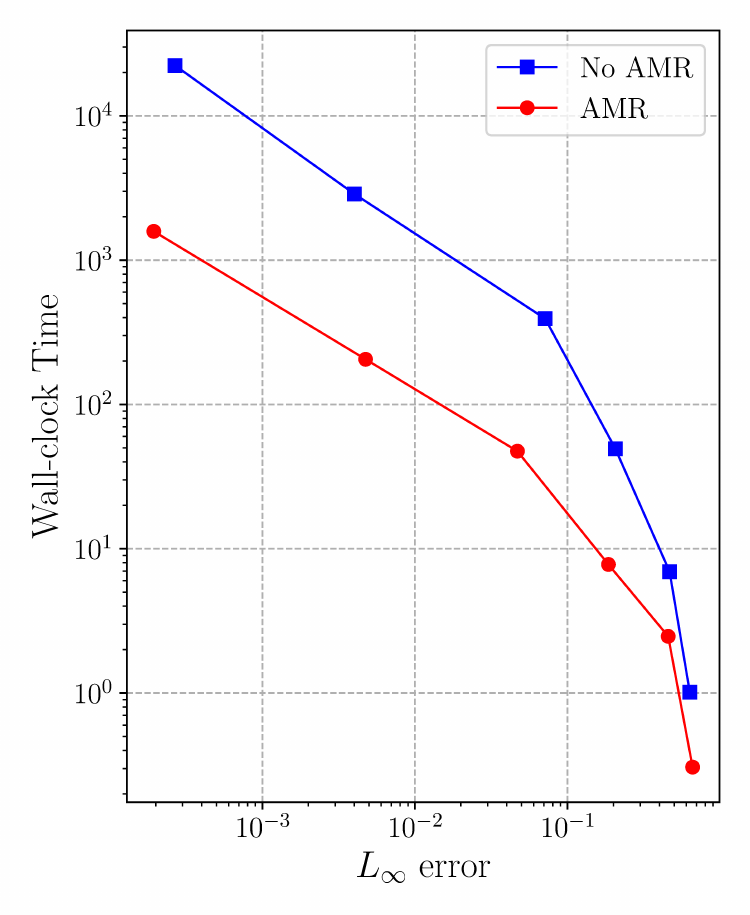}
    \caption{$L_\infty$ error versus wall-clock time in seconds.}
    \label{fig: linf_error_wctime}
    \end{subfigure}
    \caption{Isentropic vortex test: Results with density $D$ at $t=0.2$.}
    \label{fig: isentropic_vortex}
\end{figure}

We now consider the isentropic vortex test~\cite{ling2019physical, duan2019high} for the case of ID-EOS~\eqref{eq: ID_eos}, where a vortex travels diagonally. The initial condition is given by
\begin{align*}
    \rho(x,y) &= \left(1 - a\exp(1-x_0^2 - y_0^2)\right)^{\frac{1}{\gamma -1}}, \quad p(x,y) = \rho(x,y)^\gamma,\\
    v_d(x,y) &= \frac{1}{1-\frac{w(v_1'+v_2')}{\sqrt{2}}}\left(\frac{v_d'}{\phi} - \frac{w}{\sqrt{2}}+\frac{\phi w^2}{2(\phi +1)}(v_1'+v_2') \right),
\end{align*}
with $d=1,2$, and
\begin{align*}
    &a = \frac{(\gamma - 1)}{\gamma}\frac{\epsilon^2}{8\pi^2},\quad
    x_0 = x + \frac{\phi -1}{2}(x+y), \quad y_0 = y + \frac{\phi -1}{2}(x+y), \quad \phi = \frac{1}{\sqrt{1-w^2}},\\
    &(v_1',v_2')=(-y_0, x_0) \sqrt{\frac{b}{1+b(x_0^2 + y_0^2)}},\qquad b=\frac{2\gamma a \exp(1-x_0^2 - y_0^2)}{2\gamma - 1 -\gamma a \exp(1-x_0^2 - y_0^2)}.
\end{align*}
The parameters $\epsilon, w$ are taken as $ 5, 0.5\sqrt{2}$, respectively, which defines the vortex strength. The domain of computation (Figure~\ref{fig: IV_domain}) is taken by transforming the square $[0,1]^2$ with the map $[0,1]^2 \to \Omega$ given by
\[
    \begin{pmatrix}
        x\\
        y
    \end{pmatrix} = 
    \begin{pmatrix}
    40 \hat{x} -  4 \sin(2 \pi \hat{y})-20\\
    40 \hat{y} +  4 \sin(2 \pi \hat{x})-20
    \end{pmatrix}.
\]
The simulations are \revb{carried out} with periodic boundaries until time $t=80$. At this point, the solution matches the initial solution profile. The time step computation is done as in~\cite{basak2025bound}. The grid-convergence study is presented in Tables~\ref{table: IV_N3_ID53}-\ref{table: IV_N4_ID53} with $\gamma= \frac{5}{3}$, and it is observed that \reva{the scheme attains the numerical order of accuracy close to $4,5$ as we refine the mesh, for polynomial degrees $N=3,4$, respectively.}

\begin{table}[!htbp]
    \centering
       \begin{tabular}{|l|c|c|c|c|}
    \hline
\multirow{2}{*}{Elements}
& \multicolumn{2}{|c|}{Error} 
& \multicolumn{2}{|c|}{Order} \\
\cline{2-5}
& $L_2$ & $L_\infty$ & $L_2$ & $L_\infty$ \\
\hline
$10^2$     &       $2.03338e-02$       & $6.33295e-01$      &     -           &           -    \\ 
$20^2$     &       $1.61197e-02$       & $4.66832e-01$      &  $0.33506$        &        $0.43997$ \\ 
$40^2$     &       $7.39973e-03$       &  $2.05937e-01$     &  $1.12328$        &         $1.18070$ \\ 
$80^2$     &       $1.26624e-03$       &   $7.13690e-02$    &  $2.54693$        &          $1.52884$ \\ 
$160^2$     &       $7.34580e-05$       &  $4.01054e-03$    &   $4.10748$        &         $4.15343$ \\ 
$320^2$     &       $3.62265e-06$       &   $2.67802e-04$   &   $4.34180$        &           $3.90456$ \\  
    \hline
    \end{tabular}
    \vspace{0.5em}
    \centering
    \caption{$\rho$ using $N=3$ and ID-EOS~\eqref{eq: ID_eos} with $\gamma = \frac{5}{3}$.}
    \label{table: IV_N3_ID53}
    
    \begin{tabular}{|l|c|c|c|c|}
    \hline
\multirow{2}{*}{Elements}
& \multicolumn{2}{|c|}{Error} 
& \multicolumn{2}{|c|}{Order} \\
\cline{2-5}
& $L_2$ & $L_\infty$ & $L_2$ & $L_\infty$ \\
\hline
$10^2$     &       $1.95474e-02$      &           $5.90090e-01$       &         -           &         -    \\ 
$20^2$     &       $1.08039e-02$       &            $2.88201e-01$       &   $0.85542$        &        $1.03386$ \\ 
$40^2$     &       $3.29655e-03$       &             $2.11061e-01$       &   $1.71253$        &      $0.44942$ \\ 
$80^2$     &       $1.72426e-04$       &            $1.05593e-02$       &  $4.25691$        &        $4.32107$ \\ 
$160^2$    &       $3.11463e-06$       &          $3.66146e-04$       &   $5.79077$        &         $4.84995$ \\  
    \hline
    \end{tabular}
    \vspace{0.5em}
    \centering
    \caption{$\rho$ using $N=4$ and ID-EOS~\eqref{eq: ID_eos} with $\gamma = \frac{5}{3}$.}
    \label{table: IV_N4_ID53}
\end{table}

The solution in this test case is constant except in a small region around the vortex. Thus, a uniform mesh to resolve the vortex wastes a lot of cells in the constant regions where a fine mesh is not required. Hence, we use this test case to show benefits of adaptive mesh refinement in terms of wall-clock time for the simulation. The Figures~\ref{fig: l2_error_wctime},~\ref{fig: linf_error_wctime} show the graph of $L_2$ and $L_\infty$ errors versus wall-clock times (in seconds) of simulations, respectively; performed with uniform mesh and adaptive mesh with the AMR indicator~\eqref{eq: amr_indicator}. It is clear from the figures that, to achieve lower error, the wall-clock time for the simulations with a uniform mesh is larger compared to those with AMR, as the refinement levels increase. For example, on the finest meshes used, the simulation time for uniform mesh is $22313$ seconds while that with AMR is $1583$ seconds, showing a speedup of almost $14$ is achieved with AMR.

\subsection{Blast test}
\begin{figure}[]
    \centering
    \begin{subfigure}{0.31\textwidth}
    \includegraphics[width=\linewidth]{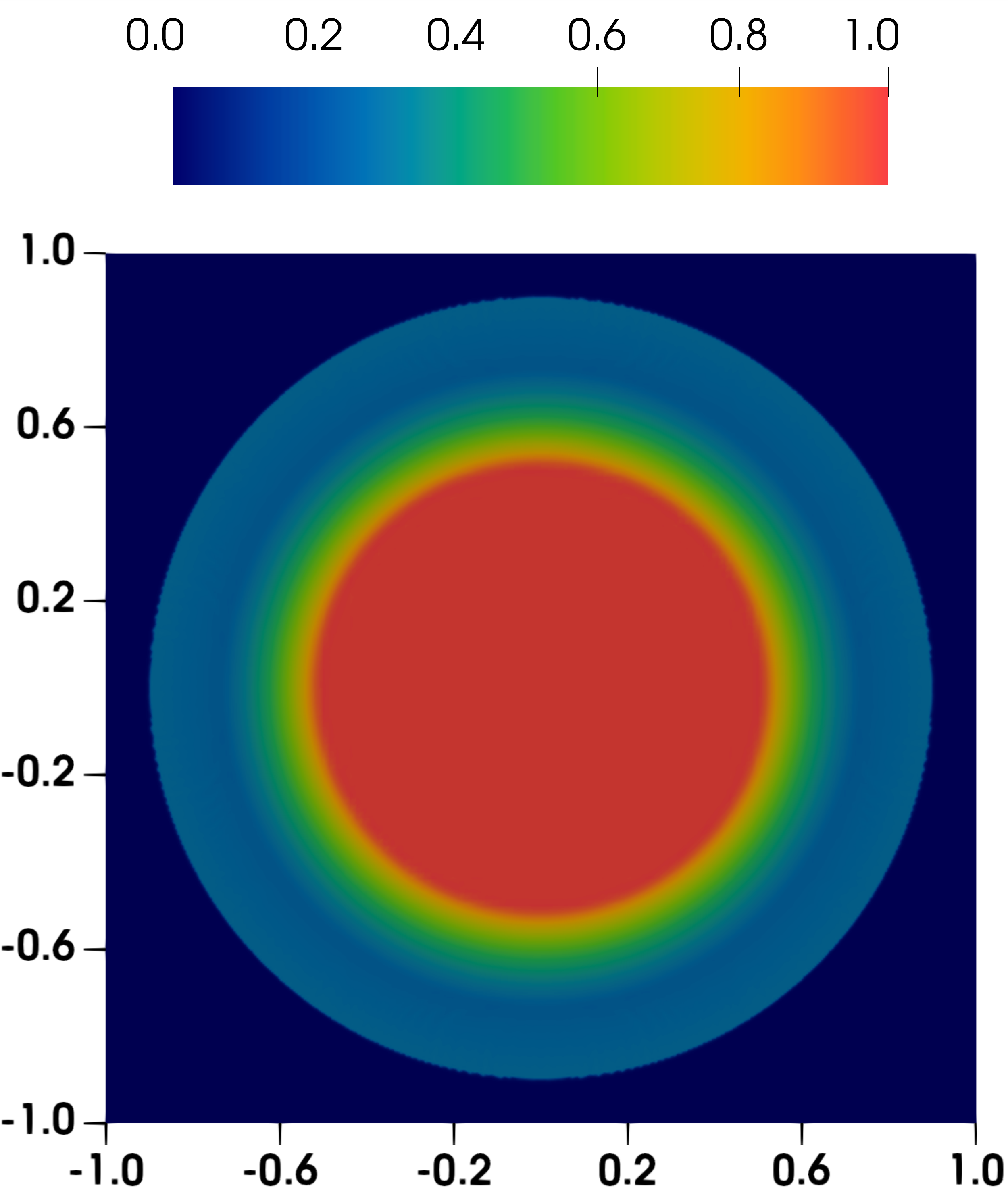}
    \caption{$\rho$ with AMR.}
    \label{fig:blast_rho}
    \end{subfigure}
    \vspace{0.1cm}
    \quad
    \begin{subfigure}{0.31\textwidth}
    \includegraphics[width=\linewidth]{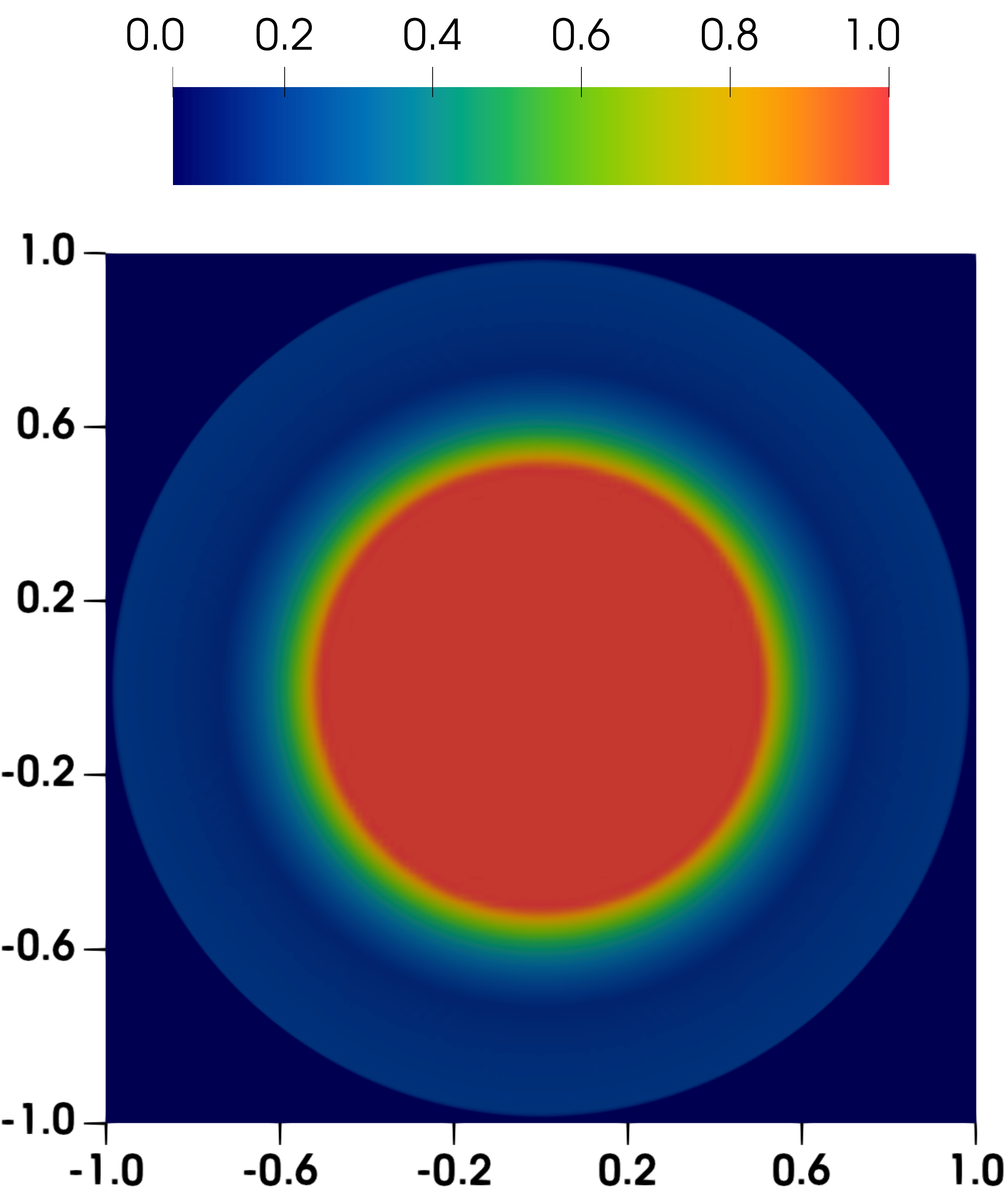}
    \caption{$p$ with AMR.}
    \label{fig:blast_p}
    \end{subfigure}
    \vspace{0.1cm}
    \quad
    \begin{subfigure}{0.31\textwidth}
    \includegraphics[width=\linewidth]{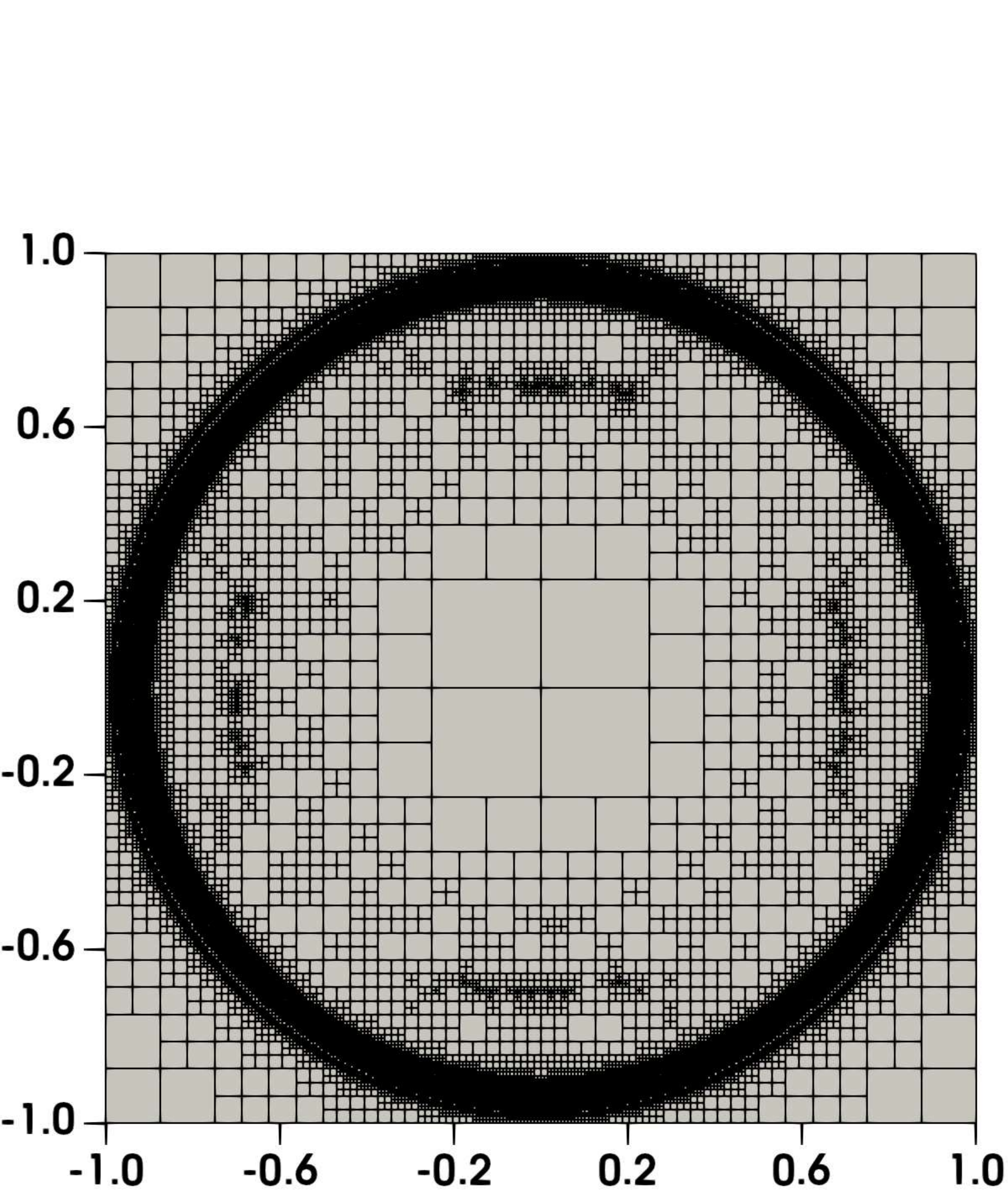}
    \caption{Adaptively refined mesh.}
    \label{fig:blast_mesh}
    \end{subfigure}
    \begin{subfigure}{0.31\textwidth}
    \includegraphics[width=\linewidth]{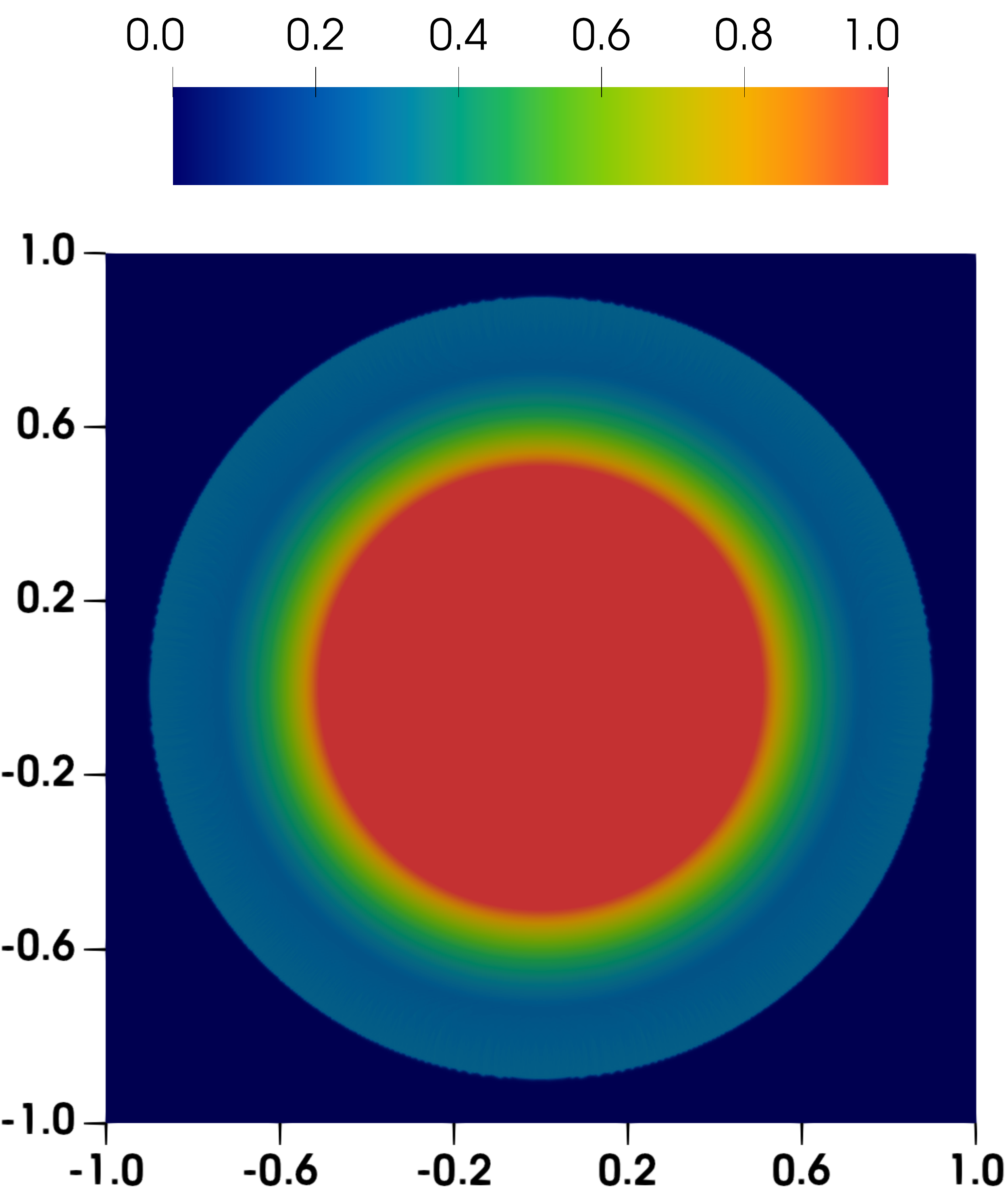}
    \caption{$\rho$ with uniform mesh.}
    \label{fig:blast_rho_uniform_mesh}
    \end{subfigure}
    \quad
    \begin{subfigure}{0.31\textwidth}
    \includegraphics[width=\linewidth]{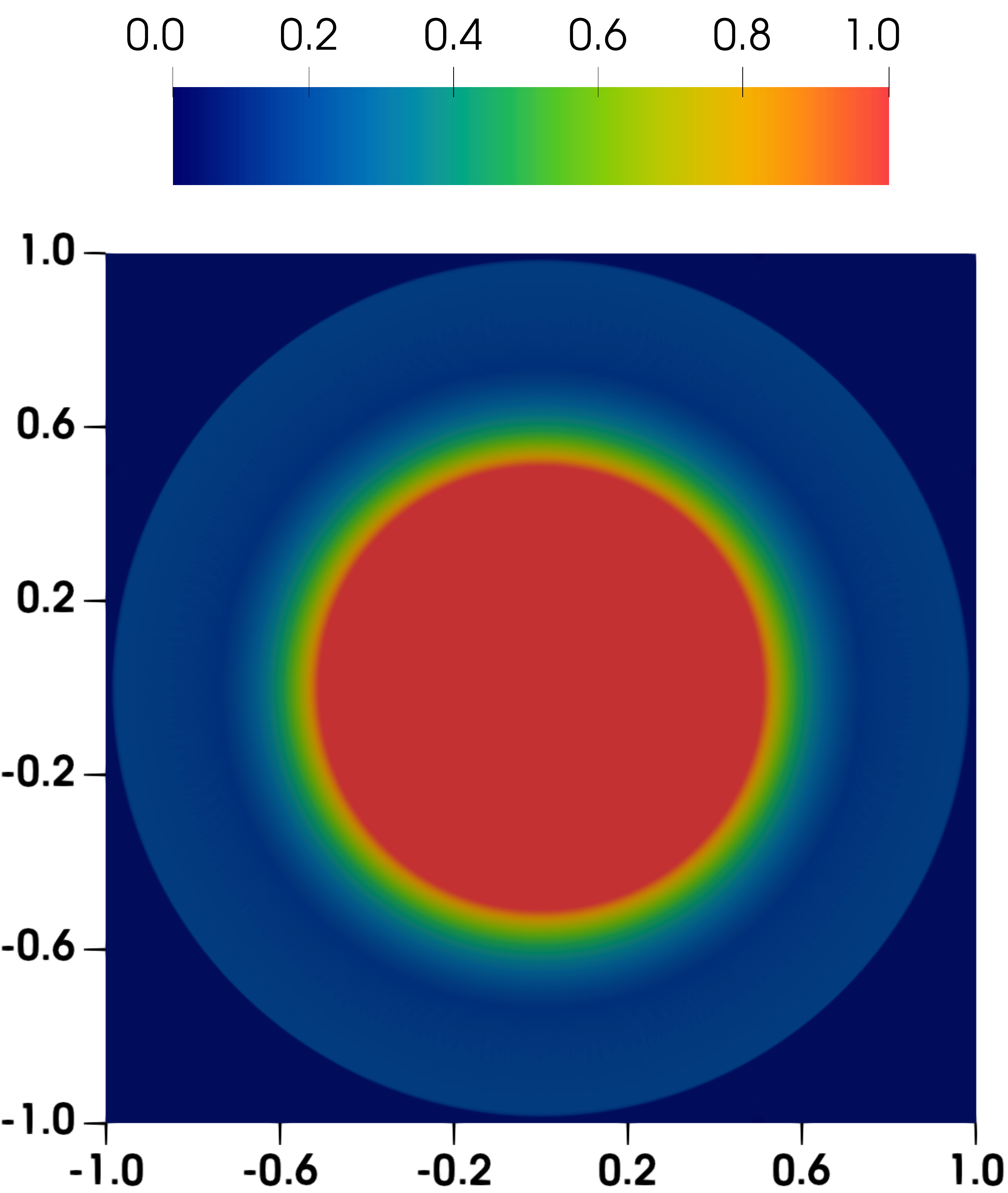}
    \caption{$p$ with uniform mesh.}
    \label{fig:blast_p_uniform_mesh}
    \end{subfigure}
    \quad 
    \begin{subfigure}{0.28\textwidth}
    \includegraphics[width=\linewidth]{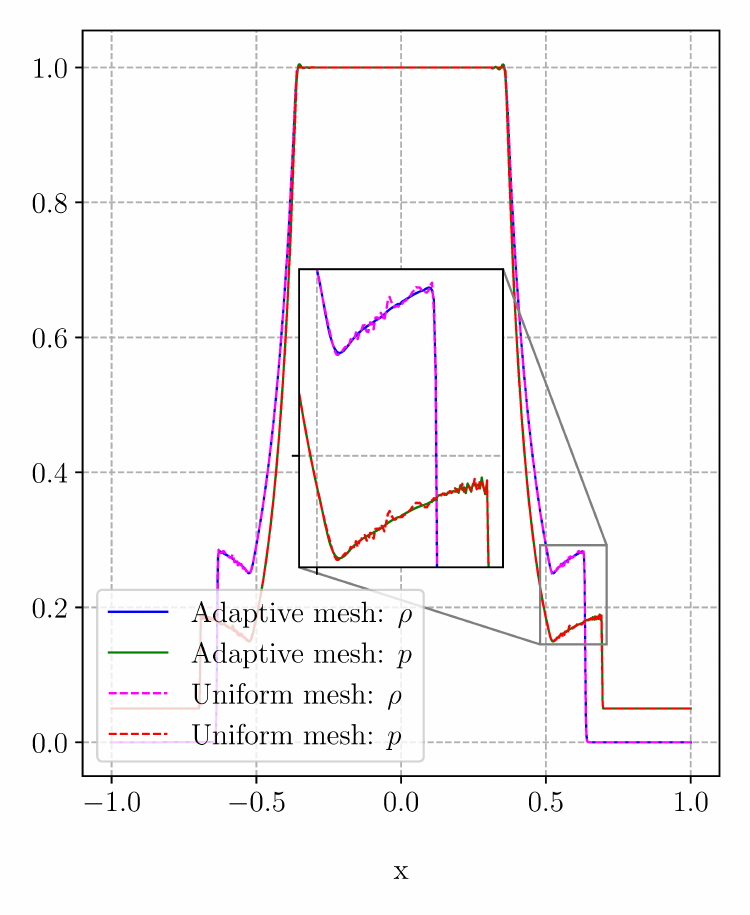}
    \caption{Cut plot from lower left to upper right.}
    \label{fig: blast_cut_plot}
    \end{subfigure}
    \caption{Blast test: Results at time $t=0.35$.}
    \label{fig:blast}
\end{figure}

This test is similar to the tests in~\cite{ling2019physical}. We use it to demonstrate the mesh adaptivity capability of the scheme, as it has a shock, a rarefaction, and a contact discontinuity in a  circular shape. Initially, a discontinuity is placed at a radius of $0.5$ with the initial conditions of the form
\[
    (\rho, v_1, v_2, p) = \begin{cases}
        (1.0, 0.0, 0.0, 1.0), \quad &\text{if } x^2 +y^2 < 0.5\\
        (\revb{10^{-12}}, 0.0, 0.0, 0.05) , \quad &\text{if } x^2 +y^2 > 0.5
    \end{cases}
\]
in the domain $[-1, 1]^2$. All the boundaries of the computational domain are taken as artificial boundaries (Section~\ref{sec: artificial_boundary}) with the components in $\mb{u}^{\text{set}}$ as $(\revb{10^{-12}}, 0.0, 0.0, 0.05)$. The rarefaction wave in the solution moves towards the center with time, and the shock and contact discontinuity move away from the center. The simulation is \revb{performed} with AMR as well as with uniform mesh \revb{in the time interval $[0, 0.35]$} with RC-EOS~\eqref{eq: RC_eos}, and $N=4$. \revb{It is observed that at final time, the minimum densities in the domain reach approximately $9.94 \times 10^{-13}$, and $3.39 \times 10^{-13}$, with AMR and uniform mesh, respectively.} The results are shown in Figure~\ref{fig:blast}. For the AMR case, the mesh is allowed to be refined and coarsened at each temporal step, with the AMR indicator~\eqref{eq: amr_indicator} and setting 
\[
    (\texttt{base\_level}, \quad \texttt{med\_level}, \quad \texttt{max\_level}) = (0, 3, 9), \qquad (\epsilon_1, \epsilon_2) = (0.01, 0.1),
\]
in the AMR controller~\eqref{eq: amr_controller}. Here, $\epsilon_1, \epsilon_2$ are the thresholds in the AMR controller and $\texttt{base\_level}= 0$ corresponds to $1\times 1$ mesh. The uniform mesh result is with $2^9 \times 2^9$ elements, which is equivalent to the maximum refinement level in AMR. It is clear from the figure that the AMR indicator is working effectively, capturing the moving contact discontinuity and shock wave with a highly refined mesh, and the smooth region at the center with a coarser mesh. The results with AMR are in good agreement with the uniform mesh result. To show the comparison between the AMR result and the uniform mesh result explicitly, we have presented a cut-plot along the diagonal from the lower left to the upper right corner of the domain by projecting onto the $x$-axis in Figure~\ref{fig: blast_cut_plot}. The final mesh with AMR has $\revb{43093}$ elements in the domain, while the uniform mesh has $262144$ elements. It is also worth mentioning that the wall-clock time for the simulation with AMR is $\revb{1349}$ seconds, while with the uniform mesh it is $\revb{5970}$ seconds. This shows a significant improvement of using AMR over the uniform mesh.

\subsection{Shock vortex test}
In this problem, a vortex similar to the isentropic vortex in Section~\ref{sec: isentropic_vortex}  interacts with a shock. This test was first used in a non-relativistic setting in~\cite{pao1981numerical}. Similar tests for the relativistic cases can be found in~\cite{balsara2016subluminal, duan2019high, cai2024provably}.

The vortex is the same as in Section~\ref{sec: isentropic_vortex} except with strength $w=0.9$, and a shock wave is placed at $x=-6$ with post-shock state 
\[
    (\rho, v_1, v_2, p) = (4.891497310766981, -0.388882958251919, 0.0, 11.894863258311670).
\]
We take the computational domain as $[-17, 3] \times [-5, 5]$ with $y=\pm 5$ as solid walls. The boundaries $x=-17, 3$ are taken as artificial boundaries (Section~\ref{sec: artificial_boundary}) with $\mb{u}^{\text{set}}$ as the solution near the corresponding boundary (inside the domain) at initial time. The simulation is \revb{performed} with ID-EOS~\eqref{eq: ID_eos} using $\gamma = 1.4$ and $N=4$ \revb{up to} $t=19$. The mesh is adapted according to the AMR indicator~\eqref{eq: amr_indicator} with the three-level controller~\eqref{eq: amr_controller} having
\[
    (\texttt{base\_level}, \quad \texttt{med\_level}, \quad \texttt{max\_level}) = (0,7,9), \qquad (\epsilon_1, \epsilon_2) = (0.002, 0.003).
\]
Here, $\epsilon_1, \epsilon_2$ are the thresholds in the AMR controller, and $\texttt{base\_level} = 0$ corresponds to $2\times 1$ mesh, which is allowed to be refined and coarsened at each time step.

\begin{figure}[]
    \centering
    \begin{subfigure}{0.48\textwidth}
    \includegraphics[width=\linewidth]{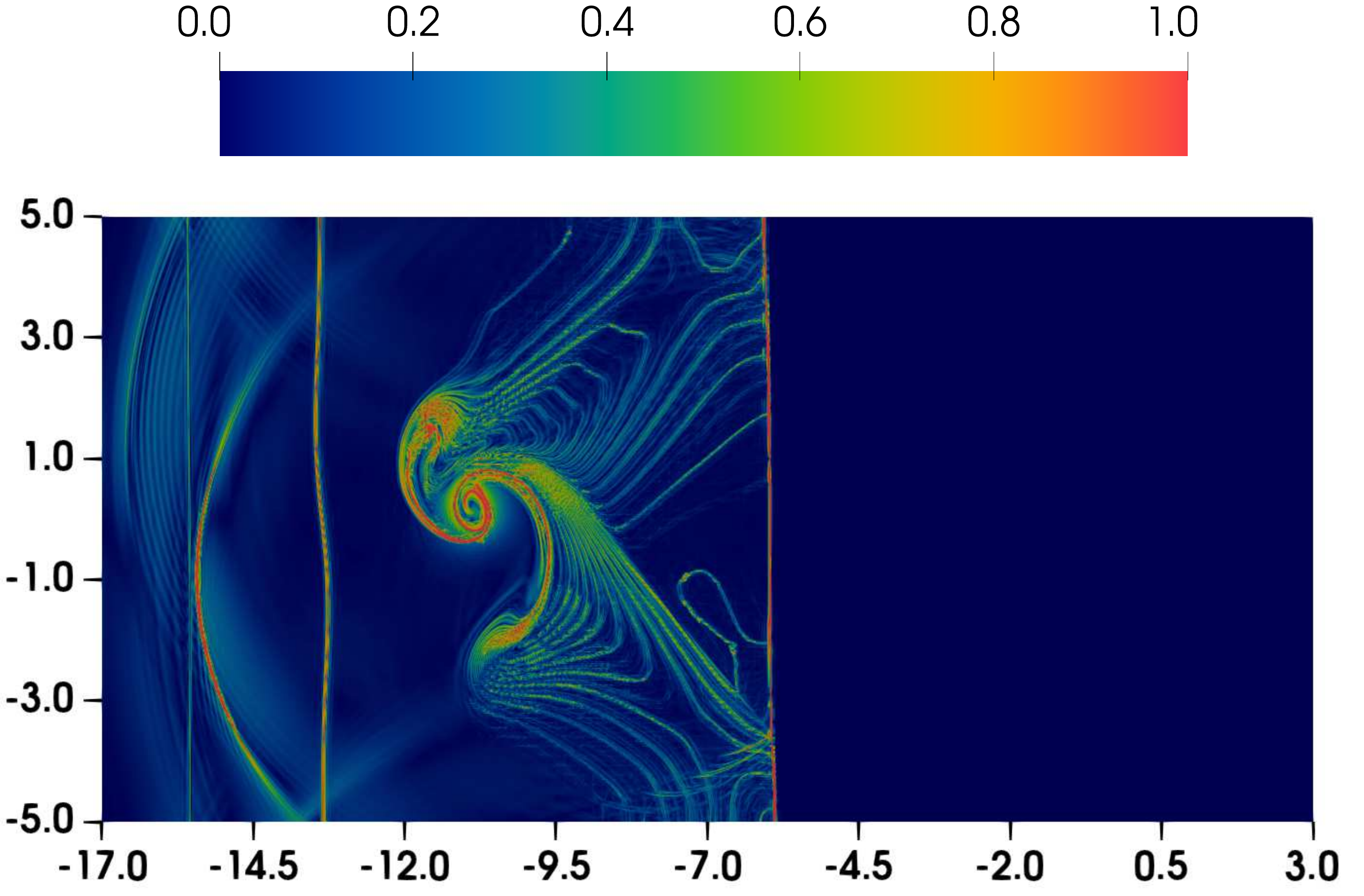}
    \caption{$\log_{10}(1+\nabla\rho)$ with AMR.}
    \label{fig:sv_rho_adaptive}
    \end{subfigure}
    \vspace{0.1cm}
    \quad
    \begin{subfigure}{0.48\textwidth}
    \includegraphics[width=\linewidth]{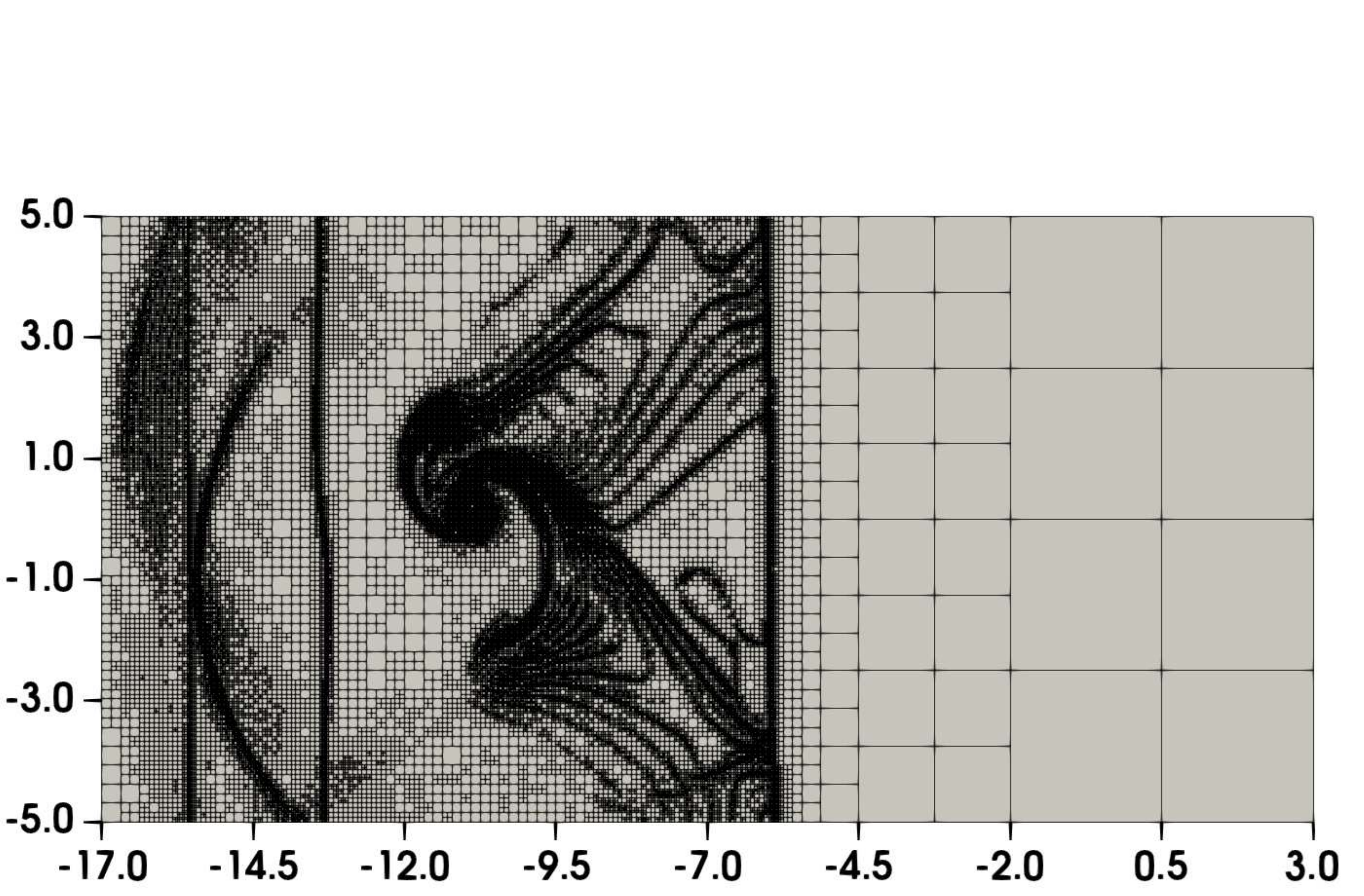}
    \caption{Adaptively refined mesh.}
    \label{fig:sv_mesh_adaptive}
    \end{subfigure}
    \vspace{0.1cm}
    \quad
    \begin{subfigure}{0.48\textwidth}
    \includegraphics[width=\linewidth]{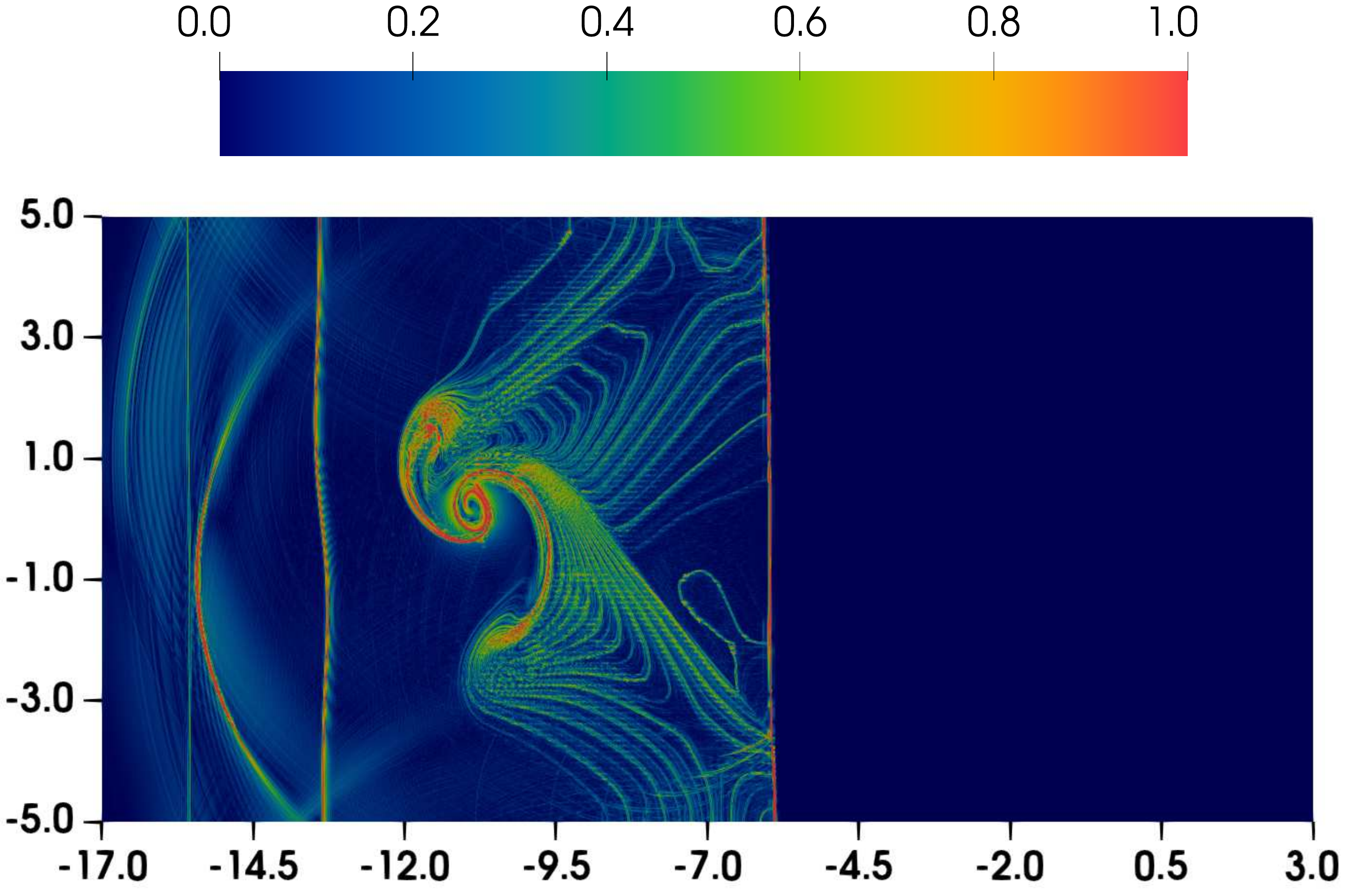}
    \caption{$\log_{10}(1+\nabla\rho)$ with uniform mesh.}
    \label{fig:sv_rho_uniform}
    \end{subfigure}
    \caption{Shock vortex test: Results at time $t=19$.}
    \label{fig:sv}
\end{figure}

Figure~\ref{fig:sv_rho_adaptive} shows the $\log_{10}(1+\nabla\rho)$ profile with AMR at the final time, and the corresponding mesh is shown in Figure~\ref{fig:sv_mesh_adaptive}. We have also shown the result with a uniform mesh with resolution equivalent to $\texttt{max\_level} = 9$ in Figure~\ref{fig:sv_rho_uniform}. We observe from the figure that the results with AMR and uniform mesh are in good agreement, capturing the linear and non-linear waves propagating in the domain. The AMR indicator is performing effectively, making the mesh resolution higher near the sharp structures in the solution, while coarsening the mesh in the smoother regions. The final mesh with AMR has $71708$ elements while the uniform mesh has $524288$ elements. Consequently, the computational cost decreases significantly with AMR. The wall-clock time for the simulation with AMR is $130520$ seconds, which is substantially less than the wall-clock time taken on the uniform mesh, which is $833895$ seconds.

\subsection{Bubble shock test}
\begin{figure}[]
    \centering
    \begin{subfigure}{0.49\textwidth}
    \includegraphics[width=\linewidth]{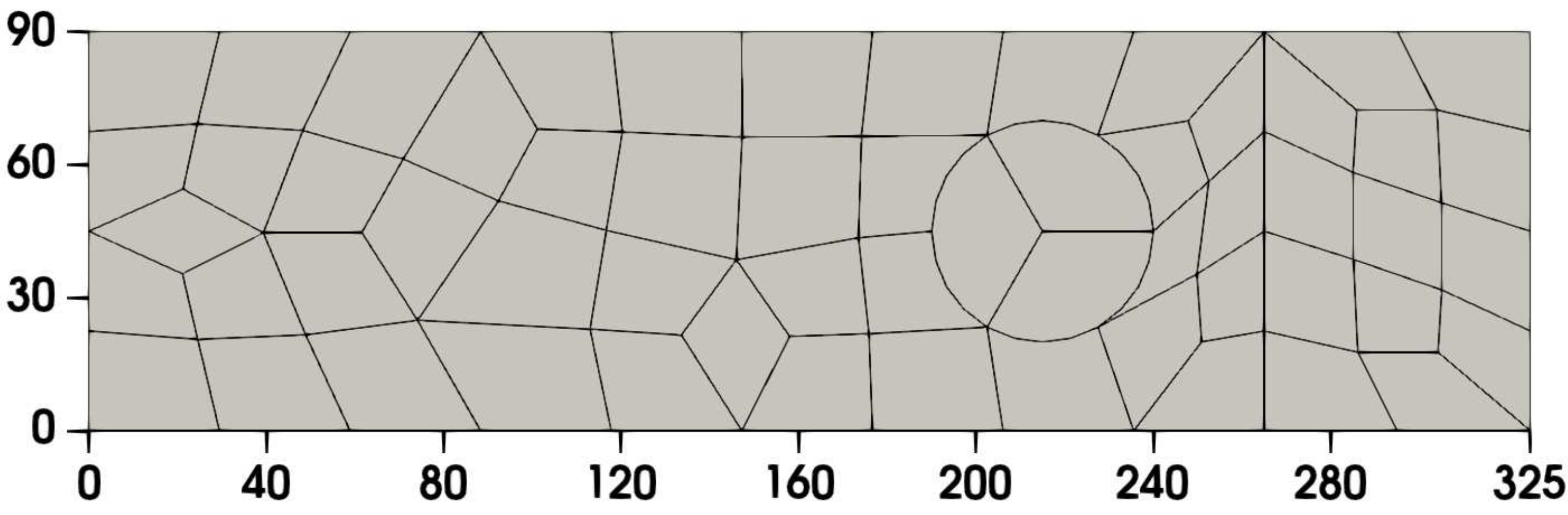}
    \caption{Mesh with $\texttt{base\_level} = 0$.}
    \label{fig: bs_mesh_bsl0}
    \end{subfigure}
    \qquad
       \vspace*{0.4cm} \hspace*{0.15cm}
          \begin{subfigure}{0.49\textwidth}
    \includegraphics[width=\linewidth]{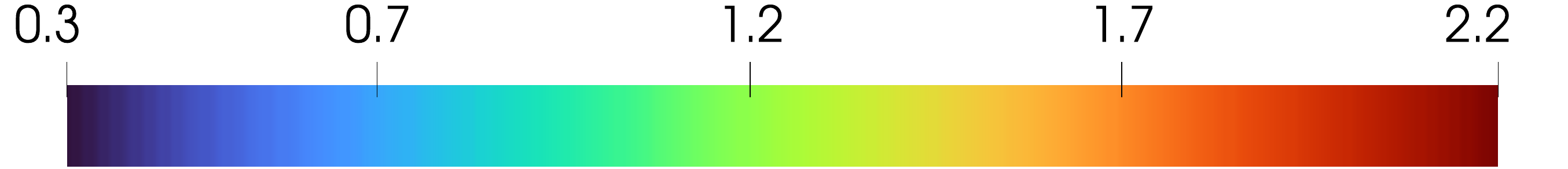}
    \end{subfigure}
        \begin{subfigure}{0.49\textwidth}
    \includegraphics[width=\linewidth]{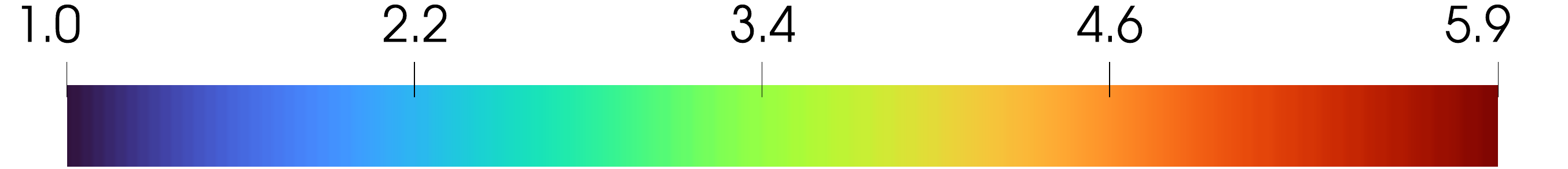}
    \end{subfigure}
    \begin{subfigure}{0.49\textwidth}
    \includegraphics[width=\linewidth]{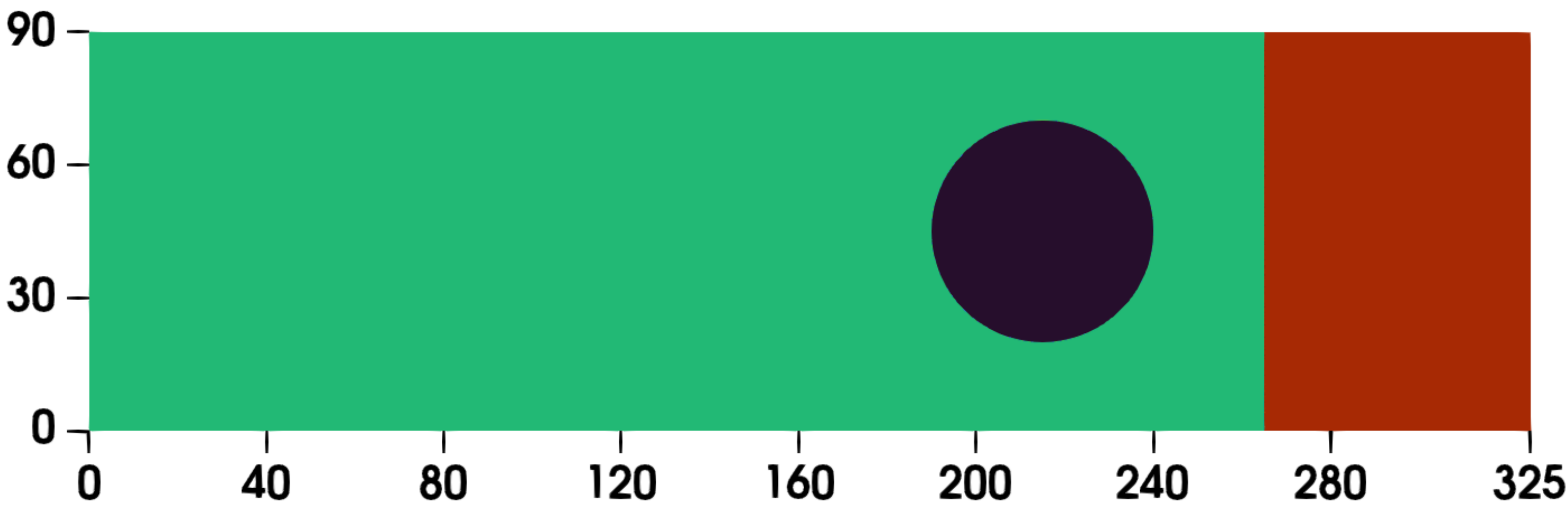}
    \caption{$\rho$ at $t=0$ for Case-I.}
    \label{fig: bs1_rho_t0}
    \end{subfigure}
    \begin{subfigure}{0.49\textwidth}
    \includegraphics[width=\linewidth]{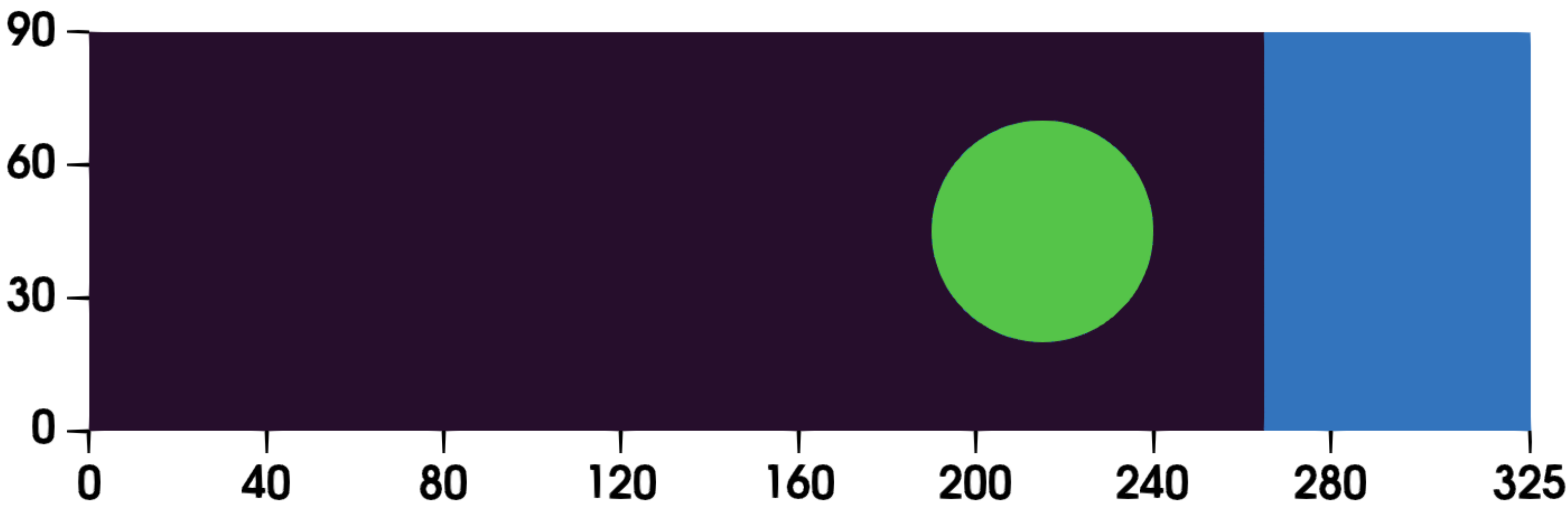}
    \caption{$\rho$ at $t=0$ for Case-II.}
    \label{fig: bs2_rho_t0}
    \end{subfigure}
    \begin{subfigure}{0.49\textwidth}
    \includegraphics[width=\linewidth]{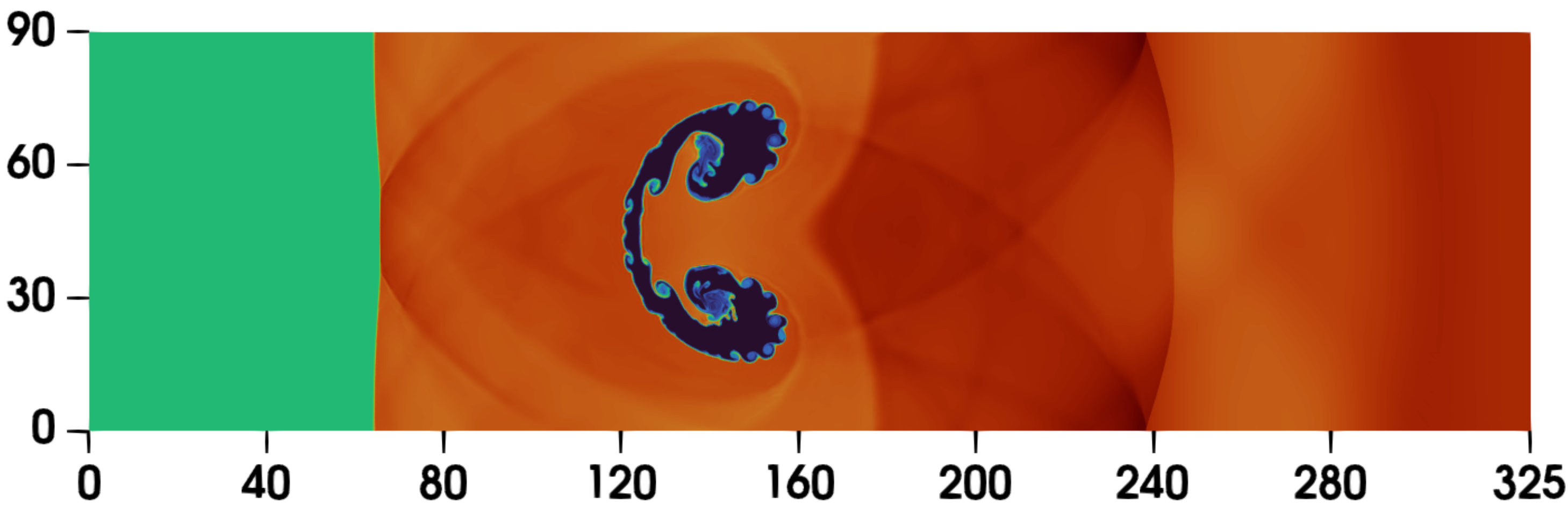}
    \caption{$\rho$ at $t=450$ with AMR for Case-I.}
    \label{fig: bs1_rho_tf_amr}
    \end{subfigure}
        \begin{subfigure}{0.49\textwidth}
    \includegraphics[width=\linewidth]{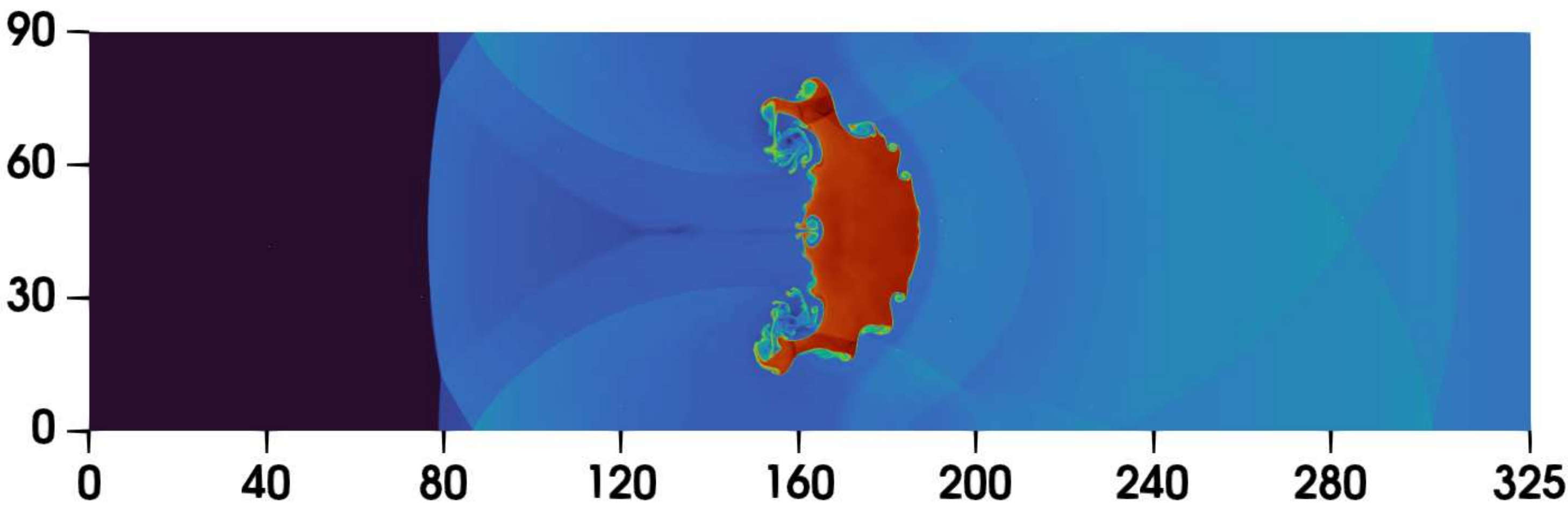}
    \caption{$\rho$ at $t=450$ with AMR for Case-II.}
    \label{fig: bs2_rho_tf_amr}
    \end{subfigure}
           \begin{subfigure}{0.49\textwidth}
    \includegraphics[width=\linewidth]{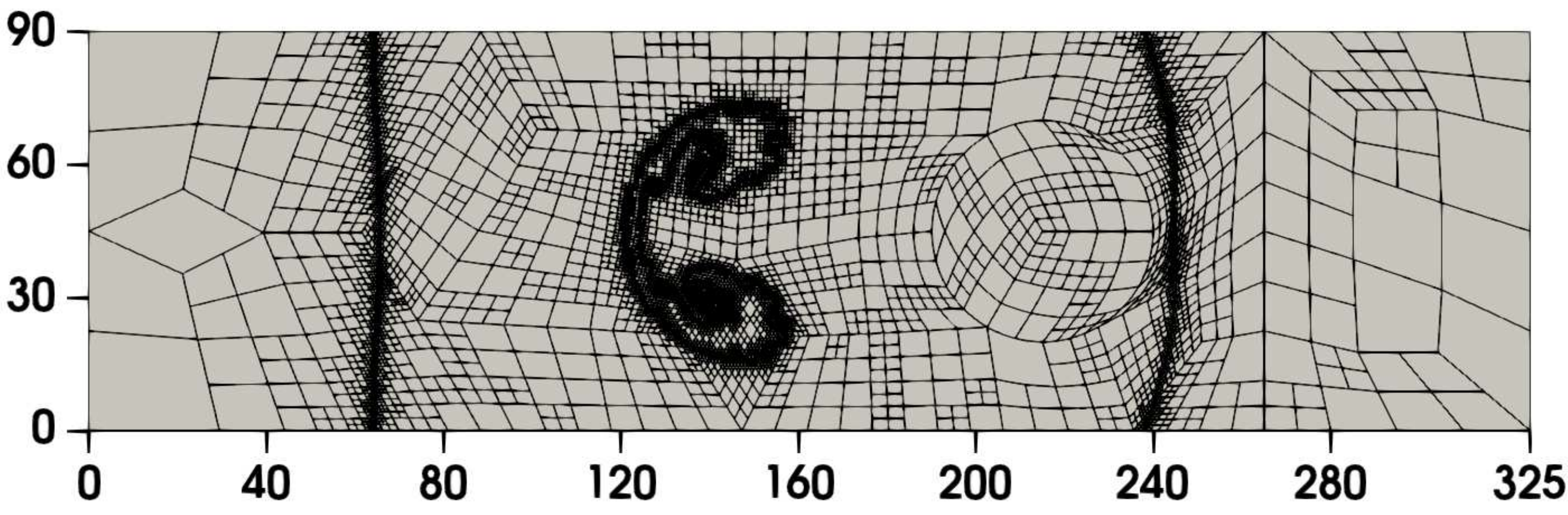}
    \caption{Adaptively refined mesh at $t=450$ for Case-I.}
    \label{fig: bs1_mesh_amr}
    \end{subfigure}
        \begin{subfigure}{0.49\textwidth}
    \includegraphics[width=\linewidth]{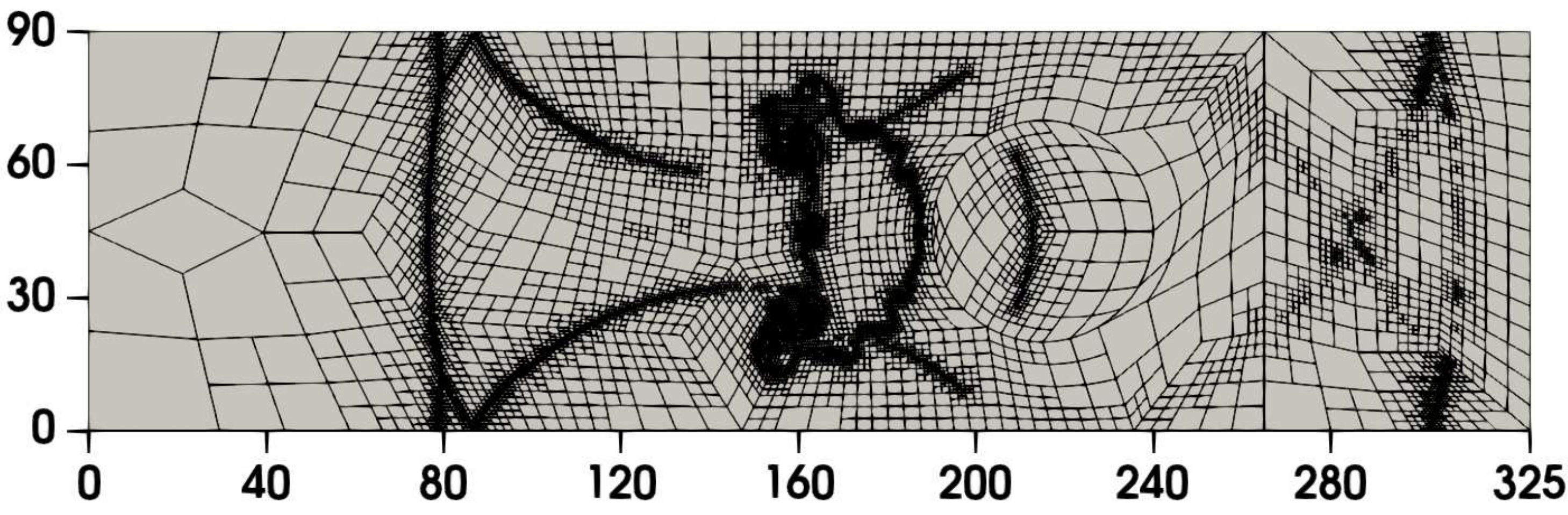}
    \caption{Adaptively refined mesh at $t=450$ for Case-II.}
    \label{fig: bs2_mesh_amr}
    \end{subfigure}
        \begin{subfigure}{0.49\textwidth}
    \includegraphics[width=\linewidth]{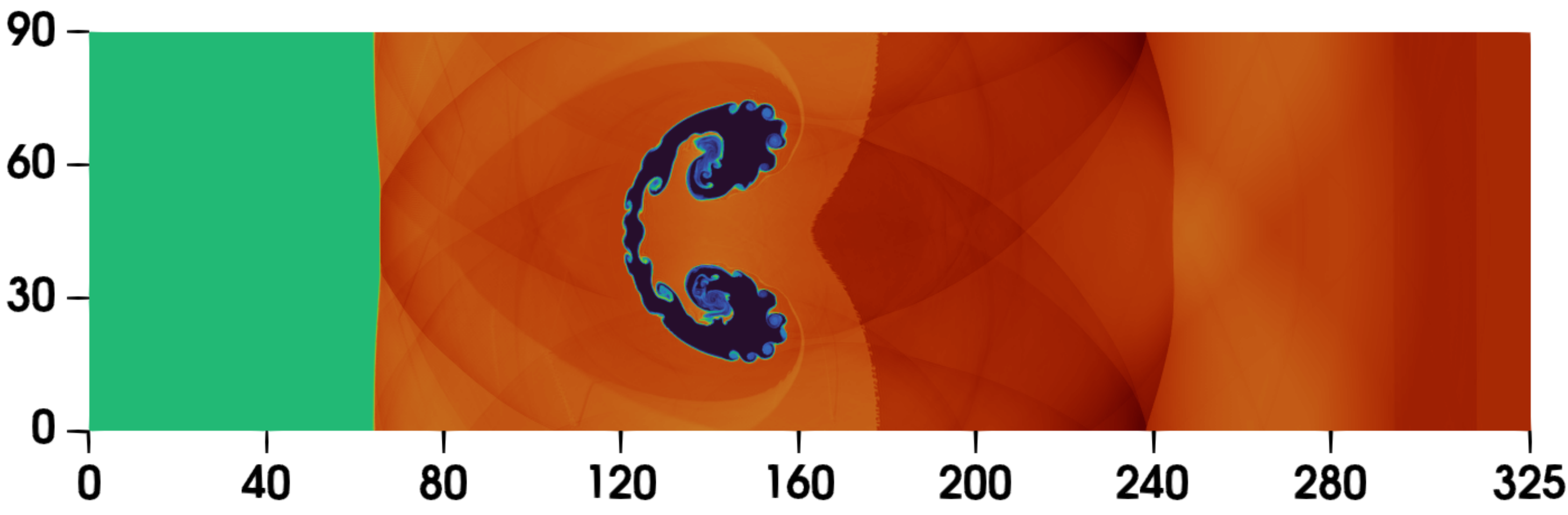}
    \caption{$\rho$ at $t=450$ without AMR for Case-I.}
    \label{fig: bs1_rho_tf_uniform}
    \end{subfigure}
        \begin{subfigure}{0.49\textwidth}
    \includegraphics[width=\linewidth]{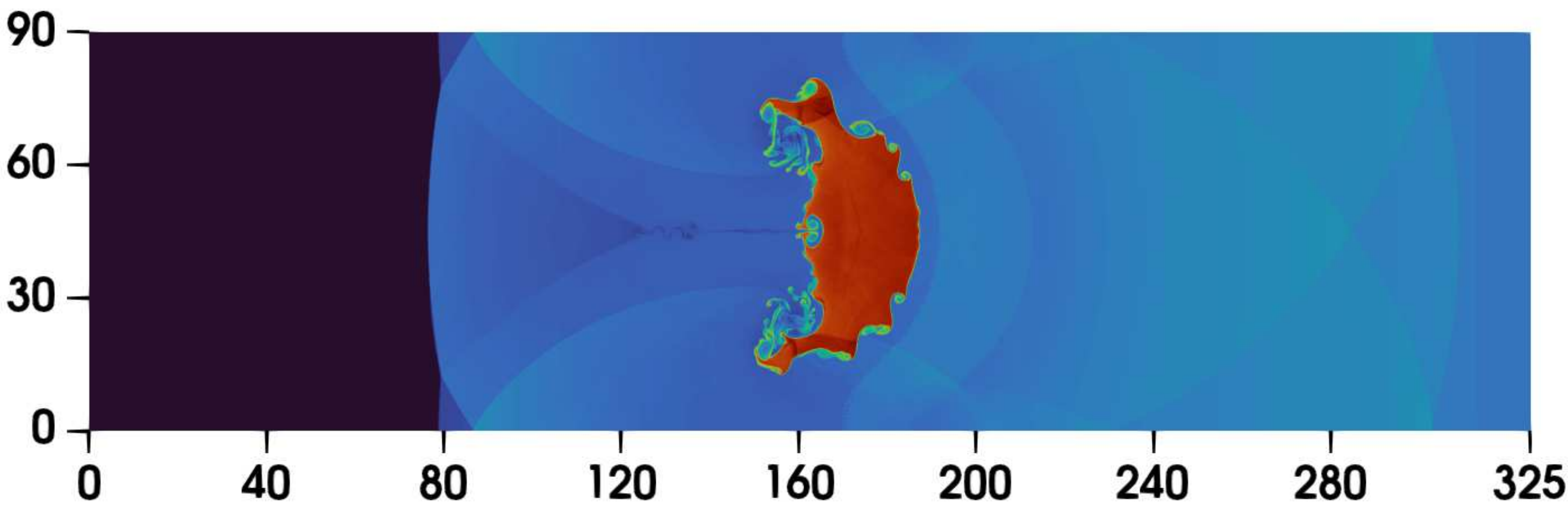}
    \caption{$\rho$ at $t=450$ without AMR for Case-II.}
    \label{fig: bs2_rho_tf_uniform}
    \end{subfigure}
    \caption{Bubble shock test: Density $\rho$ profiles and corresponding meshes.}
    \label{fig: bs}
\end{figure}

This test is used in~\cite{xu2024high} with RC-EOS~\eqref{eq: RC_eos}, and the result with IP-EOS~\eqref{eq: IP_eos} can be found in~\cite{basak2025constraints}. A bubble interacts with a moving shock and forms various waves in the domain. Initially, a shock is placed at $x=265$, and a circular bubble having radius $25$ with different density than its ambiance is placed with center at $(215, 45)$. The domain of computation is taken as $[0, 325] \times [0, 90]$ with $y=0, 90$ as reflective boundaries. The boundaries $x=0, 325$ are taken as artificial boundaries (Section~\ref{sec: artificial_boundary}) with components of $\mb{u}^{\text{set}}$ as in~\eqref{eq: bs_initial_data}. The pre- and post-shock states are defined as
\begin{align}\label{eq: bs_initial_data}
    (\rho&, v_1, v_2, p)= \begin{cases}
        (1, 0, 0, 0.05) & \text{if}\ x < 265\\
        (1.941272902134272, -0.200661045980881, 0,  0.15) & \text{if}\ x > 265.
    \end{cases}
\end{align}
The density inside the bubble is taken by dividing it into two cases,
\[
\text{case~I: } \rho = 0.1358, \quad \text{case~II: } \rho = 3.1538,
\]
with pressure $p$ same as the ambient pressure. We do the simulation \revb{up to} $t=450$ with $N=4$ using IP-EOS~\eqref{eq: IP_eos}. For the AMR case, the mesh is allowed to be refined and coarsened at each time step with the AMR indicator~\eqref{eq: amr_indicator} using the three-level controller~\eqref{eq: amr_controller} with
\begin{equation}\label{eq: case_bs_amr_params}
        (\texttt{base\_level}, \quad \texttt{med\_level}, \quad \texttt{max\_level}) = (0, 3, 6).
\end{equation}
The thresholds in the AMR controller are set as
\begin{equation*}
    (\epsilon_1, \epsilon_2) = (0.03, 0.1) \text{ for case~I, and } (\epsilon_1, \epsilon_2) = (0.02, 0.09) \text{ for case~II.}
\end{equation*}
For $\texttt{base\_level}=0$, the mesh, created with Gmsh~\cite{geuzaine2009gmsh}, is shown in Figure~\ref{fig: bs_mesh_bsl0}.

The density profiles with AMR are shown in Figure~\ref{fig: bs1_rho_t0}-Figure~\ref{fig: bs2_rho_tf_amr} along with the final meshes in Figures~\ref{fig: bs1_mesh_amr},~\ref{fig: bs2_mesh_amr}. For comparison purposes, we have also shown the results without AMR, using a mesh corresponding to the $\texttt{max\_level}=6$ in the whole domain in Figures~\ref{fig: bs1_rho_tf_uniform},~\ref{fig: bs2_rho_tf_uniform}. We observe that the results with AMR and without AMR are in good agreement. It is clear from the figure that the AMR indicator is working effectively, refining the mesh in some regions to capture the waves in the solution while coarsening it near the smooth regions. The number of elements in the final mesh with AMR is $10585$ and $13183$, for case~I and case~II, respectively; while without AMR the mesh has $237568$ elements throughout the whole simulation. We have also noted the wall-clock time for all the simulations:
\begin{align*}
    &\text{case~I with AMR: } \revb{4525} \text{ seconds}, \qquad \text{case~I without AMR: } 66729 \text{ seconds},\\
    &\text{case~II with AMR: } 4233 \text{ seconds}, \qquad \text{case~II without AMR : } 55669 \text{ seconds}.
\end{align*}
The simulations with AMR take significantly less wall-clock time compared to the simulations without AMR for both cases.

\subsection{Relativistic jet test}
The relativistic jet test is used in~\cite{wu2016physical} with RC-EOS~\eqref{eq: RC_eos} to show the robustness of the scheme. In this test, a jet with very high-speed, moving close to the speed of light enters the domain. Along with the high-speed jet, the solution also has relativistic shock waves, shear waves, ultra-relativistic regions, and interface instabilities~\cite{marti1994analytical, duncan1994simulations, komissarov1998large}. We take the domain of computation as $[-16, 16] \times [0, \reva{34}]$ with all the boundaries as artificial boundaries (Section~\ref{sec: artificial_boundary}) except the bottom boundary ($x=-16$).

Initially, the domain is filled with a fluid of unit density $\rho$. At the bottom boundary, a beam of density $\rho = 0.01$ is allowed to enter the domain with velocity $v_2 = \reva{0.9999999}$ through $\{(x,y):|x| < 0.5, y = 0 \}$ and with pressure \reva{calculated assuming the classical Mach number to be $1.74$. The ambient pressure throughout the domain is set equal to the pressure at the beam.} The rest of the bottom boundary is taken as mixed flow (Section~\ref{sec: mixed_flow}). For the artificial boundaries, the solution at the boundary is taken as $\rho = 1$, $(v_1, v_2) = (0,0)$, and $p$ \reva{same as the ambient pressure in the initial condition.}

\begin{figure}[]
    \centering
    \begin{subfigure}{0.31\textwidth}
    \includegraphics[width=\linewidth]{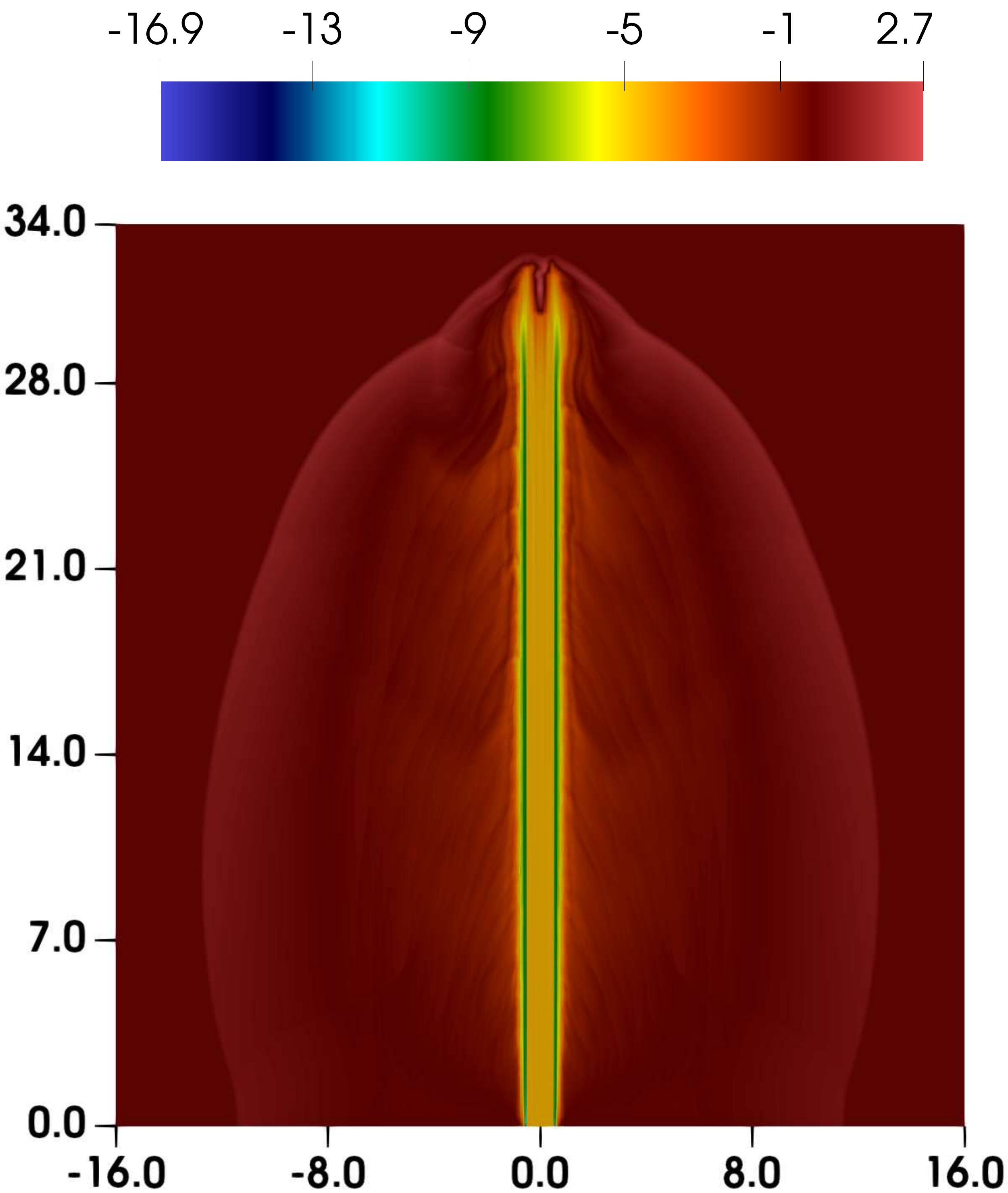}
    \caption{$\ln \rho$ with AMR.}
    \label{fig:reljet_rho_adaptive}
    \end{subfigure}
    \vspace{0.1cm}
    \quad
    \begin{subfigure}{0.31\textwidth}
    \includegraphics[width=\linewidth]{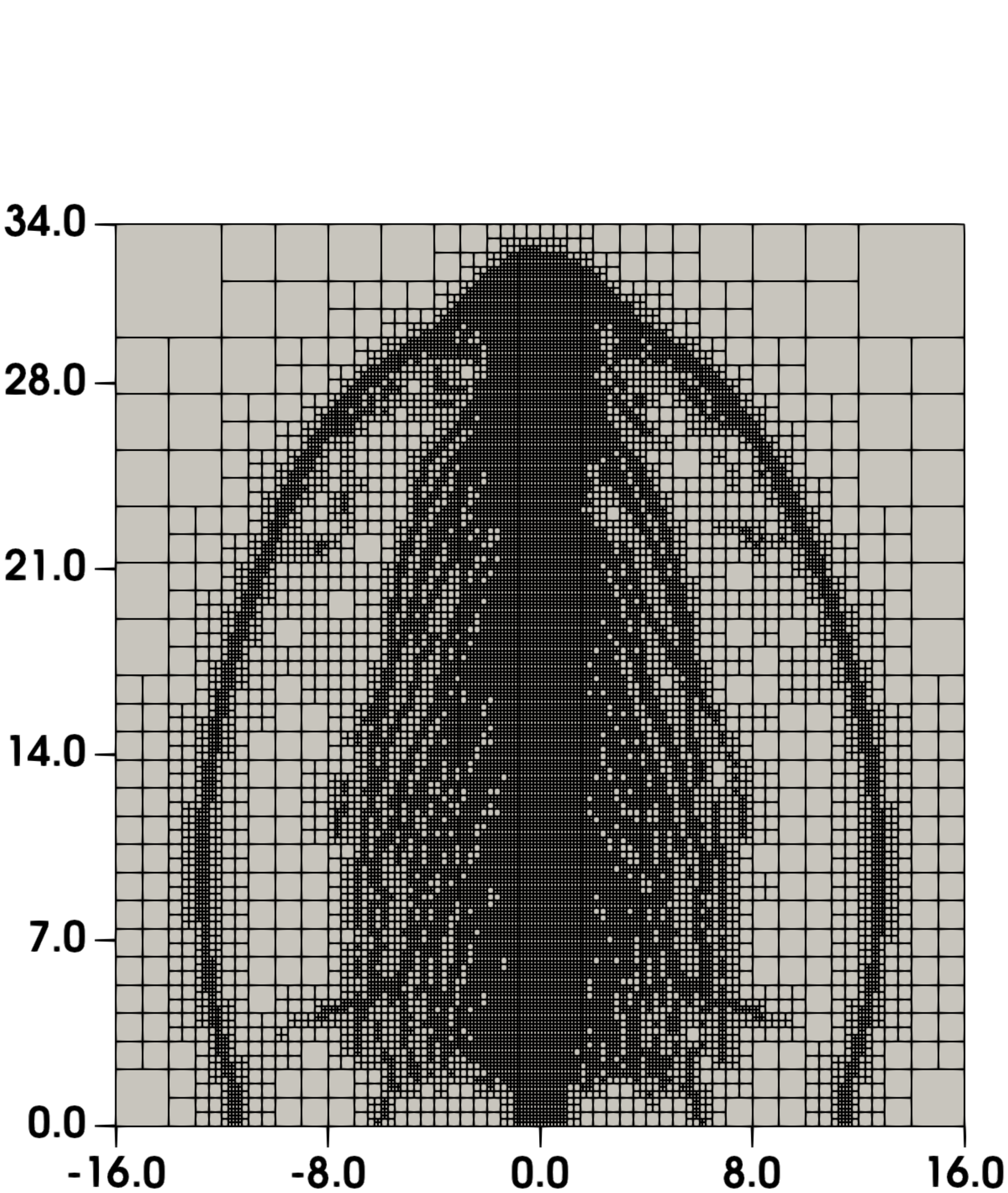}
    \caption{Adaptively refined mesh.}
    \label{fig:reljet_mesh_adaptive}
    \end{subfigure}
    \vspace{0.1cm}
    \quad
    \begin{subfigure}{0.31\textwidth}
    \includegraphics[width=\linewidth]{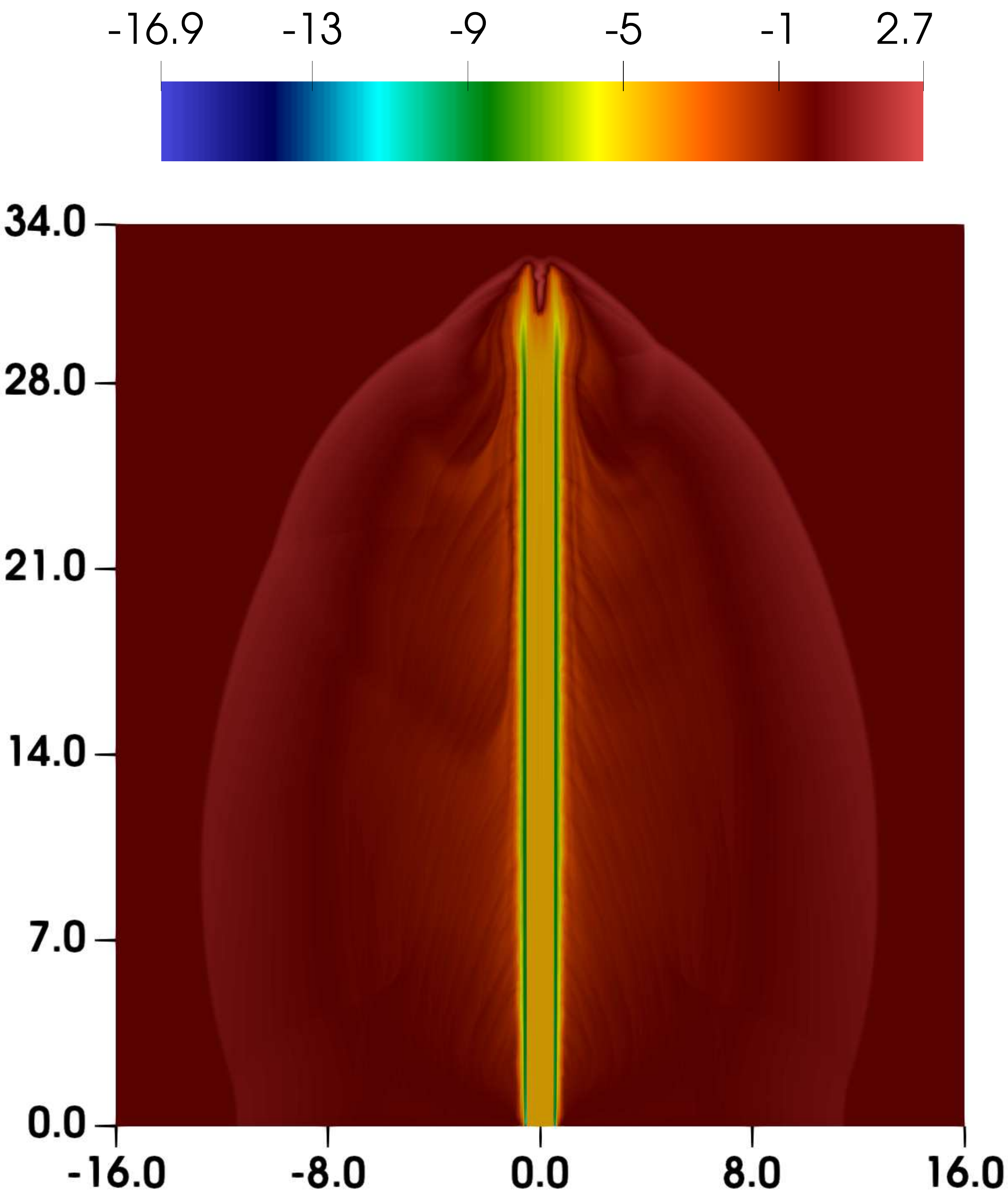}
    \caption{$\ln \rho$ with uniform mesh.}
    \label{fig:reljet_rho_uniform}
    \end{subfigure}
    \caption{Relativistic jet test: Results at time $t=30$.}
    \label{fig:reljet}
\end{figure}

The simulation is \revb{performed} with $N=4$, and RC-EOS~\eqref{eq: RC_eos} with both AMR and uniform mesh \revb{up to} $t=30$. \reva{At the final time, the maximum value of the Lorentz factor reaches approximately $2278$ with uniform mesh and $2394$ with AMR.} The density profile with AMR is shown in Figure~\ref{fig:reljet_rho_adaptive}, and the corresponding adaptively refined mesh at the final time in Figure~\ref{fig:reljet_mesh_adaptive}. The mesh is allowed to be refined and coarsened at each time step with the AMR indicator~\eqref{eq: amr_indicator} implemented with the three-level controller~\eqref{eq: amr_controller} with
\begin{equation}\label{eq: case_reljet_amr_params}
        (\texttt{base\_level}, \quad \texttt{med\_level}, \quad \texttt{max\_level}) = (0, 4, 8), \qquad
        (\epsilon_1, \epsilon_2) = (0.05, 0.1).
\end{equation}
The $\texttt{base\_level} = 0$ in~\eqref{eq: case_reljet_amr_params} refers to the $1\times 1$ mesh and $\epsilon_1, \epsilon_2$ are the thresholds used in the AMR controller. The result with the uniform mesh corresponding to the highest level of refinement, that is $2^8 \times 2^8$ elements, is shown in Figure~\ref{fig:reljet_rho_uniform}. It is clear from the figure that the results with AMR and uniform mesh are in good agreement with each other in capturing the tip of the beam, along with the other structures. It is also clear from the figure that in AMR, the mesh is getting refined and coarsened according to flow features and is able to capture the waves effectively. The final mesh with AMR has $\reva{22186}$ elements while the uniform mesh has $65536$ elements. We also note that the wall-clock time for the simulation with the uniform mesh is $\reva{36168}$ seconds, while with AMR it is $\reva{9266}$ seconds, showing the benefit of using AMR.

\subsection{Riemann problem 1}
\begin{figure}[]
    \centering
    \begin{subfigure}{0.33\textwidth}
    \includegraphics[width=\linewidth]{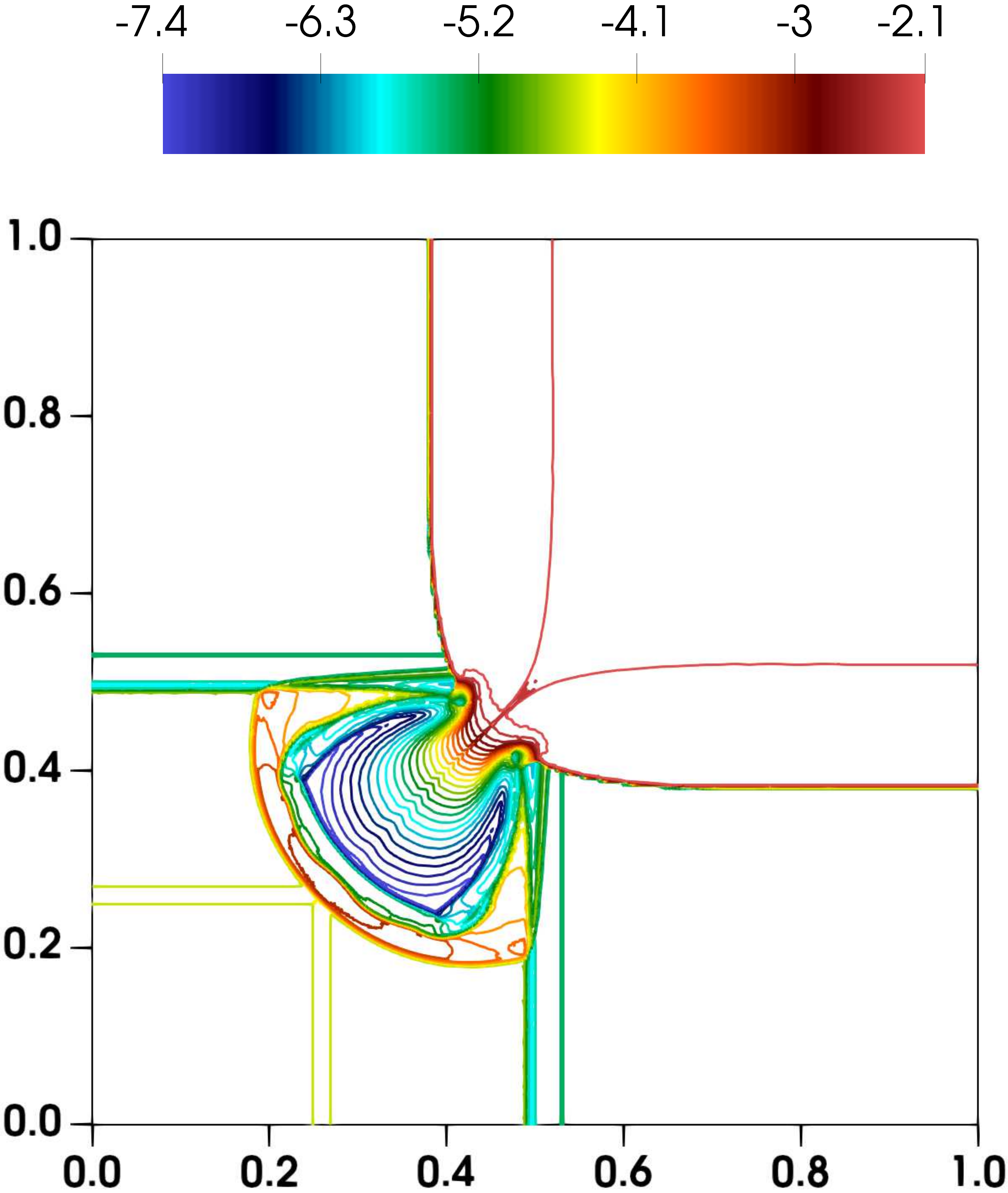}
    \caption{$\ln\rho$ with AMR.}
    \label{fig: lnrho_wurp2_amr}
    \vspace*{10pt}
    \end{subfigure}
    \begin{subfigure}{0.33\textwidth}
    \includegraphics[width=\linewidth]{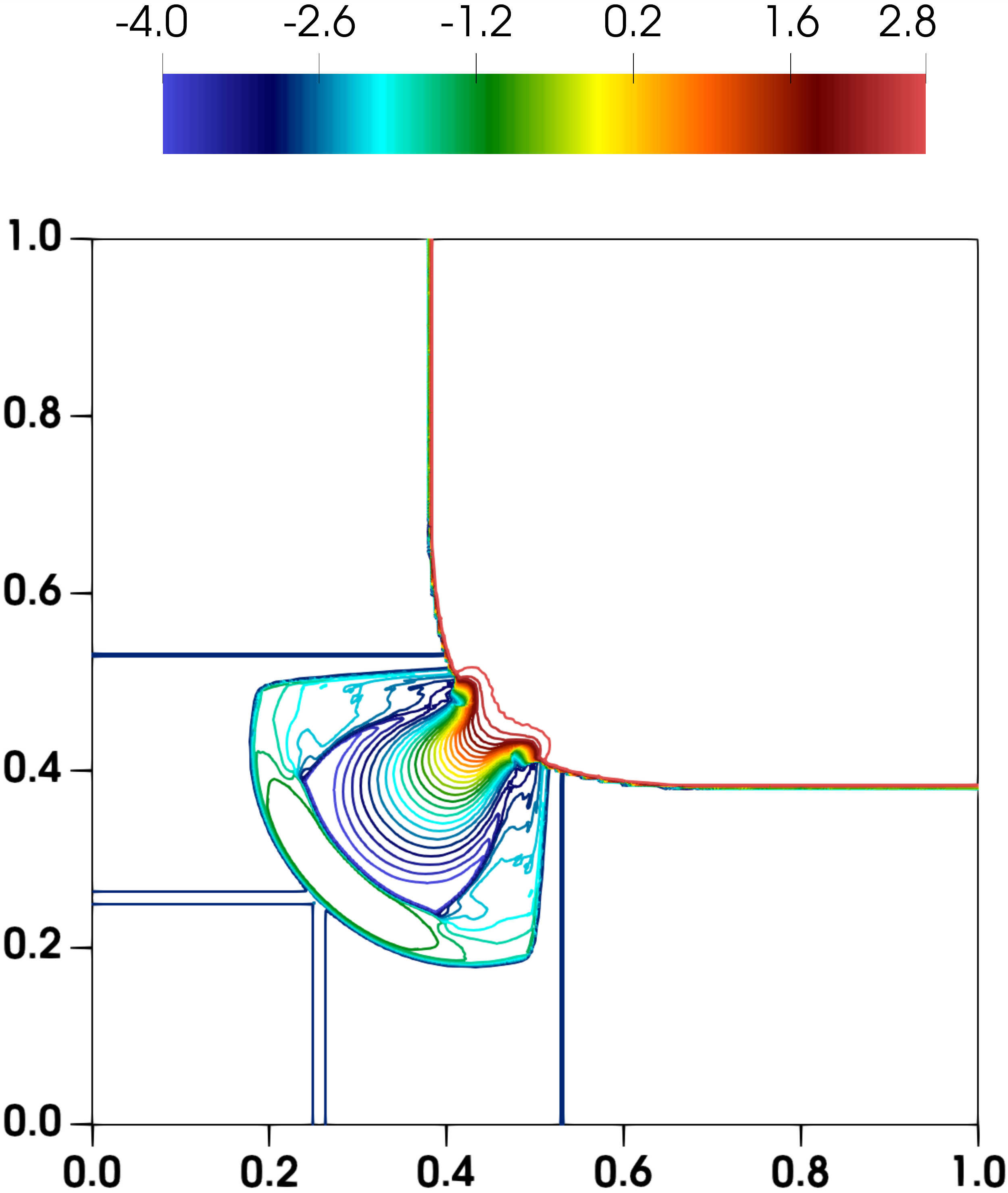}
    \caption{$\ln p$ with AMR.}
    \label{fig: lnp_wurp2_amr}
    \vspace*{10pt}
    \end{subfigure}
    \begin{subfigure}{0.33\textwidth}
    \includegraphics[width=\linewidth]{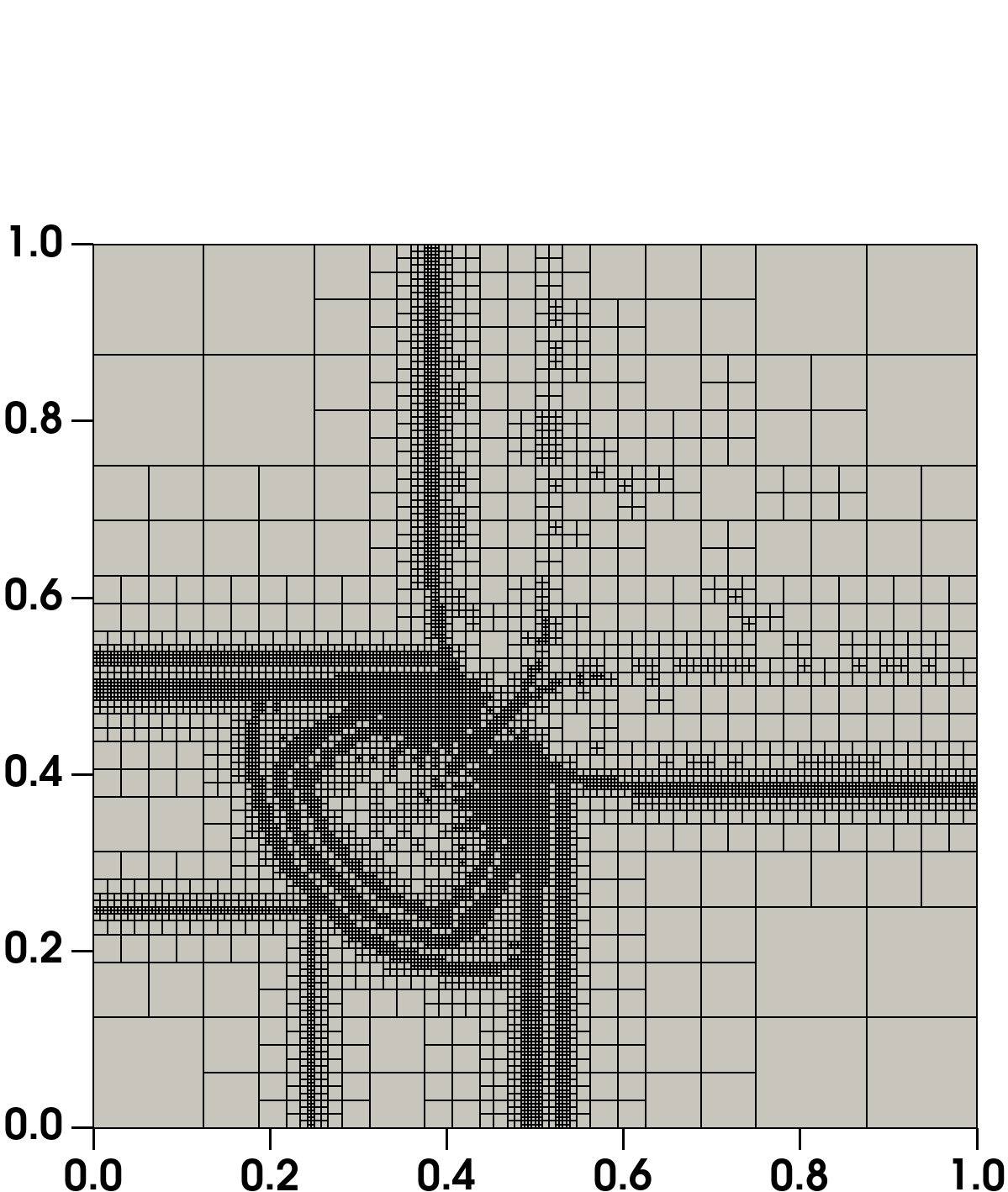}
    \caption{Adaptively refined mesh.}
    \label{fig: mesh_wurp2_amr}
    \vspace*{10pt}
    \end{subfigure}
        \begin{subfigure}{0.33\textwidth}
    \includegraphics[width=\linewidth]{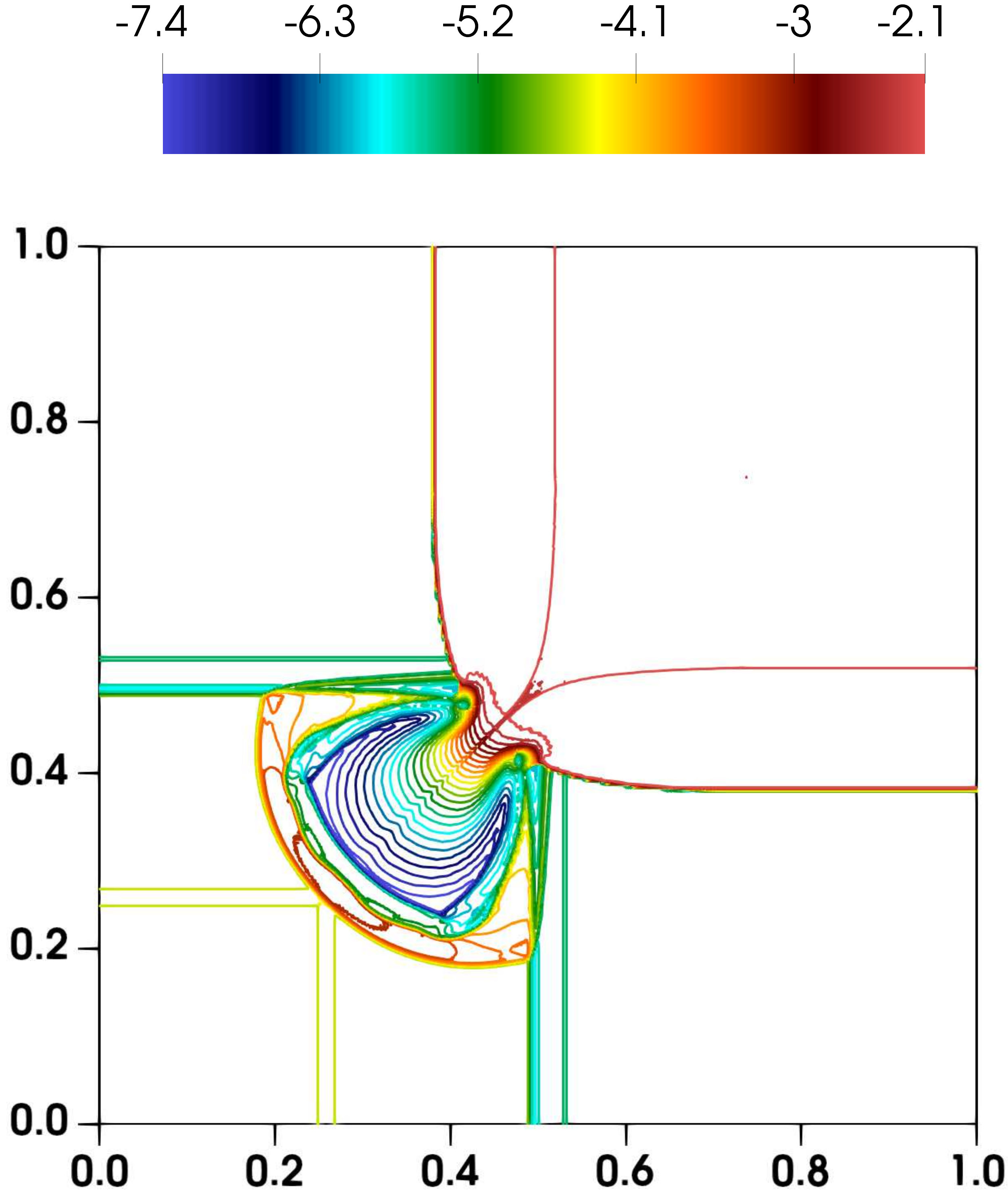}
    \caption{$\ln\rho$ with uniform mesh.}
    \label{fig: lnrho_wurp2_without_amr}
    \end{subfigure}
           \begin{subfigure}{0.33\textwidth}
    \includegraphics[width=\linewidth]{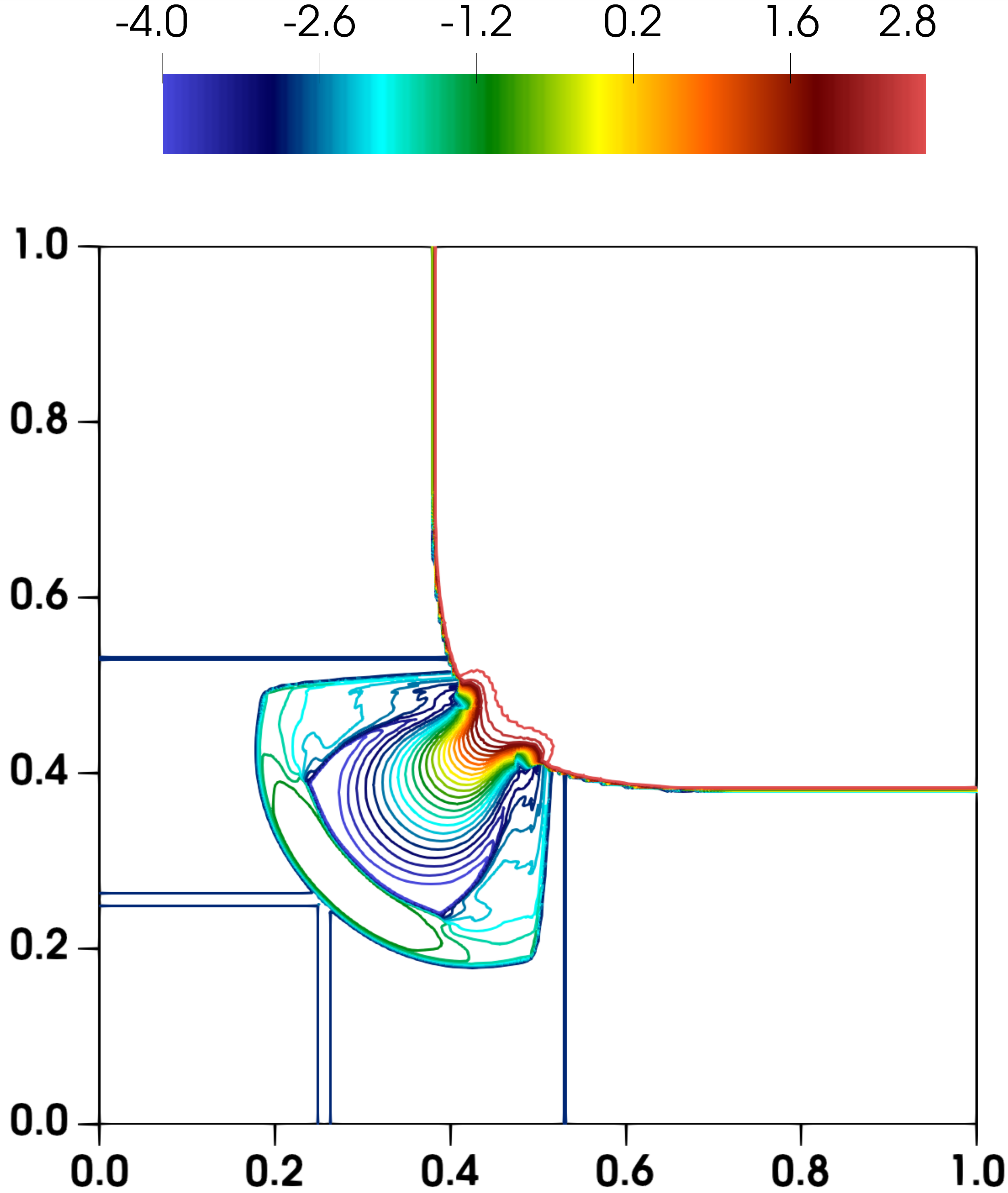}
    \caption{$\ln p$ with uniform mesh.}
    \label{fig: lnp_wurp2_without_amr}
    \end{subfigure}
    \caption{Riemann problem 1: Results at time $t=0.4$ in $[0,1] \times [0,1]$. The contour plots have 25 uniform contours in the specified regions in corresponding color bars.}
    \label{fig: wurp2}
\end{figure}

This Riemann problem is taken from~\cite{wu2015high}. At the initial time, four constant states in the four quadrants of $[0,1]\times[0,1]$ are given by
\begin{align*}
    (\rho, v_1, v_2, p) = \begin{cases}
        (0.1, 0, 0, 20) & \text{if}\ x > 0.5,\ y > 0.5\\
        (0.00414329639576, 0.9946418833556542, 0, 0.05) & \text{if}\ x < 0.5,\ y>0.5\\
        (0.01, 0, 0, 0.05) & \text{if}\ x < 0.5,\ y < 0.5\\
        (0.00414329639576, 0,0.9946418833556542, 0.05) & \text{if}\ x > 0.5,\ y<0.5.\\
    \end{cases}
\end{align*}
The simulation is \revb{carried out} with $N=4$ and IP-EOS~\eqref{eq: IP_eos} \revb{up to} $t=0.4$. For this simulation, we consider an extended domain beyond $[0,1]\times [0,1]$ as we do not have the information at the boundary from the outside, which we need for the artificial boundaries (Section~\ref{sec: artificial_boundary}). This idea is also used in the literature, for example, in Figure 2 of~\cite{chan2025artificial}. The solution profiles $\ln \rho$, $\ln p$ with AMR are shown in Figure~\ref{fig: lnrho_wurp2_amr},\ref{fig: lnp_wurp2_amr} along with the final mesh in Figure~\ref{fig: mesh_wurp2_amr}. Here as well, we have used the AMR indicator~\eqref{eq: amr_indicator}, that is implemented with the three-level controller~\eqref{eq: amr_controller} with the parameters
\begin{equation}\label{eq: case_rp1_amr_params}
        (\texttt{base\_level}, \quad \texttt{med\_level}, \quad \texttt{max\_level}) = (0, 3, 8), \qquad
        (\epsilon_1, \epsilon_2) = (0.07, 0.08).
\end{equation}
\revb{The} $\epsilon_1, \epsilon_2$ are the thresholds used in~\eqref{eq: amr_controller}, and $\texttt{base\_level}=0$ is equivalent to a mesh with $1\times 1$ element in $[0,1] \times [0,1]$, which is allowed to be refined and coarsened in each time step. We also present the solution profiles in Figure~\ref{fig: lnrho_wurp2_without_amr},~\ref{fig: lnp_wurp2_without_amr} with a uniform mesh having the same resolution corresponding to $\texttt{max\_level} = 8$ in the adaptive mesh. It is observed from the figures that the contours are in good agreement with each other in the results with AMR and the uniform mesh. The solution profile initially has a pair of contact discontinuities and a pair of shock waves, which collide to form mushroom-shaped waves in the bottom-left quadrant. From the figure~\ref{fig: mesh_wurp2_amr}, it is clear that the AMR is working effectively, capturing all the structures in the solution, while coarsening the mesh in the smoother regions. We also note that, inside the domain $[0,1]\times [0,1]$, the number of elements in the final mesh with AMR is $9304$ while without AMR it is $65536$. Consequently, the computational cost decreases significantly with AMR.

\subsection{Riemann problem 2}
\begin{figure}[]
    \centering
    \begin{subfigure}{0.33\textwidth}
    \includegraphics[width=\linewidth]{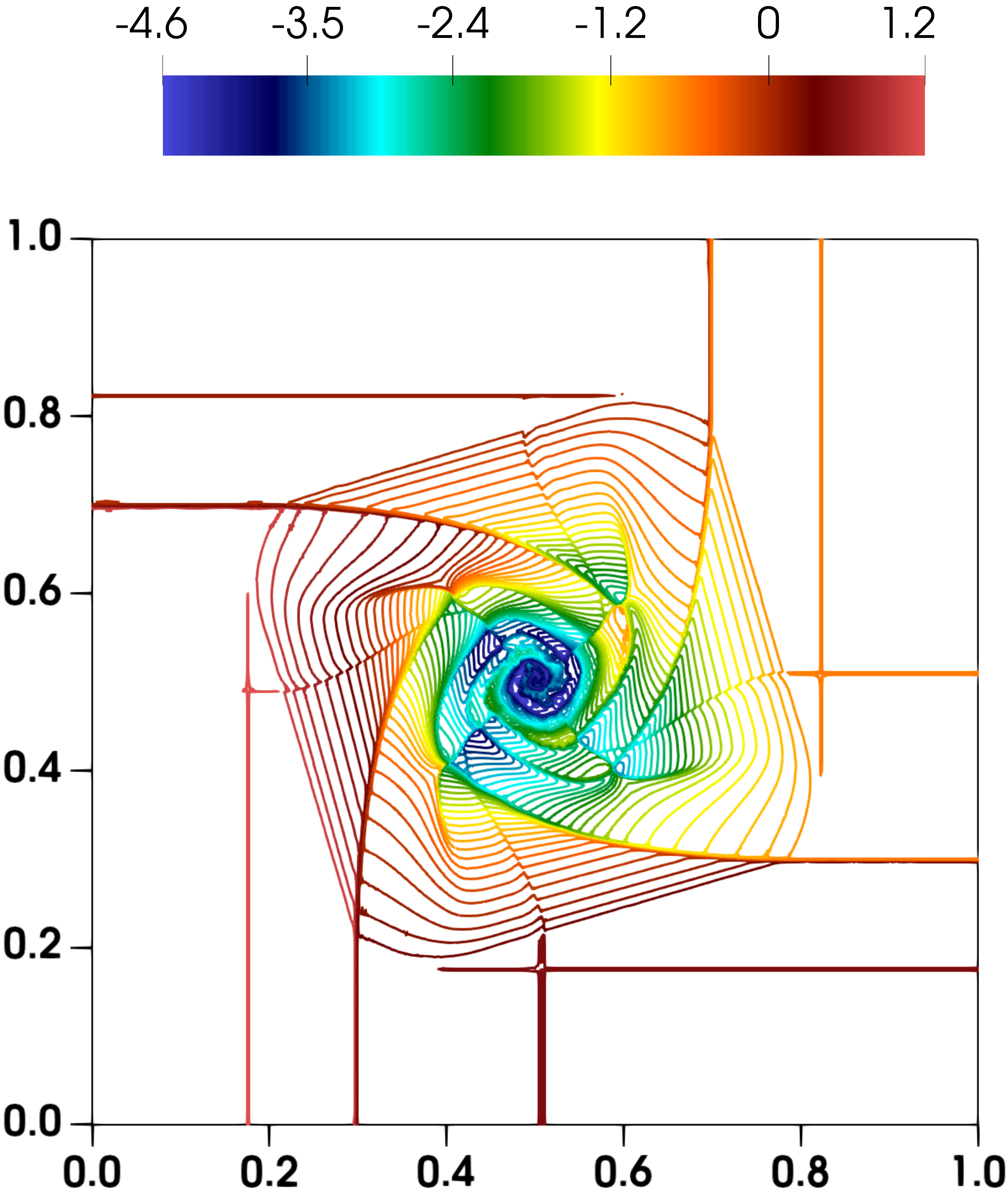}
    \caption{$\ln\rho$ with AMR.}
    \label{fig: lnrho_wu2rp1_amr}
    \vspace*{5pt}
    \end{subfigure}
    \begin{subfigure}{0.33\textwidth}
    \includegraphics[width=\linewidth]{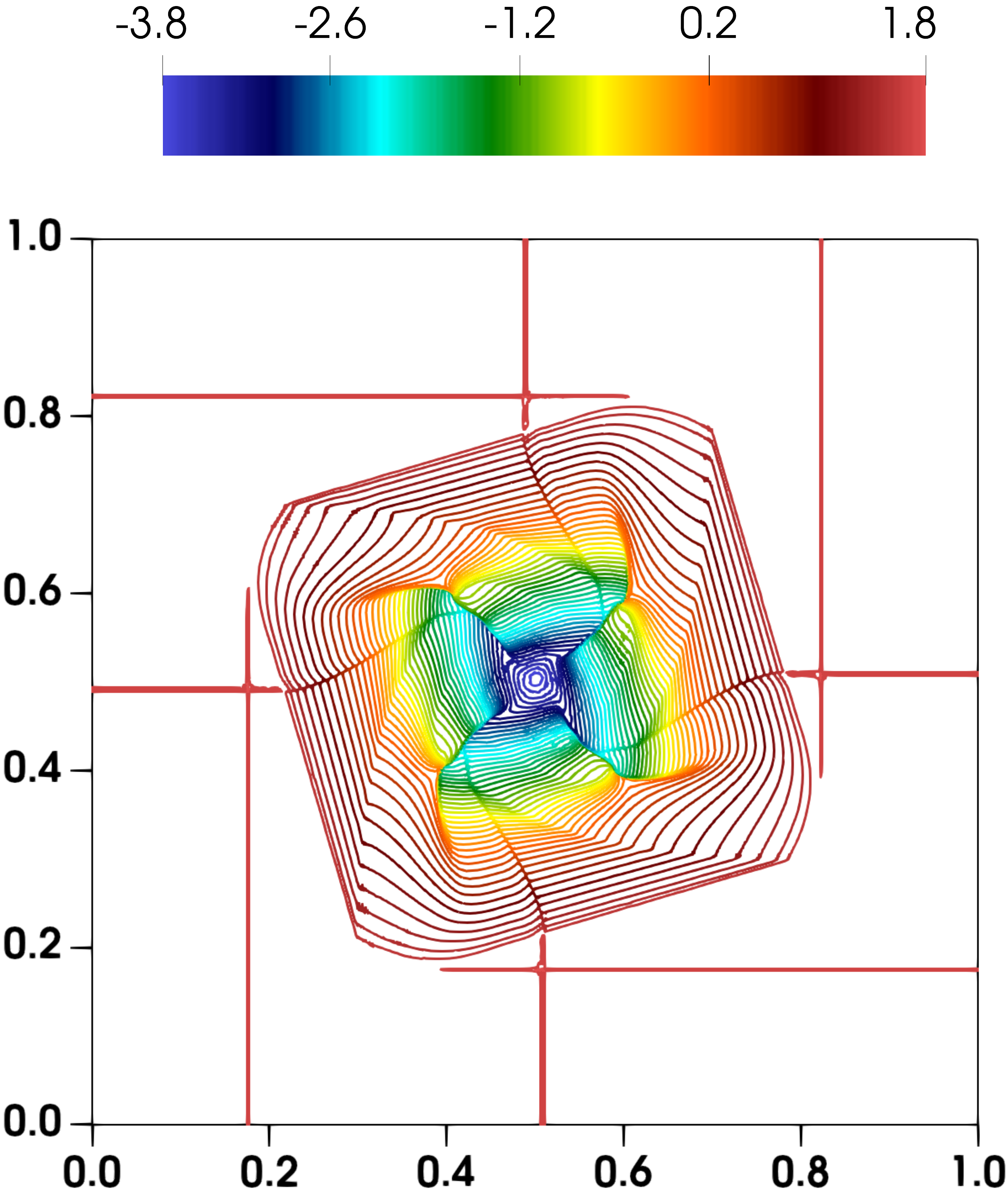}
    \caption{$\ln p$ with AMR.}
    \label{fig: lnp_wu2rp1_amr}
    \vspace*{5pt}
    \end{subfigure}
    \begin{subfigure}{0.33\textwidth}
    \includegraphics[width=\linewidth]{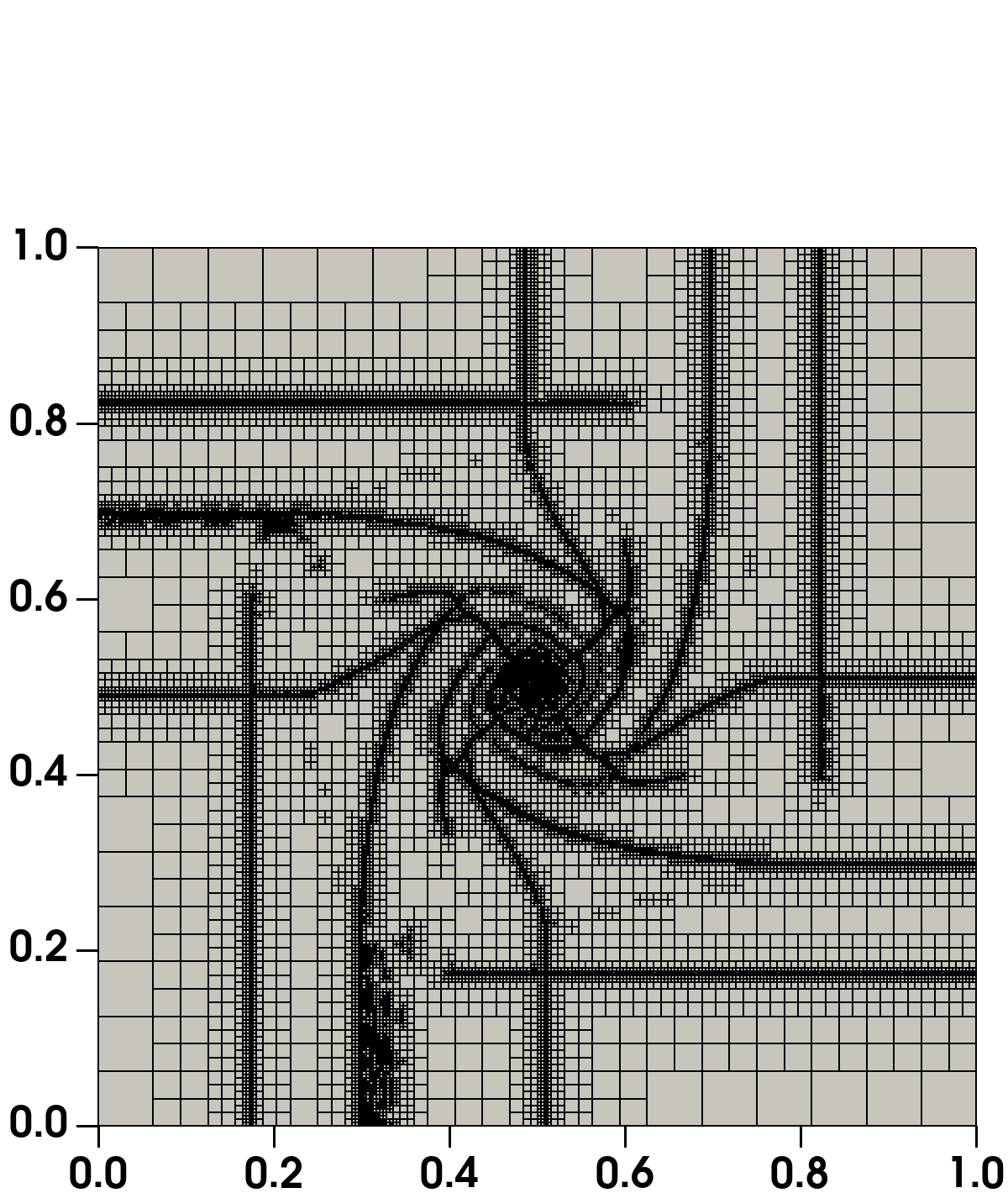}
    \caption{Adaptively refined mesh.}
    \label{fig: mesh_wu2rp1_amr}
    \vspace*{5pt}
    \end{subfigure}
        \begin{subfigure}{0.33\textwidth}
    \includegraphics[width=\linewidth]{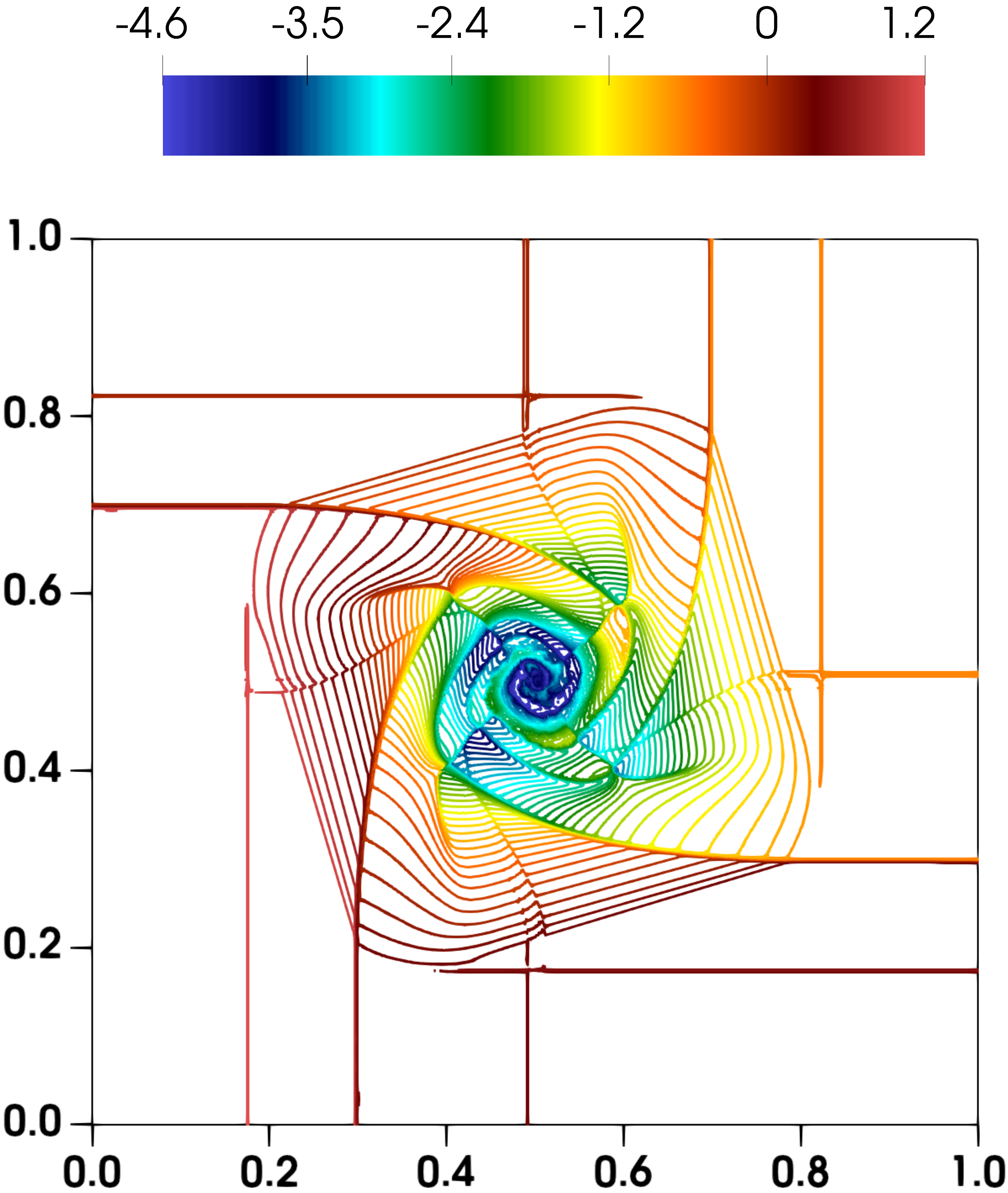}
    \caption{$\ln\rho$ with uniform mesh.}
    \label{fig: lnrho_wu2rp1_without_amr}
    \end{subfigure}
           \begin{subfigure}{0.33\textwidth}
    \includegraphics[width=\linewidth]{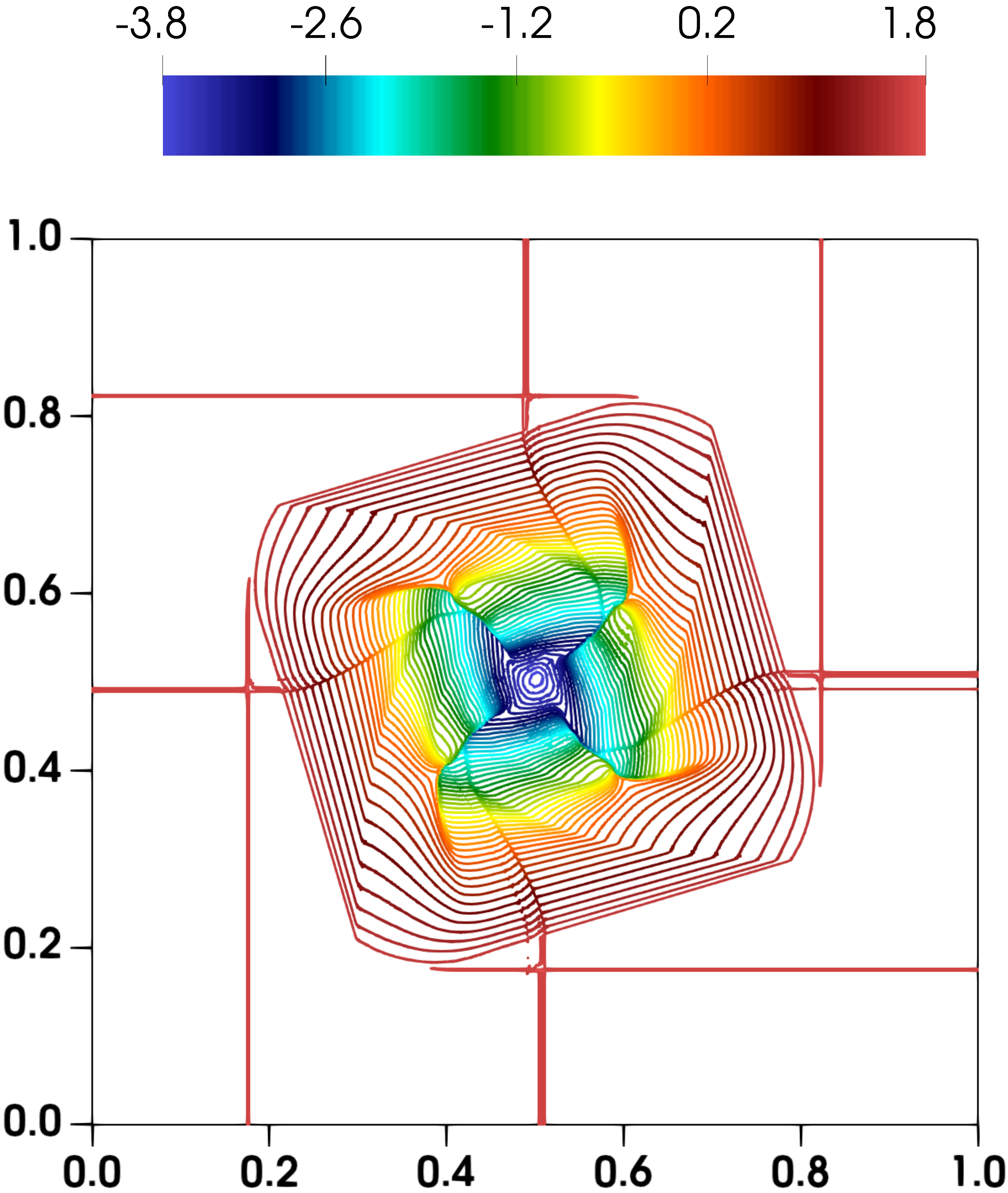}
    \caption{$\ln p$ with uniform mesh.}
    \label{fig: lnp_wu2rp1_without_amr}
    \end{subfigure}
    \caption{Riemann problem 2: Results at time $t=0.4$  in $[0,1] \times [0,1]$. The contour plots have 50 uniform contours in the specified regions in corresponding color bars.}
    \label{fig: wu2rp1}
\end{figure}

\begin{figure}[]
    \centering
    \begin{subfigure}{0.48\textwidth}
    \includegraphics[width=\linewidth]{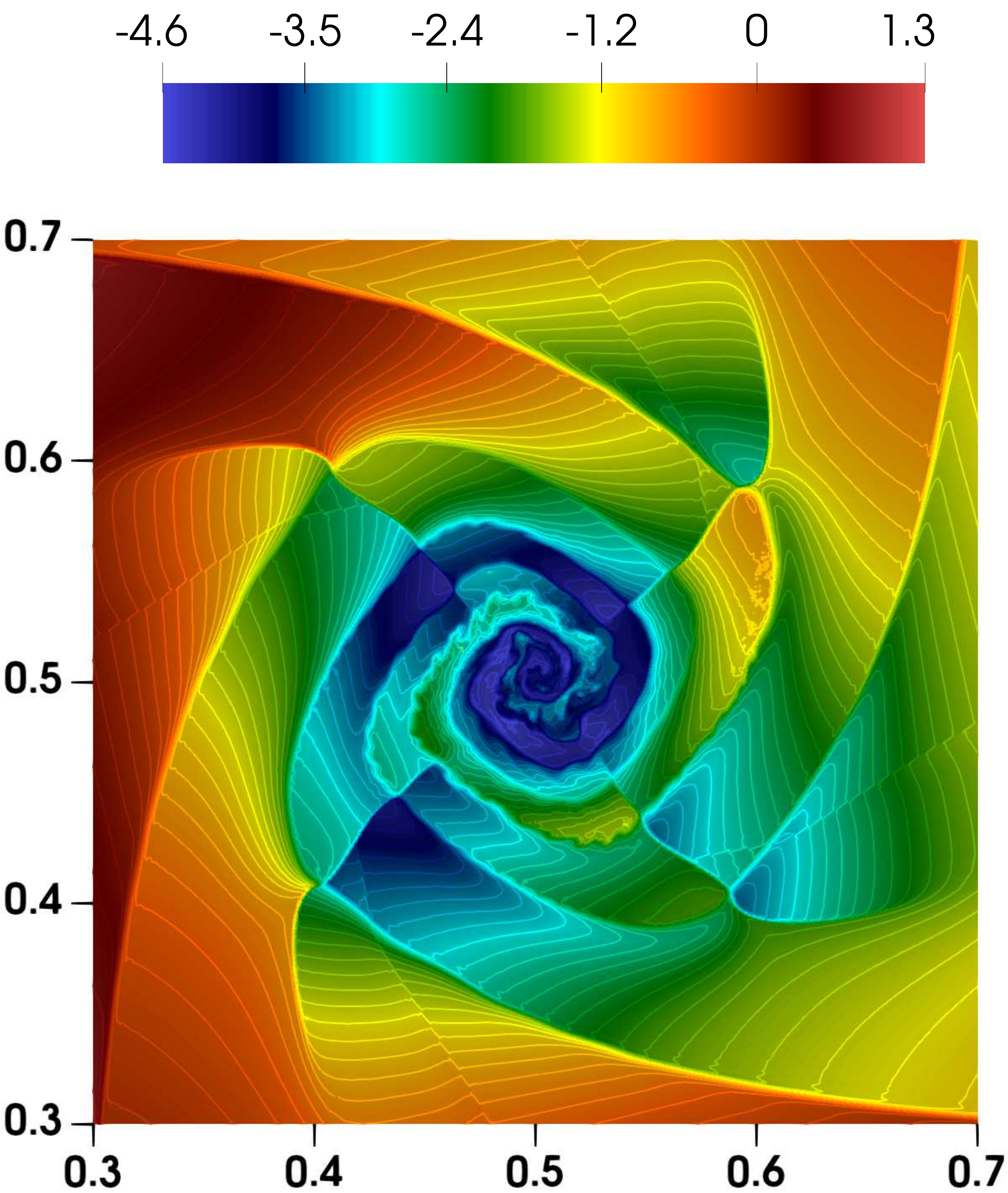}
    \caption{$\ln\rho$ with AMR.}
    \label{fig: lnrho_wu2rp1_with_amr_zoom}
    \end{subfigure}
           \begin{subfigure}{0.48\textwidth}
    \includegraphics[width=\linewidth]{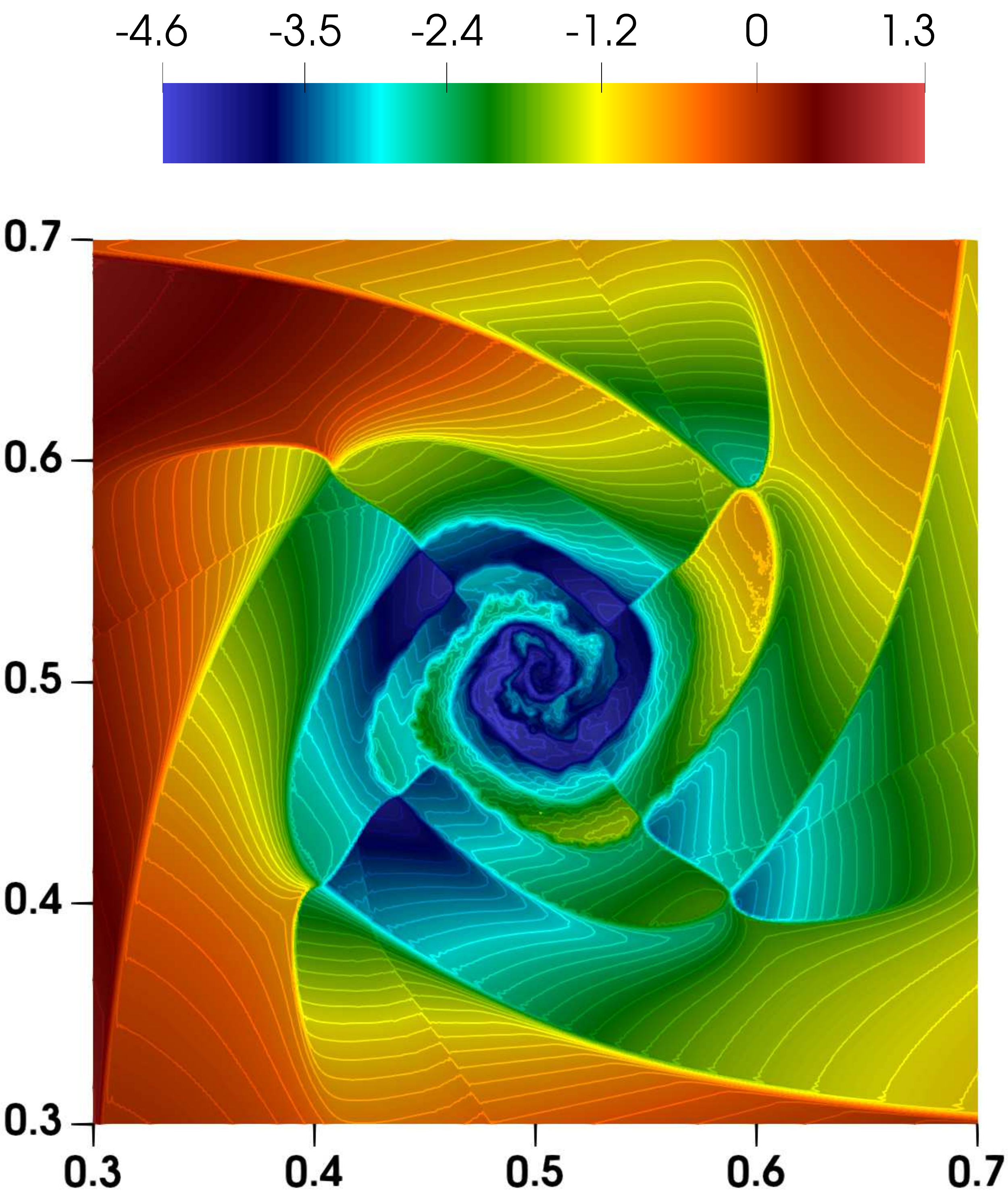}
    \caption{$\ln \rho$ with uniform mesh.}
    \label{fig: lnrho_wu2rp1_without_amr_zoom}
    \end{subfigure}
    \caption{Riemann problem 2: Results with $50$ contours in $[-4.6, 1.3]$ at time $t=0.4$ by zooming in $[0.3, 0.7] \times [0.3, 0.7]$.}
    \label{fig: wu2rp1_zoom}
\end{figure}

This Riemann problem can be found in~\cite{he2012adaptive, nunez2016xtroem} simulated with ID-EOS~\eqref{eq: ID_eos}, and in~\cite{basak2025constraints} simulated with TM-EOS~\eqref{eq: TM_eos}. At the initial time, the solution profile is given by
\[
(\rho, v_1, v_2, p) = \begin{cases}
        (0.5, 0.5, -0.5, 5) & \text{if}\ x > 0.5,\ y > 0.5\\
        (1, 0.5, 0.5, 5) & \text{if}\ x < 0.5,\ y>0.5\\
        (3, -0.5, 0.5, 5) & \text{if}\ x < 0.5,\ y < 0.5\\
        (1.5, -0.5, -0.5, 5) & \text{if}\ x > 0.5,\ y<0.5
    \end{cases}
\]
in the domain $[0,1]\times [0,1]$. Similar to the last Riemann problem, here as well we extend the domain beyond $[0,1] \times [0,1]$ for the simulation. The equation of state is taken as TM-EOS~\eqref{eq: TM_eos}, final time as $0.4$, and $N=4$. The time step is found as in~\cite{basak2025bound}. We consider the adaptive mesh with the AMR indicator~\eqref{eq: amr_indicator} along with the three-level controller~\eqref{eq: amr_controller} with
\begin{equation}\label{eq: case_rp2_amr_params}
        (\texttt{base\_level}, \quad \texttt{med\_level}, \quad \texttt{max\_level}) = (0, 3, 9), \qquad
        (\epsilon_1, \epsilon_2) = (0.07, 0.08).
\end{equation}
\revb{The quantities} $\epsilon_1, \epsilon_2$ are the thresholds used in the AMR controller, and $\texttt{base\_level} = 0$ is equivalent to a mesh with $1\times 1$ element in $[0,1]\times [0,1]$, that is allowed to be refined and coarsened in each time step. The results with the adaptive mesh are presented in Figures~\ref{fig: lnrho_wu2rp1_amr},~\ref{fig: lnp_wu2rp1_amr}, along with the final mesh in Figure~\ref{fig: mesh_wu2rp1_amr}. In Figures~\ref{fig: lnrho_wu2rp1_without_amr},~\ref{fig: lnp_wu2rp1_without_amr}, we present the result with a uniform mesh corresponding to $\texttt{max\_level} = 9$ for comparison purposes, and observe that the results are in good agreement with AMR results. The $\ln\rho$ profile has a spiral structure in the domain, with the center having the lowest value. We also present a zoomed view near the center of the domain for the $\ln\rho$ profile with colored contours in Figure~\ref{fig: lnrho_wu2rp1_with_amr_zoom},\ref{fig: lnrho_wu2rp1_without_amr_zoom} with AMR and uniform mesh, respectively. We observe that AMR works effectively for this problem as well, capturing all the sharp structures in the domain with a refined mesh while coarsening it in the smoother regions. The mesh near the center in the AMR case is refined to the $\texttt{max\_level}$ to capture the lowest pressure value. Both AMR and uniform mesh cases give the same lowest pressure in the center. It is also worth mentioning that in the domain $[0,1] \times [0,1]$, the final mesh with AMR has $26167$ elements while the mesh without AMR has $262144$ elements. Hence, there is a significant improvement in the computational cost when using AMR.

\subsection{Kelvin-Helmholtz instability test}
This test is taken from~\cite{beckwith2011second,radice2012thc,zanotti2015high}, which develops some fluid instabilities that are hard to capture with diffusive schemes or a low-resolution mesh. We take the computational domain as $[-1,1]\times [-0.5, 0.5]$ with periodic boundaries. At the initial time, the fluid state is given by dividing the domain into two parts:
\begin{align*}
    &\rho = \begin{cases}
        0.505 - 0.495 \tanh \left(\frac{x+0.5}{a}\right), \quad \text{if } x<0\\
        0.505 + 0.495 \tanh \left(\frac{x-0.5}{a}\right), \quad \text{if } x>0,
    \end{cases}\\
    &v_1 = \begin{cases}
        -\eta_0 v_s \sin(2\pi y) \exp\left(\frac{-(x+0.5)^2}{\sigma} \right), \quad \text{if } x<0\\
        \eta_0 v_s \sin(2\pi y) \exp\left(\frac{-(x-0.5)^2}{\sigma} \right), \quad \text{if } x>0,
    \end{cases}\\
    &v_2 = \begin{cases}
        -v_s \tanh\left(\frac{x+0.5}{a}\right), \quad \text{if } x<0\\
        v_s \tanh\left(\frac{x-0.5}{a}\right), \quad \text{if } x>0,
    \end{cases}\\
    &p = 1.0.
\end{align*}
The characteristic size $a$ is taken as $0.01$,  and $v_s = 0.5$. The other parameters are taken as $\eta_0 = 0.1$, $\sigma = 0.1$, which define the small perturbation in the velocity $v_1$.

\begin{figure}[]
    \centering
        \begin{subfigure}{0.6\textwidth}
    \includegraphics[width=\linewidth]{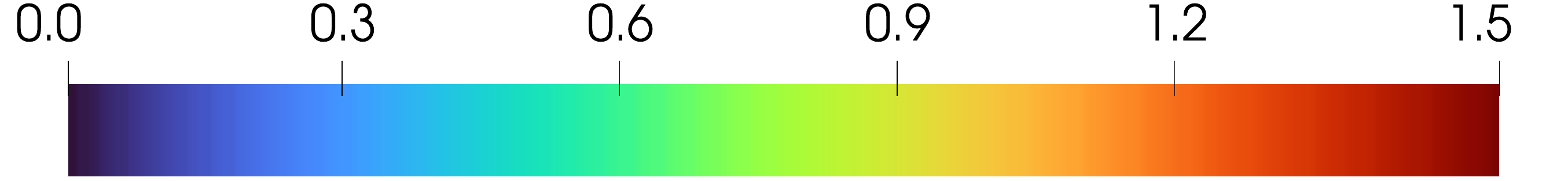}
    \vspace{0.01cm}
    \end{subfigure}
    \qquad
    \begin{subfigure}{0.49\textwidth}
    \includegraphics[width=\linewidth]{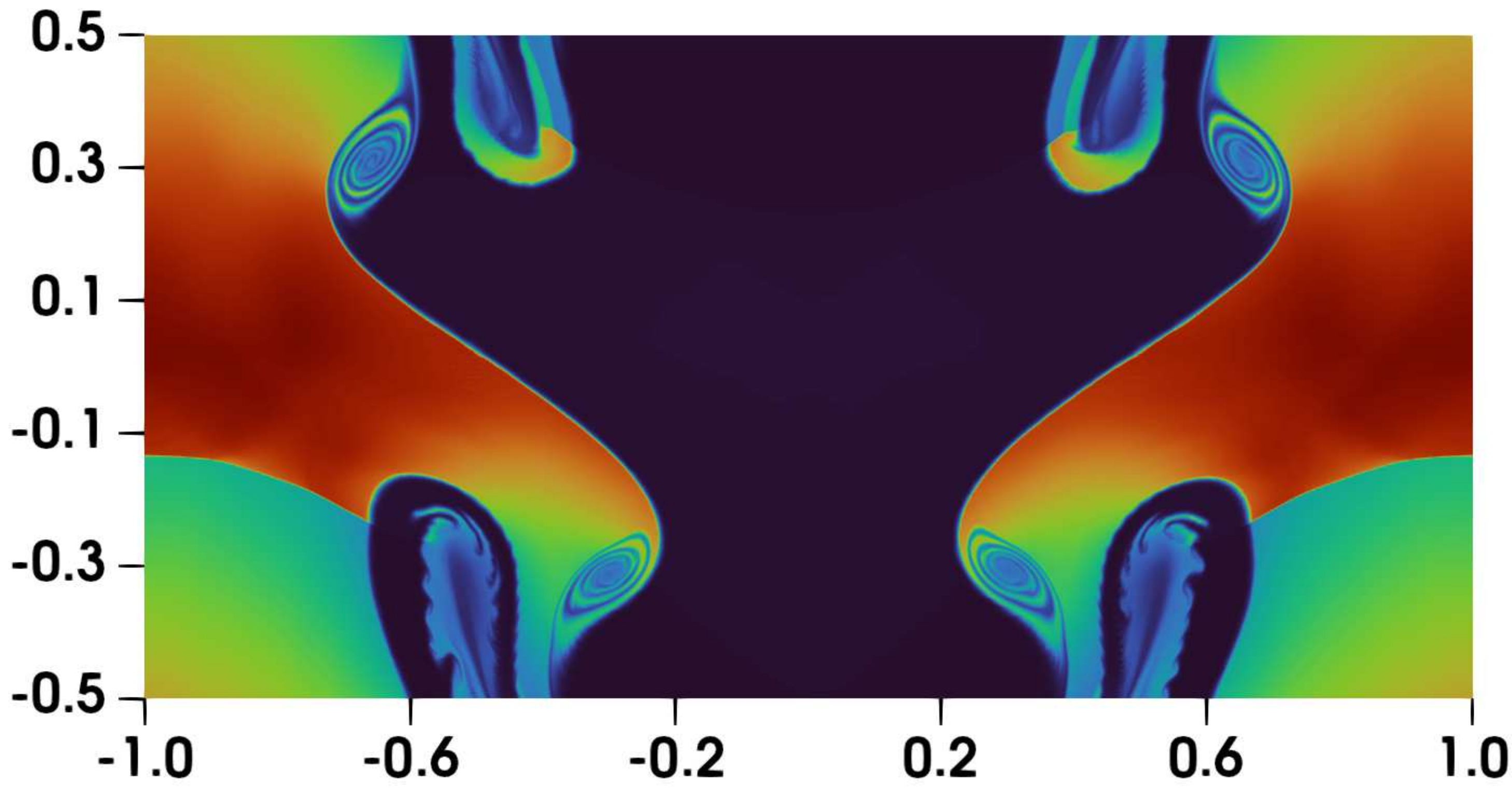}
    \caption{$\rho$ with AMR.}
    \label{fig:kh_rho_adaptive}
    \end{subfigure}
    \vspace{0.1cm}
    \begin{subfigure}{0.49\textwidth}
    \includegraphics[width=\linewidth]{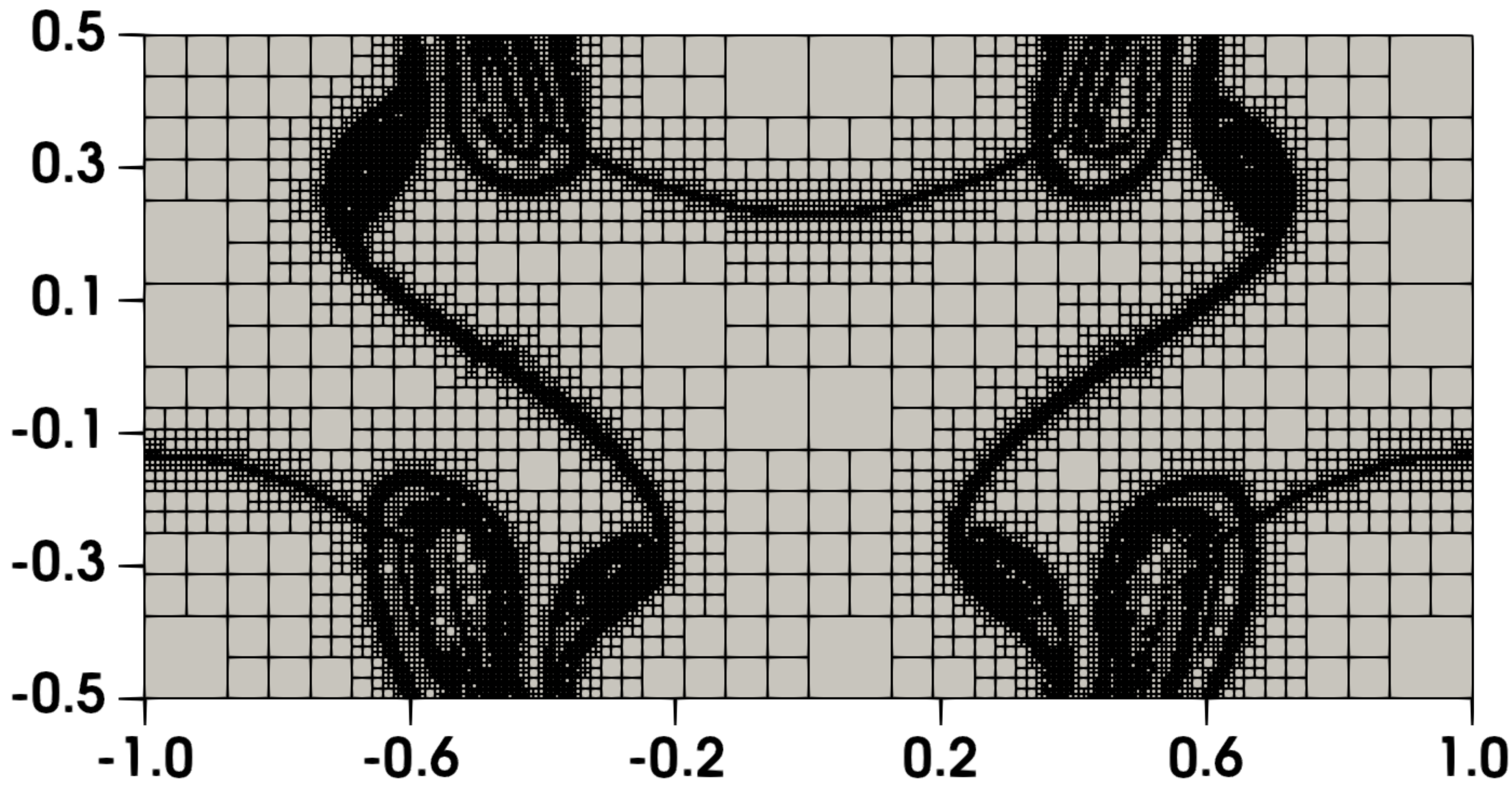}
    \caption{Adaptively refined mesh.}
    \label{fig:kh_mesh_adaptive}
    \end{subfigure}
    \vspace{0.1cm}
    \begin{subfigure}{0.49\textwidth}
    \includegraphics[width=\linewidth]{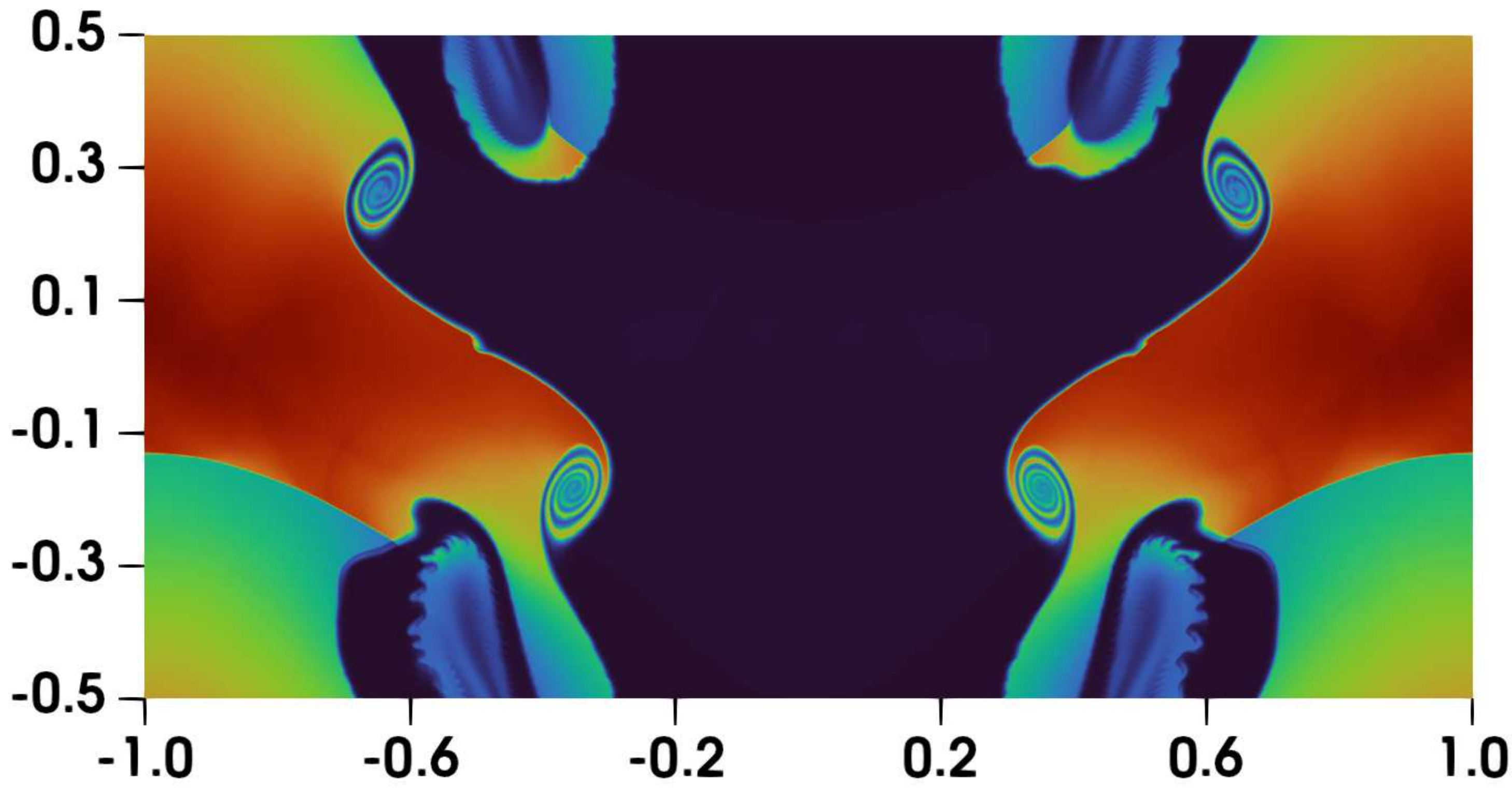}
    \caption{$\rho$ with uniform mesh $2^9 \times 2^8$.}
    \label{fig:kh_rho_uniform_8}
    \end{subfigure}
    \begin{subfigure}{0.49\textwidth}
    \includegraphics[width=\linewidth]{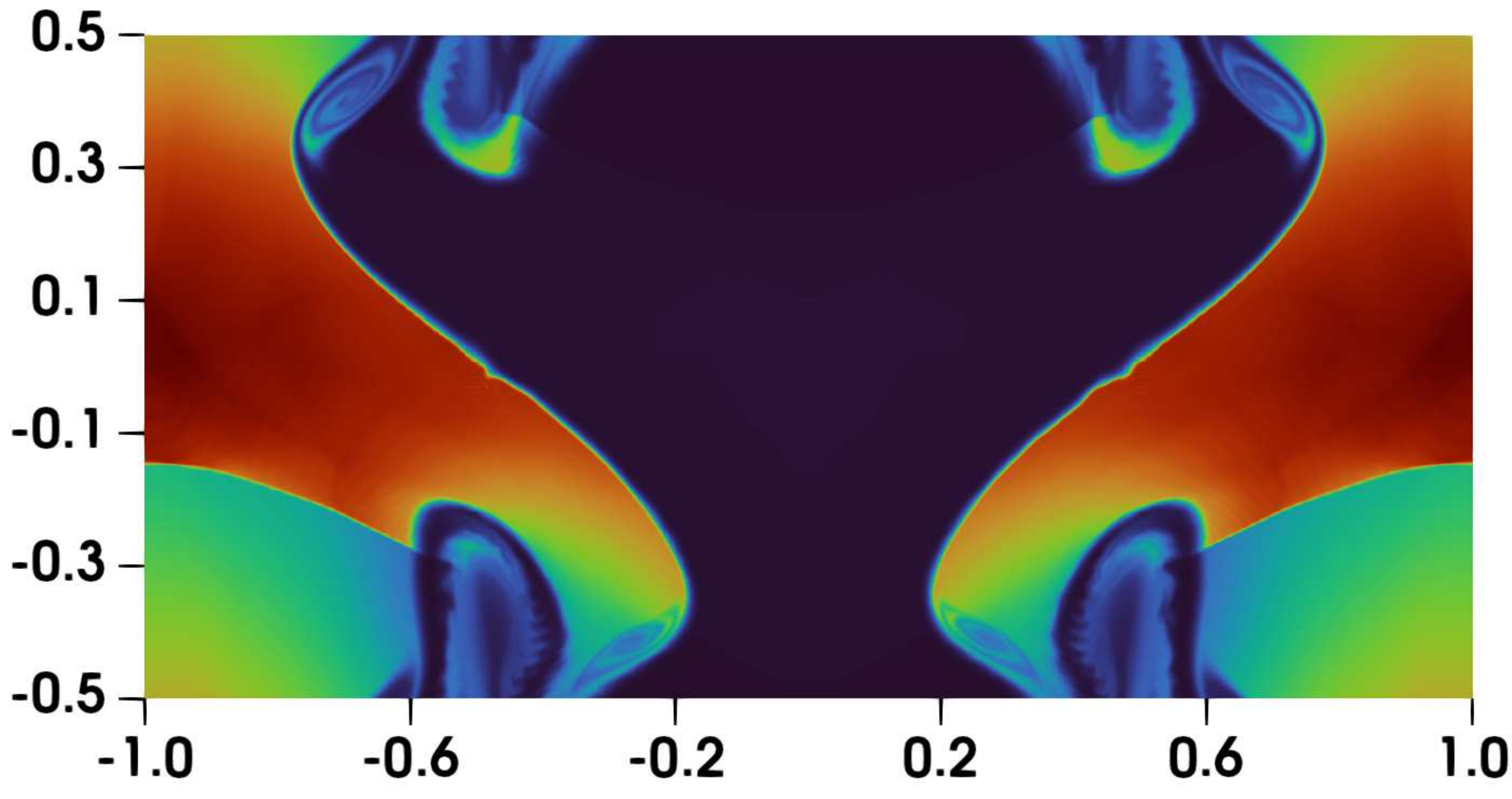}
    \caption{$\rho$ with uniform mesh $2^8 \times 2^7$.}
    \label{fig:kh_rho_uniform_7}
    \end{subfigure}
    \caption{Kelvin-Helmholtz instability test: Results at time $t=3$.}
    \label{fig:kh}
\end{figure}

We \revb{simulate} this test with $N=4$, and RC-EOS~\eqref{eq: RC_eos} \revb{until} time $t=3$ with AMR and a uniform mesh. For the AMR, we refine and coarsen the mesh in each time step, with the AMR indicator~\eqref{eq: amr_indicator} along with the three-level controller~\eqref{eq: amr_controller} having 
\begin{equation}\label{eq: case_kh_amr_params}
        (\texttt{base\_level}, \quad \texttt{med\_level}, \quad \texttt{max\_level}) = (0, 4, 8), \qquad
        (\epsilon_1, \epsilon_2) = (0.07, 0.1).
\end{equation}
The mesh for $\texttt{base\_level} = 0$ is equivalent to a mesh with $2\times 1$ elements, and $\epsilon_1, \epsilon_2$ are the thresholds in the AMR controller. The density $\rho$ profile with AMR is presented in Figure~\ref{fig:kh_rho_adaptive} along with the final mesh in Figure~\ref{fig:kh_mesh_adaptive}. We observe that the instabilities are captured by the scheme with AMR, by refining the mesh near the instabilities while coarsening it at other smoother regions. We also present the solution with two uniform mesh configurations with $2^9 \times 2^8 = 131072$ elements (equivalent to the $\texttt{max\_level}=8$), and $2^8\times 2^7=32768$ elements in Figure~\ref{fig:kh_rho_uniform_8}, and Figure~\ref{fig:kh_rho_uniform_7}, respectively. We see that the position and size of the instabilities vary with the mesh resolution, which is also reported in~\cite{radice2012thc}. The final mesh with AMR has $23969$ elements in the domain. Here as well, we note the wall-clock times for the simulations:
\begin{align*}
    &\text{with adaptive mesh: } 4709 \text{ seconds,}\\ &\text{with uniform mesh } 2^9 \times 2^8\text{: } 18003 \text{ seconds,}\\
    &\text{with uniform mesh }2^8 \times 2^7\text{: } 2390 \text{ seconds}.
\end{align*}
It is clear that the wall-clock time for the simulation with AMR is significantly less compared to that with the uniform mesh with $2^9 \times 2^8$ elements, which is equivalent to the mesh with $\texttt{max\_level} = 8$ in AMR.

\section{Conclusion}\label{sec: conclusion}
In~\cite{basak2025bound,basak2025constraints}, the Lax-Wendroff flux reconstruction method was developed for uniform Cartesian meshes to solve the RHD equations. In this work, the system of RHD equations is transformed to a reference element in such a way that it can handle meshes with curved elements. A general challenge in solving the RHD equations comes in capturing the sharp features in the domain, which are moving with time, and need a mesh with high resolution to be captured properly. A uniform mesh with such high resolution makes the simulations expensive. We use adaptive mesh refinement (AMR) to capture these structures with a finer mesh while making the mesh coarser in the smooth regions. This helps in reducing the computational cost. The mesh is allowed to be refined and coarsened in every time step with an AMR indicator~\cite{lohner1987adaptive} implemented with a three-level controller~\cite{schlottkelakemper2025trixi}. To compute the flux at the non-conformal faces, we incorporate the idea of the mortar element method~\cite{kopriva1996conservative, kopriva2002computation}. In~\cite{basak2025constraints}, the idea of approximate Lax-Wendroff procedure is used for evaluating the time average fluxes, which makes the scheme Jacobian-free. However, it requires flux evaluation of some quantities which are not guaranteed to be admissible, and this issue is handled by an additional scaling in~\cite{basak2025constraints}. Here, we use automatic differentiation (AD)~\cite{babbar2025automatic} to compute the temporal derivatives in contravariant time average fluxes, which eliminates the need for flux evaluation of non-admissible quantities, along with making the scheme Jacobian-free.  In this work, the computational boundaries are handled taking the nature of the flow into consideration, that is, the fluxes at the boundaries are computed with the information from inside as well as outside the domain, depending on the direction of the characteristics. For controlling the Gibbs oscillations in the solution, the high-order scheme is blended with a scheme of low-order at the subcell level. The first-order finite volume scheme with Rusanov flux is chosen as the low-order scheme, which preserves the admissibility of the solution~\cite{wu2017design,basak2025bound,basak2025constraints} of the RHD equations for the cases of all the equations of state used here. Then the blending limiter helps in ensuring the admissibility of the final high resolution solution. Finally, we have simulated a wide class of test cases to verify the effectiveness of the numerical scheme. We have also presented the results using a mesh with uniform resolution for all the test cases to compare them with the AMR results. It is observed that AMR results capture all the sharp structures in the domain, and are similar to the results obtained on a fine uniform mesh.  However, there is a significant difference in the wall-clock time for the simulations which shows the benefits of using AMR for such problems.

\section*{Acknowledgments}
Sujoy Basak would like to acknowledge the support of the Prime Minister's Research Fellowship, with PMRF ID 1403239. Praveen Chandrashekar would like to acknowledge the support of the Department of Atomic Energy, Government of India, under Project Identification No. 1303/9/2025-R\&D-II-DAE/TIFR-17312. Harish Kumar acknowledges support from the Vajra award (VJR/2018/00129).

\section*{Data availability}
The code and data used in this work will be made publicly available at~\cite{RHDTrixiLW}, after the publication of the paper.

\section*{Additional data}
\begin{samepage}
Animations illustrating the temporal evolutions of the simulations are available at
\begin{center}
   \href{https://www.youtube.com/playlist?list=PLrZ1LUocyVaT5ZlB2J_b9DfbAtQVf94zV}{https://www.youtube.com/playlist?list=PLrZ1LUocyVaT5ZlB2J\_b9DfbAtQVf94zV} 
\end{center}
\end{samepage}

\appendix
\section{Transformation of the conservation law into the reference element}\label{sec: appendix}

Considering the reference co-ordinates $\xi, \eta$ as a function of physical co-ordinates $x,y$, for the flux functions $\mb{f}, \mb{g}$ in~\eqref{eq: RHD_equation_2d} by the chain rule we have
\begin{align}\label{eq: chain_rule}
    \frac{\p \mb{f}}{\p x} = \frac{\p \mb{f}}{\p \xi} \frac{\p \xi}{\p x} + \frac{\p \mb{f}}{\p \eta} \frac{\p \eta}{\p x}, \qquad \frac{\p \mb{g}}{\p y} = \frac{\p \mb{g}}{\p \xi} \frac{\p \xi}{\p y} + \frac{\p \mb{g}}{\p \eta} \frac{\p \eta}{\p y}.
\end{align}
The determinant of the Jacobian~\eqref{eq: jacobian} of the transformation map~\eqref{eq: transform map} is given by
\begin{equation}
    |\mb{J}| = \frac{\p x}{\p \xi}\frac{\p y}{\p \eta} - \frac{\p x}{\p \eta} \frac{\p y}{\p \xi}.
\end{equation}

Multiplying both sides of the equations in~\eqref{eq: chain_rule} by $|\mb{J}|$ we get
\begin{align*}
   |\mb{J}|\frac{\p \mb{f}}{\p x} = \frac{\p \mb{f}}{\p \xi} \frac{\p y}{\p \eta} - \frac{\p \mb{f}}{\p \eta} \frac{\p y}{\p \xi}
    \implies \frac{\p \mb{f}}{\p x} = \frac{1}{|\mb{J}|} \left(\frac{\p \mb{f}}{\p \xi} \frac{\p y}{\p \eta} - \frac{\p \mb{f}}{\p \eta} \frac{\p y}{\p \xi} \right),
\end{align*}
and
\begin{align*}
   |\mb{J}|\frac{\p \mb{g}}{\p y} = -\frac{\p \mb{g}}{\p \xi} \frac{\p x}{\p \eta} + \frac{\p \mb{g}}{\p \eta} \frac{\p x}{\p \xi}
    \implies \frac{\p \mb{g}}{\p y} = \frac{1}{|\mb{J}|} \left(-\frac{\p \mb{g}}{\p \xi} \frac{\p x}{\p \eta} + \frac{\p \mb{g}}{\p \eta} \frac{\p x}{\p \xi} \right).
\end{align*}
Using the metric identities~\cite{kopriva2006metric}
\begin{equation}\label{eq: metric_identity}
    \frac{\p y}{\p \xi \p \eta} = \frac{\p y}{\p \eta \p \xi}, \qquad \frac{\p x}{\p \xi \p \eta} = \frac{\p x}{\p \eta \p \xi},
\end{equation}
we get
\[
    \frac{\p \mb{f}}{\p x} = \frac{1}{|\mb{J}|}\left[ \frac{\p }{\p \xi}\left(\frac{\p y}{\p \eta} \mb{f} \right) - \frac{\p }{\p \eta}\left(\frac{\p y}{\p \xi} \mb{f} \right)\right], \qquad \frac{\p \mb{g}}{\p y} = \frac{1}{|\mb{J}|}\left[ -\frac{\p }{\p \xi}\left(\frac{\p x}{\p \eta} \mb{g} \right) + \frac{\p }{\p \eta}\left(\frac{\p x}{\p \xi} \mb{g} \right)\right],
\]
and hence, the conservation law~\eqref{eq: RHD_equation_2d} becomes
\begin{align}
    &\frac{\p \mb{u}}{\p t} + \frac{1}{|\mb{J}|} \left[ 
    \frac{\p}{\p \xi} \left(\frac{\p y}{\p \eta} \mb{f} - \frac{\p x}{\p \eta} \mb{g} \right) 
    + \frac{\p}{\p \eta} \left( - \frac{\p y}{\p \xi} \mb{f} + \frac{\p x}{\p \xi} \mb{g}\right) \right]
    =\mb{0}\nonumber
\end{align}
This can be written as
\begin{align}
    \frac{\p \contra{\mb{u}}}{\p t} + \frac{\p\contra{\mb{f}}}{\p \xi} + \frac{\p\contra{\mb{g}}}{\p \eta} = \mb{0},\label{eq: transformed_cons_law_app}
\end{align}
where
\begin{equation}
    \contra{\mb{u}} = |\mb{J}| \mb{u}, \quad 
    \contra{\mb{f}} = \frac{\p y}{\p \eta} \mb{f} - \frac{\p x}{\p \eta} \mb{g}, \quad 
    \contra{\mb{g}} = - \frac{\p y}{\p \xi} \mb{f} + \frac{\p x}{\p \xi} \mb{g}.
\end{equation}

Although the above relations are enough to explain all the formulations, we show a differential geometry approach, which makes the analysis simpler when dealing with higher dimensions in a curvilinear mesh environment~\cite{kopriva2006metric, kopriva2009implementing}. We take the two-dimensional formulation of the covariant basis vectors~\cite{kopriva2006metric, kopriva2009implementing} as
\begin{equation}\label{eq: covariant}
    \mb{a}_1 = \left(\frac{\p x}{\p \xi}, \frac{\p y}{\p \xi}, 0 \right)^\top, \quad \mb{a}_2 = \left(\frac{\p x}{\p \eta}, \frac{\p y}{\p \eta}, 0 \right)^\top, \quad 
    \mb{a}_3 = (0, 0, 1)^\top,
\end{equation}
and contravariant basis vectors~\cite{kopriva2006metric, kopriva2009implementing} as
\begin{equation}\label{eq: contravariant}
    \mb{a}^1 = \left( \frac{\p \xi}{\p x}, \frac{\p \xi}{\p y}, 0\right)^\top, \quad \mb{a}^2 = \left( \frac{\p \eta}{\p x}, \frac{\p \eta}{\p y}, 0\right)^\top, \quad
    \mb{a}^3 = (0, 0, 1)^\top,
\end{equation}
which is possible by taking the reference map~\eqref{eq: transform map} as
\[
    \Theta(\xi_i, \eta_j, \zeta) = (x_i, y_j, \zeta), \qquad (\xi_i, \eta_j) \in \contra{\Omega}, (x_i, y_j) \in \Omega.
\]
We have the following relation between the covariant and the contravariant basis vectors
\begin{equation}
    |\mb{J}| \mb{a}^i = \left(\frac{\p x}{\p \xi}\frac{\p y}{\p \eta} - \frac{\p x}{\p \eta} \frac{\p y}{\p \xi}\right) \left(\frac{\p \xi}{\p x}, \frac{\p \xi}{\p y}, 0 \right) = \mb{a}_j \times \mb{a}_k,
\end{equation}
where $i,j,k \in \{1, 2, 3\}$ with $(i, j, k)$ cyclic. In the form of covariant basis vectors, the $|\mb{J}|$ can be expressed as
\begin{equation}
    |\mb{J}| = \frac{\p x}{\p \xi}\frac{\p y}{\p \eta} - \frac{\p x}{\p \eta} \frac{\p y}{\p \xi} = \mb{a}_1 \cdot (\mb{a}_2 \times \mb{a}_3).
\end{equation}
Thus, in the transformed conservation law~\eqref{eq: transformed_cons_law_app}, the flux functions are given in terms of the contravariant basis vectors~\eqref{eq: contravariant} as
\begin{align*}
    \contra{\mb{f}} &= \frac{\p y}{\p \eta} \mb{f} - \frac{\p x}{\p \eta} \mb{g} = |\mb{J}|\left(\frac{\p \xi}{\p x} \mb{f} + \frac{\p \xi}{\p y} \mb{g} \right) = |\mb{J}| \mb{a}^1 \cdot (\mb{f}, \mb{g}, \mb{0})^\top,\\
    \contra{\mb{g}} &= - \frac{\p y}{\p \xi} \mb{f} + \frac{\p x}{\p \xi} \mb{g} = |\mb{J}|\left(\frac{\p \eta}{\p x} \mb{f} + \frac{\p \eta}{\p y} \mb{g} \right) = |\mb{J}| \mb{a}^2 \cdot (\mb{f}, \mb{g}, \mb{0})^\top.
\end{align*}
Note that, the flux functions $\contra{\mb{f}}, \contra{\mb{g}}$ are called contravariant fluxes. Finally, using the contravariant basis vectors~\eqref{eq: contravariant} we can rewrite the metric identities in~\eqref{eq: metric_identity} as
\begin{equation}\label{eq: metric_identity_2}
    \frac{\p (|\mb{J}| \mb{a}^1)}{\p \xi} + \frac{\p (|\mb{J}| \mb{a}^2)}{\p \eta} = 
    \mb{0},
\end{equation}
which is similar to the metric identities in~\cite{kopriva2006metric, kopriva2009implementing}.

\begin{remark}
    As proved for the LWFR scheme in~\cite{babbar2025lax} and for the standard FR schemes in~\cite{kopriva2006metric}, satisfying the metric identity~\eqref{eq: metric_identity} or~\eqref{eq: metric_identity_2} is necessary and sufficient for a scheme to preserve the free-stream. To satisfy it discretely, some additional restrictions need to be imposed on the degree of the reference map $\Theta$~\eqref{eq: transform map}~\cite{kopriva2006metric}. For the two-dimensional case, it is sufficient to take the degree of the reference map to be no greater than the degree of the FR basis. In this work, we choose the degree of the reference map to be the same as the degree of the solution polynomial, satisfying the condition. The conditions for the higher-dimensional case can be found in~\cite{kopriva2006metric}.
\end{remark}

\bibliographystyle{abbrv}
\bibliography{reference}

@article{begelman1984theory,
  title={Theory of extragalactic radio sources},
  author={Begelman, Mitchell C and Blandford, Roger D and Rees, Martin J},
  journal={Reviews of Modern Physics},
  volume={56},
  number={2},
  pages={255},
  year={1984},
  publisher={APS}
}

@book{bottcher2012relativistic,
  title={Relativistic jets from active galactic nuclei},
  author={B{\"o}ttcher, Markus and Harris, Daniel E and Krawczynski, Henric},
  year={2012},
  publisher={John Wiley and Sons}
}

@article{mirabel1999sources,
  title={Sources of relativistic jets in the galaxy},
  author={Mirabel, I Felix and Rodriguez, Luis F},
  journal={Annual Review of Astronomy and Astrophysics},
  volume={37},
  number={1},
  pages={409--443},
  year={1999},
  publisher={Annual Reviews 4139 El Camino Way, PO Box 10139, Palo Alto, CA 94303-0139, USA}
}

@article{zensus1997parsec,
  title={Parsec-scale jets in extragalactic radio sources},
  author={Zensus, J Anton},
  journal={Annual Review of Astronomy and Astrophysics},
  volume={35},
  number={1},
  pages={607--636},
  year={1997},
  publisher={Annual Reviews 4139 El Camino Way, PO Box 10139, Palo Alto, CA 94303-0139, USA}
}

@book{godlewski1991hyperbolic,
  title={Hyperbolic systems of conservation laws},
  author={Godlewski, Edwige and Raviart, Pierre-Arnaud},
  year={1991},
  publisher={Ellipses}
}

@article{wilson1972numerical,
  title={Numerical study of fluid flow in a Kerr space},
  author={Wilson, James R},
  journal={The Astrophysical Journal},
  volume={173},
  pages={431},
  year={1972}
}

@article{marti1991numerical,
  title={Numerical relativistic hydrodynamics: Local characteristic approach},
  author={Mart{\'\i}, Jos{\'e} Ma and Ib{\'a}nez, Jos{\'e} Ma and Miralles, Juan A},
  journal={Physical Review D},
  volume={43},
  number={12},
  pages={3794},
  year={1991},
  publisher={APS}
}

@article{marti1994analytical,
  title={The analytical solution of the Riemann problem in relativistic hydrodynamics},
  author={Mart{\'\i}, Jos{\'e} Ma and M{\"u}ller, Ewald},
  journal={Journal of Fluid Mechanics},
  volume={258},
  pages={317--333},
  year={1994},
  publisher={Cambridge University Press}
}

@article{dai1997iterative,
  title={An iterative Riemann solver for relativistic hydrodynamics},
  author={Dai, Wenlong and Woodward, Paul R},
  journal={SIAM Journal on Scientific Computing},
  volume={18},
  number={4},
  pages={982--995},
  year={1997},
  publisher={SIAM}
}

@article{ibanez1999riemann,
  title={Riemann solvers in relativistic astrophysics},
  author={Ib{\'a}{\~n}ez, J M{\textordfeminine} and Mart{\i}́, J M{\textordfeminine}},
  journal={Journal of computational and applied mathematics},
  volume={109},
  number={1-2},
  pages={173--211},
  year={1999},
  publisher={Elsevier}
}

@article{marti1996extension,
  title={Extension of the piecewise parabolic method to one-dimensional relativistic hydrodynamics},
  author={Mart{\i}, Jos{\'e} Ma and M{\"u}ller, Ewald},
  journal={Journal of Computational Physics},
  volume={123},
  number={1},
  pages={1--14},
  year={1996},
  publisher={Elsevier}
}

@article{aloy1999genesis,
  title={GENESIS: A high-resolution code for three-dimensional relativistic hydrodynamics},
  author={Aloy, Miguel A and Ib{\'a}nez, J Ma and Mart{\'\i}, J Ma and M{\"u}ller, E},
  journal={The Astrophysical Journal Supplement Series},
  volume={122},
  number={1},
  pages={151},
  year={1999},
  publisher={IOP Publishing}
}

@article{mignone2005piecewise,
  title={The piecewise parabolic method for multidimensional relativistic fluid dynamics},
  author={Mignone, A and Plewa, T and Bodo, G},
  journal={The Astrophysical Journal Supplement Series},
  volume={160},
  number={1},
  pages={199},
  year={2005},
  publisher={IOP Publishing}
}

@article{dolezal1995relativistic,
  title={Relativistic hydrodynamics and essentially non-oscillatory shock capturing schemes},
  author={Dolezal, A and Wong, SSM},
  journal={Journal of Computational Physics},
  volume={120},
  number={2},
  pages={266--277},
  year={1995},
  publisher={Elsevier}
}

@article{del2002efficient,
  title={An efficient shock-capturing central-type scheme for multidimensional relativistic flows-I. Hydrodynamics},
  author={Del Zanna, Luca and Bucciantini, Niccolo},
  journal={Astronomy \& Astrophysics},
  volume={390},
  number={3},
  pages={1177--1186},
  year={2002},
  publisher={EDP Sciences}
}

@article{tchekhovskoy2007wham,
  title={WHAM: a WENO-based general relativistic numerical scheme--I. Hydrodynamics},
  author={Tchekhovskoy, Alexander and McKinney, Jonathan C and Narayan, Ramesh},
  journal={Monthly Notices of the Royal Astronomical Society},
  volume={379},
  number={2},
  pages={469--497},
  year={2007},
  publisher={Blackwell Publishing Ltd Oxford, UK}
}

@article{radice2011discontinuous,
  title={Discontinuous Galerkin methods for general-relativistic hydrodynamics: Formulation and application to spherically symmetric spacetimes},
  author={Radice, David and Rezzolla, Luciano},
  journal={Physical Review D},
  volume={84},
  number={2},
  pages={024010},
  year={2011},
  publisher={APS}
}

@article{zhang2006ram,
  title={RAM: a relativistic adaptive mesh refinement hydrodynamics code},
  author={Zhang, Weiqun and MacFadyen, Andrew I},
  journal={The Astrophysical Journal Supplement Series},
  volume={164},
  number={1},
  pages={255},
  year={2006},
  publisher={IOP Publishing}
}

@article{hughes2002three,
  title={Three-dimensional hydrodynamic simulations of relativistic extragalactic jets},
  author={Hughes, Philip A and Miller, Mark A and Duncan, G Comer},
  journal={The Astrophysical Journal},
  volume={572},
  number={2},
  pages={713},
  year={2002},
  publisher={IOP Publishing}
}

@article{wu2015high,
  title={High-order accurate physical-constraints-preserving finite difference WENO schemes for special relativistic hydrodynamics},
  author={Wu, Kailiang and Tang, Huazhong},
  journal={Journal of Computational Physics},
  volume={298},
  pages={539--564},
  year={2015},
  publisher={Elsevier}
}

@article{wu2016physical,
  title={Physical-constraint-preserving central discontinuous Galerkin methods for special relativistic hydrodynamics with a general equation of state},
  author={Wu, Kailiang and Tang, Huazhong},
  journal={The Astrophysical Journal Supplement Series},
  volume={228},
  number={1},
  pages={3},
  year={2016},
  publisher={IOP Publishing}
}

@article{wu2017design,
  title={Design of provably physical-constraint-preserving methods for general relativistic hydrodynamics},
  author={Wu, Kailiang},
  journal={Physical Review D},
  volume={95},
  number={10},
  pages={103001},
  year={2017},
  publisher={APS}
}

@article{qin2016bound,
  title={Bound-preserving discontinuous Galerkin methods for relativistic hydrodynamics},
  author={Qin, Tong and Shu, Chi-Wang and Yang, Yang},
  journal={Journal of Computational Physics},
  volume={315},
  pages={323--347},
  year={2016},
  publisher={Elsevier}
}

@article{bhoriya2020entropy,
  title={Entropy-stable schemes for relativistic hydrodynamics equations},
  author={Bhoriya, Deepak and Kumar, Harish},
  journal={Zeitschrift f{\"u}r angewandte Mathematik und Physik},
  volume={71},
  pages={1--29},
  year={2020},
  publisher={Springer}
}

@article{LW1,
author = {Lax, Peter and Wendroff, Burton},
title = {Systems of conservation laws},
journal = {Communications on Pure and Applied Mathematics},
volume = {13},
number = {2},
pages = {217-237},
year = {1960}
}

@article{lou2020flux,
  title={The flux reconstruction method with Lax--Wendroff type temporal discretization for hyperbolic conservation laws},
  author={Lou, Shuai and Yan, Chao and Ma, Li-Bin and Jiang, Zhen-Hua},
  journal={Journal of Scientific Computing},
  volume={82},
  pages={1--25},
  year={2020},
  publisher={Springer}
}

@inproceedings{huynh2007flux,
  title={A flux reconstruction approach to high-order schemes including discontinuous Galerkin methods},
  author={Huynh, Hung T},
  booktitle={18th AIAA computational fluid dynamics conference},
  pages={4079},
  year={2007}
}

@inproceedings{vincent2016towards,
  title={Towards green aviation with python at petascale},
  author={Vincent, Peter and Witherden, Freddie and Vermeire, Brian and Park, Jin Seok and Iyer, Arvind},
  booktitle={SC'16: Proceedings of the International Conference for High Performance Computing, Networking, Storage and Analysis},
  pages={1--11},
  year={2016},
  organization={IEEE}
}

@inproceedings{lopez2014verification,
  title={Verification and Validation of HiFiLES: a High-Order LES unstructured solver on multi-GPU platforms},
  author={L{\'o}pez, Manuel R and Sheshadri, Abhishek and Bull, Jonathan R and Economon, Thomas D and Romero, Joshua and Watkins, Jerry E and Williams, David M and Palacios, Francisco and Jameson, Antony and Manosalvas, David E},
  booktitle={32nd AIAA applied aerodynamics conference},
  pages={3168},
  year={2014}
}

@article{vandenhoeck2019implicit,
  title={Implicit high-order flux reconstruction solver for high-speed compressible flows},
  author={Vandenhoeck, Ray and Lani, Andrea},
  journal={Computer Physics Communications},
  volume={242},
  pages={1--24},
  year={2019},
  publisher={Elsevier}
}

@article{burger2017approximate,
  title={Approximate Lax--Wendroff discontinuous Galerkin methods for hyperbolic conservation laws},
  author={B{\"u}rger, Raimund and Kenettinkara, Sudarshan Kumar and Zor{\'\i}o, David},
  journal={Computers \& Mathematics with Applications},
  volume={74},
  number={6},
  pages={1288--1310},
  year={2017},
  publisher={Elsevier}
}

@article{BABBAR2022111423,
  title={Lax-Wendroff flux reconstruction method for hyperbolic conservation laws},
  author={Babbar, Arpit and Kenettinkara, Sudarshan Kumar and Chandrashekar, Praveen},
  journal={Journal of Computational Physics},
  volume={467},
  pages={111423},
  year={2022},
  publisher={Elsevier}
}

@article{babbar2024admissibility,
  title={Admissibility Preserving Subcell Limiter for Lax--Wendroff Flux Reconstruction},
  author={Babbar, Arpit and Kenettinkara, Sudarshan Kumar and Chandrashekar, Praveen},
  journal={Journal of Scientific Computing},
  volume={99},
  number={2},
  pages={31},
  year={2024},
  publisher={Springer}
}

@article{anile2005relativistic,
  title={Relativistic fluids and magneto-fluids},
  author={Anile, Angelo Marcello},
  journal={Relativistic Fluids and Magneto-fluids},
  year={2005}
}

@book{landau1987see,
  title={Fluid Mechanics. Chapter XV - Relativistic Fluid Dynamics, 2nd edn},
  author={Landau, LD and Lifshitz, EM},
  year={1987},
  publisher={Pergamon, New York}
}

@article{synge1965relativity,
  title={Relativity: the special theory},
  author={Synge, John Lighton},
  year={1965}
}

@article{ryu2006equation,
  title={Equation of state in numerical relativistic hydrodynamics},
  author={Ryu, Dongsu and Chattopadhyay, Indranil and Choi, Eunwoo},
  journal={The Astrophysical Journal Supplement Series},
  volume={166},
  number={1},
  pages={410},
  year={2006},
  publisher={IOP Publishing}
}

@article{zorio2017approximate,
  title={An approximate Lax--Wendroff-type procedure for high order accurate schemes for hyperbolic conservation laws},
  author={Zor{\'\i}o, David and Baeza, Antonio and Mulet, Pep},
  journal={Journal of Scientific Computing},
  volume={71},
  pages={246--273},
  year={2017},
  publisher={Springer}
}

@article{rusanov1962calculation,
  title={The calculation of the interaction of non-stationary shock waves and obstacles},
  author={Rusanov, Vladimir Vasil’evich},
  journal={USSR Computational Mathematics and Mathematical Physics},
  volume={1},
  number={2},
  pages={304--320},
  year={1962},
  publisher={Elsevier}
}

@article{hennemann2021provably,
  title={A provably entropy stable subcell shock capturing approach for high order split form DG for the compressible Euler equations},
  author={Hennemann, Sebastian and Rueda-Ram{\'\i}rez, Andr{\'e}s M and Hindenlang, Florian J and Gassner, Gregor J},
  journal={Journal of Computational Physics},
  volume={426},
  pages={109935},
  year={2021},
  publisher={Elsevier}
}

@article{zhang2010maximum,
  title={On maximum-principle-satisfying high order schemes for scalar conservation laws},
  author={Zhang, Xiangxiong and Shu, Chi-Wang},
  journal={Journal of Computational Physics},
  volume={229},
  number={9},
  pages={3091--3120},
  year={2010},
  publisher={Elsevier}
}

@article{nunez2016xtroem,
  title={XTROEM-FV: a new code for computational astrophysics based on very high order finite-volume methods--II. Relativistic hydro-and magnetohydrodynamics},
  author={N{\'u}{\~n}ez-de La Rosa, Jonatan and Munz, Claus-Dieter},
  journal={Monthly Notices of the Royal Astronomical Society},
  volume={460},
  number={1},
  pages={535--559},
  year={2016},
  publisher={Oxford University Press}
}

@article{he2012adaptive,
  title={An adaptive moving mesh method for two-dimensional relativistic hydrodynamics},
  author={He, Peng and Tang, Huazhong},
  journal={Communications in Computational Physics},
  volume={11},
  number={1},
  pages={114--146},
  year={2012},
  publisher={Cambridge University Press}
}

@article{biswas2022entropy,
  title={Entropy stable discontinuous Galerkin schemes for the special relativistic hydrodynamics equations},
  author={Biswas, Biswarup and Kumar, Harish and Bhoriya, Deepak},
  journal={Computers \& Mathematics with Applications},
  volume={112},
  pages={55--75},
  year={2022},
  publisher={Elsevier}
}

@article{beckwith2011second,
  title={A second-order Godunov method for multi-dimensional relativistic magnetohydrodynamics},
  author={Beckwith, Kris and Stone, James M},
  journal={The Astrophysical Journal Supplement Series},
  volume={193},
  number={1},
  pages={6},
  year={2011},
  publisher={IOP Publishing}
}

@article{radice2012thc,
  title={THC: a new high-order finite-difference high-resolution shock-capturing code for special-relativistic hydrodynamics},
  author={Radice, David and Rezzolla, Luciano},
  journal={Astronomy \& Astrophysics},
  volume={547},
  pages={A26},
  year={2012},
  publisher={EDP Sciences}
}

@article{zanotti2015high,
  title={A high order special relativistic hydrodynamic and magnetohydrodynamic code with space--time adaptive mesh refinement},
  author={Zanotti, Olindo and Dumbser, Michael},
  journal={Computer Physics Communications},
  volume={188},
  pages={110--127},
  year={2015},
  publisher={Elsevier}
}

@article{duan2019high,
  title={High-order accurate entropy stable finite difference schemes for one-and two-dimensional special relativistic hydrodynamics},
  author={Duan, Junming and Tang, Huazhong},
  journal={arXiv preprint arXiv:1905.06092},
  year={2019}
}

@article{ling2019physical,
  title={Physical-constraints-preserving Lagrangian finite volume schemes for one-and two-dimensional special relativistic hydrodynamics},
  author={Ling, Dan and Duan, Junming and Tang, Huazhong},
  journal={Journal of Computational Physics},
  volume={396},
  pages={507--543},
  year={2019},
  publisher={Elsevier}
}

@article{BALSARA1994,
title = {Riemann Solver for Relativistic Hydrodynamics},
journal = {Journal of Computational Physics},
volume = {114},
number = {2},
pages = {284-297},
year = {1994},
issn = {0021-9991},
author = {Dinshaw S. Balsara}
}

@article{sokolov2001simple,
  title={Simple and efficient Godunov scheme for computational relativistic gas dynamics},
  author={Sokolov, IV and Zhang, H-M and Sakai, JI},
  journal={Journal of Computational Physics},
  volume={172},
  number={1},
  pages={209--234},
  year={2001},
  publisher={Elsevier}
}

@article{xu2024high,
  title={High-Order Accurate Entropy Stable Schemes for Relativistic Hydrodynamics with General Synge-Type Equation of State},
  author={Xu, Linfeng and Ding, Shengrong and Wu, Kailiang},
  journal={Journal of Scientific Computing},
  volume={98},
  number={2},
  pages={43},
  year={2024},
  publisher={Springer}
}

@article{mathews1971hydromagnetic,
  title={The hydromagnetic free expansion of a relativistic gas},
  author={Mathews, William G},
  journal={Astrophysical Journal, vol. 165, p. 147},
  volume={165},
  pages={147},
  year={1971}
}

@article{gottlieb2009high,
  title={High order strong stability preserving time discretizations},
  author={Gottlieb, Sigal and Ketcheson, David I and Shu, Chi-Wang},
  journal={Journal of Scientific Computing},
  volume={38},
  number={3},
  pages={251--289},
  year={2009},
  publisher={Springer}
}

@article{chen2022physical,
  title={A physical-constraint-preserving finite volume WENO method for special relativistic hydrodynamics on unstructured meshes},
  author={Chen, Yaping and Wu, Kailiang},
  journal={Journal of Computational Physics},
  volume={466},
  pages={111398},
  year={2022},
  publisher={Elsevier}
}

@article{wu2014finite,
  title={Finite volume local evolution Galerkin method for two-dimensional relativistic hydrodynamics},
  author={Wu, Kailiang and Tang, Huazhong},
  journal={Journal of Computational Physics},
  volume={256},
  pages={277--307},
  year={2014},
  publisher={Elsevier}
}

@article{wu2014third,
  title={A third-order accurate direct Eulerian GRP scheme for one-dimensional relativistic hydrodynamics},
  author={Wu, Kailiang and Yang, Zhicheng and Tang, Huazhong},
  journal={East Asian Journal on Applied Mathematics},
  volume={4},
  number={2},
  pages={95--131},
  year={2014},
  publisher={Cambridge University Press}
}

@article{wu2021minimum,
  title={Minimum principle on specific entropy and high-order accurate invariant-region-preserving numerical methods for relativistic hydrodynamics},
  author={Wu, Kailiang},
  journal={SIAM Journal on Scientific Computing},
  volume={43},
  number={6},
  pages={B1164--B1197},
  year={2021},
  publisher={SIAM}
}

@article{cai2024provably,
  title={Provably convergent Newton--Raphson methods for recovering primitive variables with applications to physical-constraint-preserving Hermite WENO schemes for relativistic hydrodynamics},
  author={Cai, Chaoyi and Qiu, Jianxian and Wu, Kailiang},
  journal={Journal of Computational Physics},
  volume={498},
  pages={112669},
  year={2024},
  publisher={Elsevier}
}

@article{brandt1977multi,
  title={Multi-level adaptive solutions to boundary-value problems},
  author={Brandt, Achi},
  journal={Mathematics of computation},
  volume={31},
  number={138},
  pages={333--390},
  year={1977}
}

@article{berger1984adaptive,
  title={Adaptive mesh refinement for hyperbolic partial differential equations},
  author={Berger, Marsha J and Oliger, Joseph},
  journal={Journal of computational Physics},
  volume={53},
  number={3},
  pages={484--512},
  year={1984},
  publisher={Elsevier}
}

@article{berger1989local,
  title={Local adaptive mesh refinement for shock hydrodynamics},
  author={Berger, Marsha J and Colella, Phillip},
  journal={Journal of computational Physics},
  volume={82},
  number={1},
  pages={64--84},
  year={1989},
  publisher={Elsevier}
}

@article{o2005adaptive,
  title={Adaptive Mesh Refinement--Theory and Applications},
  author={O’shea, BW and Bryan, G and Bordner, J and Norman, Michael L and Abel, T and Harkness, R and Kritsuk, A and Plewa, T and Linde, T and Weirs, VG},
  journal={Lectures Notes of Computer Science Engineering},
  volume={41},
  pages={341--350},
  year={2005}
}

@article{plewa2001amra,
  title={AMRA: An Adaptive Mesh Refinement hydrodynamic code for astrophysics},
  author={Plewa, Tomek and Mueller, Ewald},
  journal={Computer Physics Communications},
  volume={138},
  number={2},
  pages={101--127},
  year={2001},
  publisher={Elsevier}
}

@article{balsara2001divergence,
  title={Divergence-free adaptive mesh refinement for magnetohydrodynamics},
  author={Balsara, Dinshaw S},
  journal={Journal of Computational Physics},
  volume={174},
  number={2},
  pages={614--648},
  year={2001},
  publisher={Elsevier}
}

@book{donmez2002general,
  title={General relativistic hydrodynamics with adaptive-mesh refinement (AMR) and modeling of accretion disks},
  author={Donmez, Orhan},
  year={2002},
  publisher={Drexel University}
}

@article{anninos2005cosmos++,
  title={Cosmos++: relativistic magnetohydrodynamics on unstructured grids with local adaptive refinement},
  author={Anninos, Peter and Fragile, P Chris and Salmonson, Jay D},
  journal={The Astrophysical Journal},
  volume={635},
  number={1},
  pages={723},
  year={2005},
  publisher={IOP Publishing}
}

@article{wang2008relativistic,
  title={Relativistic hydrodynamic flows using spatial and temporal adaptive structured mesh refinement},
  author={Wang, Peng and Abel, Tom and Zhang, Weiqun},
  journal={The Astrophysical Journal Supplement Series},
  volume={176},
  number={2},
  pages={467},
  year={2008},
  publisher={IOP Publishing}
}

@article{babbar2025lax,
  title={Lax-Wendroff flux reconstruction on adaptive curvilinear meshes with error based time stepping for hyperbolic conservation laws},
  author={Babbar, Arpit and Chandrashekar, Praveen},
  journal={Journal of Computational Physics},
  volume={522},
  pages={113622},
  year={2025},
  publisher={Elsevier}
}

@article{qiu2005discontinuous,
  title={The discontinuous Galerkin method with Lax--Wendroff type time discretizations},
  author={Qiu, Jianxian and Dumbser, Michael and Shu, Chi-Wang},
  journal={Computer methods in applied mechanics and engineering},
  volume={194},
  number={42-44},
  pages={4528--4543},
  year={2005},
  publisher={Elsevier}
}

@article{basak2025bound,
  title={Bound preserving Lax-Wendroff flux reconstruction method for special relativistic hydrodynamics},
  author={Basak, Sujoy and Babbar, Arpit and Kumar, Harish and Chandrashekar, Praveen},
  journal={Journal of Computational Physics},
  volume={527},
  pages={113815},
  year={2025},
  publisher={Elsevier}
}

@article{basak2025constraints,
  title={Constraints Preserving Lax-Wendroff Flux Reconstruction for Relativistic Hydrodynamics with General Equations of State},
  author={Basak, Sujoy and Babbar, Arpit and Kumar, Harish and Chandrashekar, Praveen},
  journal={Journal of Scientific Computing},
  volume={105},
  number={3},
  pages={70},
  year={2025},
  publisher={Springer}
}

@article{babbar2025automatic,
title = {Automatic differentiation for performing the Cauchy-Kovalevskaya procedure in Lax-Wendroff type discretizations},
journal = {Journal of Computational Physics},
volume = {563},
pages = {115101},
year = {2026},
issn = {0021-9991},
doi = {https://doi.org/10.1016/j.jcp.2026.115101},
url = {https://www.sciencedirect.com/science/article/pii/S0021999126004547},
author = {Arpit Babbar and Valentin Churavy and Michael Schlottke-Lakemper and Hendrik Ranocha}
}

@article{qiu2003finite,
  title={Finite difference WENO schemes with Lax--Wendroff-type time discretizations},
  author={Qiu, Jianxian and Shu, Chi-Wang},
  journal={SIAM Journal on Scientific Computing},
  volume={24},
  number={6},
  pages={2185--2198},
  year={2003},
  publisher={SIAM}
}

@book{griewank2008evaluating,
  title={Evaluating derivatives: principles and techniques of algorithmic differentiation},
  author={Griewank, Andreas and Walther, Andrea},
  year={2008},
  publisher={SIAM}
}

@book{kopriva2009implementing,
  title={Implementing spectral methods for partial differential equations: Algorithms for scientists and engineers},
  author={Kopriva, David A},
  year={2009},
  publisher={Springer Science \& Business Media}
}

@article{kopriva2006metric,
  title={Metric identities and the discontinuous spectral element method on curvilinear meshes},
  author={Kopriva, David A},
  journal={Journal of Scientific Computing},
  volume={26},
  number={3},
  pages={301--327},
  year={2006},
  publisher={Springer}
}

@article{moses2020instead,
  title={Instead of rewriting foreign code for machine learning, automatically synthesize fast gradients},
  author={Moses, William and Churavy, Valentin},
  journal={Advances in neural information processing systems},
  volume={33},
  pages={12472--12485},
  year={2020}
}

@article{revels2016forward,
  title={Forward-mode automatic differentiation in Julia},
  author={Revels, Jarrett and Lubin, Miles and Papamarkou, Theodore},
  journal={arXiv preprint arXiv:1607.07892},
  year={2016}
}

@phdthesis{tan2023higher,
  title={Higher-Order Automatic Differentiation and Its Applications},
  author={Tan, Songchen},
  year={2023},
  school={Massachusetts Institute of Technology}
}

@misc{tan2022taylordiff,
  title={Taylordiff. jl: fast higher-order automatic differentiation in Julia},
  howpublished={\url{https://github.com/JuliaDiff/TaylorDiff.jl}},
  author={Tan, Songchen},
  year={2022}
}

@article{di1857note,
  title={Note sur une nouvelle formule de calcul diff{\'e}rentiel},
  author={Di Bruno, F Fa{\`a}},
  journal={Quarterly J. Pure Appl. Math},
  volume={1},
  number={359-360},
  pages={12},
  year={1857}
}

@techreport{carlson2011inflow,
  title={Inflow/outflow boundary conditions with application to FUN3D},
  author={Carlson, Jan-Rene{\'e}},
  year={2011}
}

@inproceedings{kim2004integrated,
  title={Integrated simulations for multi-component analysis of gas turbines: RANS boundary conditions},
  author={Kim, Sangho and Alonso, Juan and Schluter, Jorg and Wu, Xiaohua and Pitsch, Heinz},
  booktitle={40th AIAA/ASME/SAE/ASEE Joint Propulsion Conference and Exhibit},
  pages={3415},
  year={2004}
}

@techreport{strikwerda1976initial,
  title={Initial Boundary Value Problems for Incompletely Parabolic Systems.},
  author={Strikwerda, John Charles},
  year={1976}
}

@article{kreiss1970initial,
  title={Initial boundary value problems for hyperbolic systems},
  author={Kreiss, Heinz-Otto},
  journal={Communications on Pure and Applied Mathematics},
  volume={23},
  number={3},
  pages={277--298},
  year={1970},
  publisher={Wiley Online Library}
}

@article{higdon1986initial,
  title={Initial-boundary value problems for linear hyperbolic system},
  author={Higdon, Robert L},
  journal={SIAM review},
  volume={28},
  number={2},
  pages={177--217},
  year={1986},
  publisher={SIAM}
}

@article{lohner1987adaptive,
  title={An adaptive finite element scheme for transient problems in CFD},
  author={L{\"o}hner, Rainald},
  journal={Computer methods in applied mechanics and engineering},
  volume={61},
  number={3},
  pages={323--338},
  year={1987},
  publisher={Elsevier}
}

@article{kopriva1996conservative,
  title={A conservative staggered-grid Chebyshev multidomain method for compressible flows},
  author={Kopriva, David A and Kolias, John H},
  journal={Journal of computational physics},
  volume={125},
  number={1},
  pages={244--261},
  year={1996},
  publisher={Elsevier}
}

@article{kopriva2002computation,
  title={Computation of electromagnetic scattering with a non-conforming discontinuous spectral element method},
  author={Kopriva, David A and Woodruff, Stephen L and Hussaini, M Yousuff},
  journal={International journal for numerical methods in engineering},
  volume={53},
  number={1},
  pages={105--122},
  year={2002},
  publisher={Wiley Online Library}
}

@inproceedings{pao1981numerical,
  title={A numerical study of two-dimensional shock vortex interaction},
  author={Pao, S and Salas, M},
  booktitle={14th Fluid and Plasma Dynamics Conference},
  pages={1205},
  year={1981}
}

@article{balsara2016subluminal,
  title={A subluminal relativistic magnetohydrodynamics scheme with ADER-WENO predictor and multidimensional Riemann solver-based corrector},
  author={Balsara, Dinshaw S and Kim, Jinho},
  journal={Journal of Computational Physics},
  volume={312},
  pages={357--384},
  year={2016},
  publisher={Elsevier}
}

@article{duncan1994simulations,
  title={Simulations of relativistic extragalactic jets},
  author={Duncan, G Comer and Hughes, Philip A},
  journal={arXiv preprint astro-ph/9406041},
  year={1994}
}

@article{komissarov1998large,
  title={The large-scale structure of FR-II radio sources},
  author={Komissarov, SS and Falle, SAEG},
  journal={Monthly Notices of the Royal Astronomical Society},
  volume={297},
  number={4},
  pages={1087--1108},
  year={1998},
  publisher={Blackwell Science Ltd Oxford, UK}
}

@article{geuzaine2009gmsh,
  title={Gmsh: A 3-D finite element mesh generator with built-in pre-and post-processing facilities},
  author={Geuzaine, Christophe and Remacle, Jean-Fran{\c{c}}ois},
  journal={International journal for numerical methods in engineering},
  volume={79},
  number={11},
  pages={1309--1331},
  year={2009},
  publisher={Wiley Online Library}
}

@article{ranocha2021adaptive,
  title={Adaptive numerical simulations with Trixi. jl: A case study of Julia for scientific computing},
  author={Ranocha, Hendrik and Schlottke-Lakemper, Michael and Winters, Andrew R and Faulhaber, Erik and Chan, Jesse and Gassner, Gregor J},
  journal={arXiv preprint arXiv:2108.06476},
  year={2021}
}

@misc{schlottkelakemper2025trixi,
  title={{T}rixi.jl: {A}daptive high-order numerical simulations
         of hyperbolic {PDE}s in {J}ulia},
  author={Schlottke-Lakemper, Michael and Gassner, Gregor J and
          Ranocha, Hendrik and Winters, Andrew R and Chan, Jesse
          and Rueda-Ramírez, Andrés},
  year={2025},
  howpublished={\url{https://github.com/trixi-framework/Trixi.jl}},
  doi={10.5281/zenodo.3996439}
}

@article{schlottke2021purely,
  title={A purely hyperbolic discontinuous Galerkin approach for self-gravitating gas dynamics},
  author={Schlottke-Lakemper, Michael and Winters, Andrew R and Ranocha, Hendrik and Gassner, Gregor J},
  journal={Journal of Computational Physics},
  volume={442},
  pages={110467},
  year={2021},
  publisher={Elsevier}
}

@article{chan2025artificial,
  title={An artificial viscosity approach to high order entropy stable discontinuous Galerkin methods},
  author={Chan, Jesse},
  journal={arXiv preprint arXiv:2501.16529},
  year={2025}
}

@misc{trixilw,
  title={{T}rixi{LW}.jl: {L}ax-{W}endroff {F}lux {R}econstruction on curvilinear grids},
  author={Babbar, Arpit and Chandrashekar, Praveen},
  year={2024},
  howpublished={\url{https://github.com/Arpit-Babbar/TrixiLW.jl}},
  doi={https://doi.org/10.5281/zenodo.13822912}
}

@article{godunov1959finite,
  title={Finite difference method for numerical computation of discontinuous solutions of the equations of fluid dynamics},
  author={Godunov, Sergei K and Bohachevsky, Ihor},
  journal={Matemati{\v{c}}eskij sbornik},
  volume={47},
  number={3},
  pages={271--306},
  year={1959}
}

@misc{RHDTrixiLW,
  title={{RHDT}rixi{LW}.jl: {R}elativistic {H}ydrodynamics with adaptive mesh refinement using {LWFR} scheme},
  author={Basak, Sujoy and Babbar, Arpit and Kumar, Harish and Chandrashekar, Praveen},
  year={2026},
  howpublished={\url{https://github.com/sujoy-basak/RHDTrixiLW.jl}},
  doi={https://doi.org/10.5281/zenodo.19398810}
}

@article{dumbser2014posteriori,
  title={A posteriori subcell limiting of the discontinuous Galerkin finite element method for hyperbolic conservation laws},
  author={Dumbser, Michael and Zanotti, Olindo and Loub{\`e}re, Rapha{\"e}l and Diot, Steven},
  journal={Journal of Computational Physics},
  volume={278},
  pages={47--75},
  year={2014},
  publisher={Elsevier}
}

@article{dumbser2016simple,
  title={A simple robust and accurate a posteriori sub-cell finite volume limiter for the discontinuous Galerkin method on unstructured meshes},
  author={Dumbser, Michael and Loub{\`e}re, Rapha{\"e}l},
  journal={Journal of Computational Physics},
  volume={319},
  pages={163--199},
  year={2016},
  publisher={Elsevier}
}

@article{zanotti2015solving,
  title={Solving the relativistic magnetohydrodynamics equations with ADER discontinuous Galerkin methods, a posteriori subcell limiting and adaptive mesh refinement},
  author={Zanotti, Olindo and Fambri, Francesco and Dumbser, Michael},
  journal={Monthly Notices of the Royal Astronomical Society},
  volume={452},
  number={3},
  pages={3010--3029},
  year={2015},
  publisher={Oxford University Press}
}
\end{document}